\documentclass[aps,showpacs,superscriptaddress,groupedaddress]{revtex4-1}

\usepackage{amsmath}
\usepackage{amssymb}
\allowdisplaybreaks[2]

\usepackage{graphicx}
\usepackage{subfigure}
\usepackage{longtable} %tables

\usepackage{multirow}

\usepackage{dcolumn}
\usepackage{bm}
\usepackage{array}
\usepackage{xcolor} %text color

\begin{document}

\title{Buffon-Laplace Needle Problem as a geometric probabilistic approach to filtration process}

\author{Yan-Jie Min}
\thanks{These authors contributed equally to this work.}

\author{De-Quan Zhu}
\thanks{These authors contributed equally to this work.}

\author{Jin-Hua Zhao}
\email{zhaojh@m.scnu.edu.cn}

\affiliation{CoCoLab, School of Data Science and Engineering,
South China Normal University, Shanwei 516622, China}

\date{\today}
%\date{January 29, 2024}

\begin{abstract}
%background
Buffon-Laplace Needle Problem considers a needle of a length $l$ randomly dropped on a large plane distributed with vertically parallel lines with distances $a$ and $b$ ($a \geqslant b$), respectively.
As a classical problem in stochastic probability, it serves as a mathematical basis of various physical literature, such as the efficiency of a filter and the emergence of clogging in filtration process.
Yet its potential application is limited by previous focus on its original form of the `short' needle case of $l < b$ and its analytical difficulty in a general sense.
%our contribution
Here, rather than a `short' needle embedded in two-dimensional space, we analytically solve problem versions with needles and spherocylinders of arbitrary length and radius embedded in two- and three-dimensional spaces dropped on a grid with any rectangular shape.
We further confirm our analytical theory with Monte Carlo simulation.
Our framework here helps to provide a geometric analytical perspective to filtration process, and also extend the analytical power of the needle problem into unexplored parameter regions for physical problems involving stochastic processes.

\end{abstract}

%Keywords: Buffon-Laplace Needle Problem; Buffon's Needle Problem; intersection probability; filtration process

\maketitle

%%%table of contents
%\tableofcontents
%\newpage

\section{Introduction}

%from filtration to BLNP
Filtration process of particles through a filter or a porous medium is a fundamental subject in soft and granular matter research \cite{Tien.Ramarao-2007-2e, Jaeger.Nagel.Behringer-RMP-1996}. Theoretical methods are developed to characterize its static and dynamical properties, such as the efficiency of a filter and the onset of clogging \cite{Redner.Datta-PRL-2000, To.Lai.Pak-PRL-2001, Roussel.etal-PRL-2007, Gabrielli.etal-PRL-2013, Gerber.etal-PRL-2018}. Some methods explore the geometry of particles and filters to construct probabilistic models.
In the paper \cite{Roussel.etal-PRL-2007}, a filtration process of particles is considered in two parts: (1) individual particles in a flow arrive at and make contact with a mesh structure, and (2) the interaction between the particles and the mesh after contacts leads to a formation of an arch of particles and results in clogging.
Following this logic, the paper \cite{Dell.Franklin-JSTAT-2009} considers particles as rod-like objects and adopts a generalized version of the Buffon-Laplace Needle Problem (BLNP) as a geometric probabilistic approach to the above first part.
In this paper, we follow the paper \cite{Dell.Franklin-JSTAT-2009} and focus on the analytical result of BLNP for filtration process.

%BLNP
BLNP \cite{Laplace-1812, Laplace-1820} is one of the oldest problems in geometric probability dated back more than two centuries.
BLNP considers randomly dropping a thin needle (i.e. a line segment) with a length $l$ on a two-dimensional ($2$D)  plane distributed with two sets of equidistant parallel lines with distances $a$ and $b$, respectively, which are vertical to each other, and studies the statistical property of intersections between the needle and the lines. For convenience, we set $a \geqslant b$ in this paper.
%BNP
BLNP is a natural extension of Buffon's Needle Problem (BNP) \cite{Buffon-1777}, which considers a similar setting on a plane yet with a single set of parallel lines with a distance $d$.

%other contexts for needle problems: random sequential adsorption, random-line graph
For the two needle problems, except the filtration process of particles, there are other contexts in which the scenario of needle dropping is also involved, such as random sequential adsorption which considers dropping line segments and other geometrical shapes in a non-overlapping way onto a substrate to study its jamming property \cite{Renyi-1958, Vigil.Ziff-JChemPhys-1990, Svrakic.Henkel-JPhysI-1991, Evans-RMP-1993}, and random-line graphs which are formed by intersections by line segments dropped in a square on a plane \cite{Boettcher-JCompNetw-2020}.

%summary of existing results
%(1) intersection probability
The fundamental quantity for the two needle problems is the intersection probability of a needle intersecting with at least one line.
For the intersection probability $P(l, d)$ with a needle length $l$ and a line distance $d$ in BNP, \cite{Buffon-1777} and \cite{Uspensky-1937}(page 251-252) derive its form for the `short' needle case with $l < d$. \cite{Uspensky-1937}(page 258, problem 7) also shows $P(l, d)$ for a long needle case with $l > d$. 
For the intersection probability $P(l, a, b)$ with a needle length $l$ and line distances $a$ and $b$ in BLNP, \cite{Laplace-1812, Laplace-1820} and \cite{Uspensky-1937}(page 255-257) show a prediction of the case with $l < b$.
In \cite{Dell.Franklin-JSTAT-2009}, the authors further lay down an analytical form of $P(l, a, b)$ in a specific case of $a = b$ and $l > b$ (as Eq. (6) in its main text), yet with an error in a coefficient which we will show later.
%(2) distribution of intersections
Another quantity pertinent to the two needle problems is the statistics of numbers of intersections and visited grids for a long needle.
For BNP,  \cite{Diaconis-JApplProb-1976} shows the distribution of intersection numbers in the case of $l > d$.
For BLNP, \cite{Arkhipov.Mendo-Mathematika-2023} studies the average and the maximal number of visited grids (tiles) for a needle on any rectangular grid.
 
%variants of BLNP and BNP
As classical models of randomness, variants based on BNP and BLNP are frequently discussed in both mathematical and physical literature.
A simple variation of the problems is to generalize needles and grid lines into other geometrical objects.
%BNP
For BNP, a needle can be replaced with
a planar curve \cite{Barbier-JMathPA-1860, Ramaley-AMM-1969},
polygons \cite{Johannesen-MathJ-2009},
a pivot needle (two needles sharing a tip) \cite{Baesel-ElemMath-2015}, and so on.
Sets of parallel lines can also be generalized to
parallel stripes with a finite width \cite{Chung-MathGeo-1981},
random Cantor sets \cite{Zhang-RevMatIberoam-2020, Vardakis.Volberg-JMAA-2024},
a line with randomly distributed dots \cite{Godreche-JSTAT-2022}, and so on.
%BLNP
For BLNP, a needle can be extended to an ellipse \cite{Duma.Stoka-JAppProb-1993},
a spherocylinder (a rod-like shape consistent of a cylinder and two hemispherical caps on both its ends) in three-dimensional ($3$D) space \cite{Dell.Franklin-JSTAT-2009, Pournin.etal-GranMat-2005}, and so on.
Besides, a plane in BLNP can be considered as being covered with congruent triangles \cite{Uspensky-1937} (page 258, problem 8).

%our focus
In this paper, we revisit BLNP and focus on its analytical solution of intersection probability, the most basic quantity of the problem.
%(1)previous results
For the original BLNP and its extensions, the initial `short' needle case with $l < b$ is heavily discussed in previous literature, while analytical result for the `long' needle case of $l > b$ is quite limited.
In the paper \cite{Dell.Franklin-JSTAT-2009}, after a needle in $2$D space is generalized to a spherocylinder in $3$D space, BLNP is adopted as a model for granular particles arriving at a sieve-like grid. Yet it sets $a = b$ for a square mesh, and lays down the intersection probability in the initial form of integrals and estimates the probability with Monte Carlo simulation \cite{Landau.Binder-2015-4e}.
%(2)our emphasis
Here we emphasize that, neither forcing $l < b$ nor $a = b$ is necessary for a theoretical model, and models with any permissible $l$ and $(a, b)$ represents a realistic setting when the problem serves as an approximation model in a physical context.

%our contribution
Our main contribution in this paper is that, with a single geometric probabilistic framework,  we solve BLNP with dropped objects as needles with no volume and spherocylinders with finite volume in both $2$D and $3$D spaces, and derive analytical predictions of intersection probability for any shape parameter (a length $l$ for a needle and a length $l$ and a radius $\sigma / 2$ for a spherocylinder) and any grid shape parameters $(a, b)$.
%further implication
Our results theoretically elucidate the contact probability between spherocylinders in flow in $3$D space with a mesh for the model presented in \cite{Dell.Franklin-JSTAT-2009}, and makes possible studying new complex phenomena originated from long particles and irregular grid topology. Our theory helps to construct a full geometric probabilistic picture of the filtration process further combined with the interaction between particles and mesh.
Besides, on the mathematical side, our framework naturally retrieves the result for the original BLNP with `short' needles. In all, our framework helps to release the potential of BLNP as a mathematical basis in various physical problems.

%layout
Here is the layout of the paper.
In Section \ref{sec:model}, we present BLNP and its four versions, and also lay down their representation and realization in simulation.
In Sections \ref{sec:theory-2D} to \ref{sec:theory-3DSC}, we explain the analytical framework to calculate the intersection probability when needles and spherocylinders embedded in $2$D and $3$D spaces, respectively, are dropped on a $2$D grid.
In Section \ref{sec:theory-comparison}, we compare our theoretical results with previous ones.
In Section \ref{sec:result}, we show our theoretical results, verified with Monte Carlo simulation.
In Section \ref{sec:conclusion}, we conclude the paper with some discussion.
In Appendices A--J, we leave details of calculation for the analytical framework.

%%%%%%%%%%%%%%%%%%%%%%%%%%%%%%%%
\section{Model}
\label{sec:model}

\begin{figure}[htbp]
\begin{center}
 \includegraphics[width = 0.85 \linewidth]{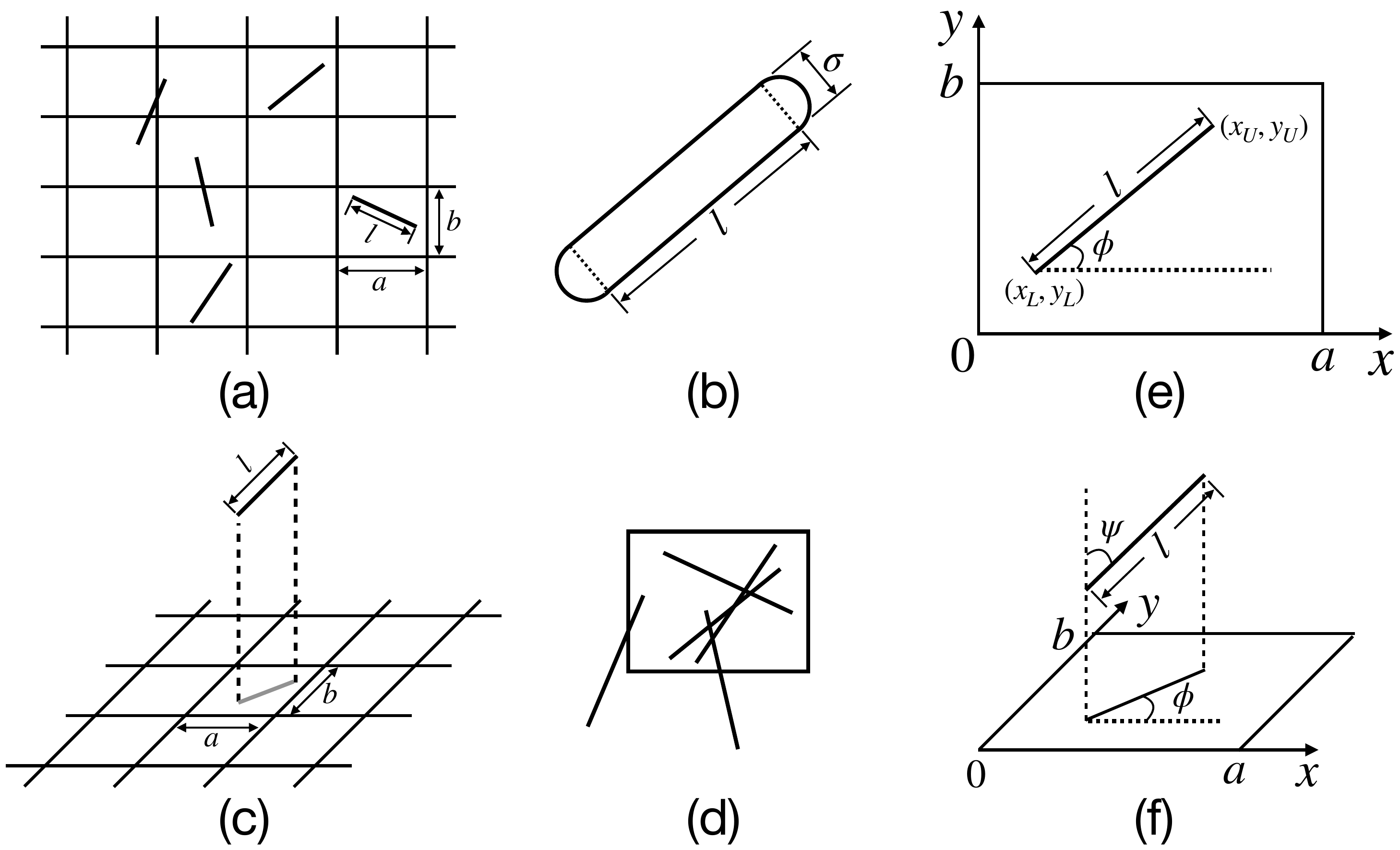}
\end{center}
\caption{
 \label{fig:model}
BLNP with needles and spherocylinders in $2$D and $3$D spaces.
(a) In the original BLNP with needles, a thin needle of a length $l$ is repeatedly dropped on a plane distributed with two vertical sets of parallel lines with distances of $a$ and $b$, respectively. The five line segments in the grid represent needle experiments, in two of which a needle intersects with lines.
(b) A $2$D needle is extended to a $2$D spherocylinder, which consists of a rectangle with a width $l$ and two semicircles with a diameter $\sigma$ on both its ends.
(c) A $2$D needle is extended to a $3$D one, in which a needle of a length $l$ in a $3$D space is dropped onto a $2$D plane. A $3$D needle intersects with lines only when its projection on the $2$D plane intersects with them.
(d) The five needle configurations in (a) are moved in a periodical way into a single grid cell, while their statuses between a needle and lines are retained respectively.
(e) In the BLNP with $2$D needles, the relative position of a needle within a grid cell with a width $a$ and a height $b$ is shown. The coordinates $(x_{\rm U}, y_{\rm U})$ and $(x_{\rm L}, y_{\rm L})$ represent respectively the upper and the lower tips of the needle, and the orientation $\phi \in [0, \pi)$ is the angle between the needle and the $x$-axis.
(f) In the BLNP with $3$D needles, a $3$D needle is firstly projected onto a $2$D grid.
The orientation $\psi \in [0, \pi/2]$ is the angle between the needle and the direction perpendicular to the grid.
The relative position of a projected needle in a grid cell shares a similar description in the case of $2$D needles shown in (e).}
\end{figure}
%

%generalized of 2D needles
In the original BLNP, the grid is $2$D with two vertical sets of parallel lines, and a needle is dropped on the plane. See Fig.\ref{fig:model}(a) for an example.
In this paper, when an object is embedded in $2$D or $3$D spaces, we call the type of the object as $2$D or $3$D one, respectively, for convenience.
We keep the grid as $2$D and generalize the dropped object from a $2$D needle into two directions.
%needles -> spherocylinders
We first generalize a thin needle with no volume into a rod-like shape with a finite volume. This generalization makes possible consideration of real-world objects in a mathematical model. A spherocylinder \cite{Dell.Franklin-JSTAT-2009, Pournin.etal-GranMat-2005} is intrinsically a $3$D object, and consists of a cylinder with a length $l$ and two semispheres with a radius $\sigma / 2$ on both its ends.
When such a spherocylinder lies on a plane, its $2$D projection on the plane consists of a rectangle with a width $l$ and two semicircles with a diameter $\sigma$ on both its ends. We simply call this shape a $2$D spherocylinder.
See Fig.\ref{fig:model}(b) for an example.
%2D -> 3D
We then generalize the embedded space of dropped objects from $2$D to $3$D, during which an extra degree of freedom is introduced for the objects. See Fig.\ref{fig:model}(c) for an example of a $3$D needle.
Beware that, an intersection between a $3$D object and a $2$D grid is defined between the projection of the $3$D object onto the grid plane and the grid lines.
Due to its simple topology, the projection of a $3$D spherocylinder on a $2$D plane is a $2$D one, yet with a shorter length except when it is parallel to the plane.

%4 versions of the problem
Summing up the above generalizations of BLNP, we consider here its four versions:
(1) BLNP with $2$D needles, which is simply the original BLNP;
(2) BLNP with $2$D spherocylinders;
(3) BLNP with $3$D needles; and
(4) BLNP with $3$D spherocylinders.

%simulation: a general Monte-Carlo simulation
%(1) grids -> a grid
We further lay down how to perform simulation for these four versions.
To count the intersection events, it is only relevant to distinguish whether or not a dropped object intersects with sets of lines, which only involves the relative position between them. Thus the picture of an object dropped on a  grid with an infinitely large number of cells can be simplified into one in which an object is dropped on a single grid cell.
See Fig.\ref{fig:model}(d) for an example when all needles are moved into a single grid cell.
We choose parallel lines with distances $b$ and $a$ to establish $x$- and $y$-axes, respectively, and then consider the relative position of a dropped object in a grid cell.
%(2)a general process
As a needle and a spherocylinder both have simple geometrical boundaries, we can use only a limited number of parameters, for example coordinates of their two tips, to describe their overall positions.
A general Monte Carlo simulation\cite{Landau.Binder-2015-4e} of an object-dropping experiment can be carried out as follows: with given sizes of objects and grid cells $(l, a, b)$ or $(l, \sigma, a, b)$, we first generate a large number $N_{\rm all}$ of object configurations with a pseudo-random number generator, say Rand(), which generates a sequence of real numbers $\in (0, 1)$ obeying a uniform distribution;
we then count the number $N_{\rm coll}$ of configurations in which an object overlaps with at least one boundary of a grid cell;  finally, we calculate the empirical intersection probability as $N_{\rm coll} / N_{\rm all}$.

%simulation: representation of configurations
In the above simulation process, only the description of an object configuration and the condition for an intersection vary in the four problem versions.
%configurations: 2D objects
For a $2$D needle and spherocylinder, their descriptions of relative positions share a similar form and consist of three ingredients: the coordinate $x_{\rm U}$ on $x$-axis and the coordinate $y_{\rm U}$ on $y$-axis for the upper tip  of a needle or the center of the upper semicircle of a spherocylinder, and the orientation $\phi$ between a needle or a spherocylinder and the $x$-axis. The `upper' tip or semicircle here corresponds to the one of tips or semicircles with a larger coordinate on the $y$-axis.
Fig.\ref{fig:model}(e) shows an illustration for the $2$D needle case.
The ranges of $(x_{\rm U}, y_{\rm U}, \phi )$ are $[0, a)$, $[0, b)$, and $[0, \pi)$, respectively. 
In a single experiment, a random configuration of a $2$D needle or spherocylinder can be generated as
\begin{equation}
\label{eq:configuration-2D-uppertip}
\left[ \begin{array}{c}
x_{\rm U} \\
y_{\rm U} \\
\phi \end{array} \right]
= \left[ \begin{array}{c}
a \times {\rm Rand()} \\
b \times {\rm Rand()} \\
\pi \times {\rm Rand()} \end{array} \right].
\end{equation}
Correspondingly, the coordinates of the lower tip of the needle or the center of the lower semicircle of the spherocylinder are
\begin{equation}
\label{eq:configuration-2D-lowertip}
\left[ \begin{array}{c}
x_{\rm L} \\
y_{\rm L}\end{array} \right]
= \left[ \begin{array}{c}
x_{\rm U} - l \cos \phi \\
y_{\rm U} - l \sin \phi \end{array} \right].
\end{equation}
%

%configurations: 3D objects
For the description of a $3$D needle or spherocylinder, except from the parameters $(x_{\rm U}, y_{\rm U}, \phi )$, a fourth parameter $\psi \in [0, \pi / 2]$ is introduced to account the orientation between an object and the axis vertical to the grid plane.
Fig.\ref{fig:model} (f) shows an illustration for the $3$D needle case.
To generate a random configuration of a $3$D needle or spherocylinder, we have
\begin{equation}
\label{eq:configuration-3D-uppertip}
\left[ \begin{array}{c}
x_{\rm U} \\
y_{\rm U} \\
\psi \\
\phi \end{array} \right]
= \left[ \begin{array}{c}
a \times {\rm Rand()} \\
b \times {\rm Rand()} \\
\frac {\pi}{2} \times {\rm Rand()} \\
\pi \times {\rm Rand()} \end{array} \right].
\end{equation}
Correspondingly, the coordinates of the lower tip of the $3$D needle or the center of the lower semicircle of the $3$D spherocylinder of their projected shapes are respectively
\begin{equation}
\label{eq:configuration-3D-lowertip}
\left[ \begin{array}{c}
x_{\rm L} \\
y_{\rm L}\end{array} \right]
= \left[ \begin{array}{c}
x_{\rm U} - l \sin \psi \cos \phi \\
y_{\rm U} - l \sin \psi \sin \phi \end{array} \right].
\end{equation}
%

%intersection constraints: needles
Then we consider the condition for an intersection between a dropped object and a grid cell.
An intersection between a $2$D needle and the grid happens simply when at least one tip is out of the grid cell. We can state as
\begin{equation}
\label{eq:intersection-needle}
\begin{aligned}
& \left( x_{\rm L} < 0 \right) \vee \left( a < x_{\rm L} \right) \\
\vee & \left(y_{\rm L} < 0 \right) \vee \left( b < y_{\rm L} \right),
\end{aligned}
\end{equation}
in which $\vee$ is a logical OR operator.
Beware that, in the numerical simulation with continuous real numbers, we ignore the cases when the two needle tips fall exactly on any boundary of a grid cell (say, $x = 0$ and $a$, and $y = 0$ and $b$), since their probability is infinitesimal.

%intersection constraints: spherocylinders
To account an intersection between a $2$D spherocylinder, also the $2$D projection of a $3$D spherocylinder, and a grid cell, the contour of a spherocylinder should be taken into consideration. It is easy to find that, an intersection happens when at least one of the furthermost tips in four directions of the contour of a spherocylinder is out of a grid cell. We have the following condition as
\begin{equation}
\label{eq:intersection-SC}
\begin{aligned}
& \left( x_{\rm U} - \frac {\sigma}{2} < 0 \right) \vee \left( a < x_{\rm U} + \frac {\sigma}{2} \right) \\
\vee & \left( y_{\rm U} - \frac {\sigma}{2} < 0 \right) \vee \left( b < y_{\rm U} + \frac {\sigma}{2} \right) \\
\vee & \left( x_{\rm L} - \frac {\sigma}{2} < 0 \right) \vee \left( a < x_{\rm L} + \frac {\sigma}{2} \right) \\
\vee & \left( y_{\rm L} - \frac {\sigma}{2} < 0 \right) \vee \left( b < y_{\rm L} + \frac {\sigma}{2} \right).
\end{aligned}
\end{equation}
%

%BLNP --> BNP
When the line distance $a \to \infty$, the constraint from parallel lines along the $y$-axis vanishes, and two sets of parallel lines effectively reduce to only one set. In this case, the BLNP with $2$D needles, $2$D spherocylinders, $3$D needles, and $3$D spherocylinders reduce to its BNP versions, respectively.
For the simulation of four BNP versions, we can simply remove $x_{\rm U}$ and $x_{\rm L}$ from the description of object configuration in Eqs.(\ref{eq:configuration-2D-uppertip})-(\ref{eq:configuration-3D-lowertip}), along with those terms with $x_{\rm U}$ and $x_{\rm L}$ in Eqs.(\ref{eq:intersection-needle}) and (\ref{eq:intersection-SC}).

%to theoretical analysis
Based on the above representations of the relative position of an object dropped on a grid, we further develop a consistent analytical framework to calculate the intersection probability with any size of objects and grid cells. In the following four sections, starting from the original BLNP version with $2$D needles, we theoretically solve problems with increasing complexity until the one with $3$D spherocylinders.
Furthermore, by setting $1 / a \to 0$, these theoretical equations naturally reduce to those for the four BNP versions, respectively.

%%%%%%%%%%%%%%%%%%%%%%%%%%%%%%%%
\section{Theory of BLNP with $2$D Needles}
\label{sec:theory-2D}

To compute the intersection probability for the four BLNP versions, our basic logic is to navigate in the parameter space of configurations of a dropped object in a grid cell and calculate the fraction of the volume of parameter space, in which an intersection happens, among the total volume of the parameter space.

Specifically, for the BLNP with $2$D needles, we calculate two volumes: $S_{\rm all}$ as one of the total parameter space of $(x_{\rm U}, y_{\rm U}, \phi)$ for a needle configuration, and $S_{\rm coll}$ as one of the proper parameter space of $(x_{\rm U}, y_{\rm U}, \phi)$ which results in an intersection between the needle and the grid cell. The intersection probability $P(l, a, b)$ can be formulated as
\begin{equation}
P(l, a, b) = \frac {S_{\rm coll}}{S_{\rm all}}.
\end{equation}
It is easy to find that
\begin{equation}
S_{\rm all} = ab\pi.
\end{equation}
The challenging part here is to derive an explicit form for $S_{\rm coll}$. We reformulate $S_{\rm coll} $ as
\begin{equation}
\label{eq:Scoll-2D}
S_{\rm coll} = S_{\rm coll}(x) + S_{\rm coll}(y) - S_{\rm coll}(x, y),
\end{equation}
in which $S_{\rm coll}(x)$, $S_{\rm coll}(y)$, and $S_{\rm coll}(x, y)$ denote the volume of the parameter space of $(x_{\rm U}, y_{\rm U}, \phi)$ when an intersection between a $2$D needle and a grid cell happens on the boundaries of the grid cell along the $y$-axis (say, $x = 0$ and $a$), the $x$-axis (say, $y = 0$ and $b$), and the $x$- and $y$-axes at the same time, respectively.
In the following equations, we denote $(x, y)$ for $(x_{\rm U}, y_{\rm U})$ for convenience.
Based on Eqs.(\ref{eq:configuration-2D-uppertip}), (\ref{eq:configuration-2D-lowertip}), and (\ref{eq:intersection-needle}), we have
\begin{align}
%Scoll-x
\label{eq:Scollx-2D}
S_{\rm coll} (x)
& = \int _{0}^{\frac{\pi}{2}} {\rm d}\phi
\int _{0}^{a} \Theta(l \cos \phi - x) {\rm d}x \int _{0}^{b} {\rm d}y
+ \int _{\frac{\pi}{2}}^{\pi} {\rm d}\phi
\int _{0}^{a} \Theta(- l \cos \phi + x - a) {\rm d}x \int _{0}^{b} {\rm d}y  \nonumber \\
& = b \int _{0}^{\pi} {\rm d}\phi
\int _{0}^{a} \Theta(|l \cos \phi| - x) {\rm d}x, \\
%Scoll-y
\label{eq:Scolly-2D}
S_{\rm coll} (y)
& = \int _{0}^{\pi} {\rm d}\phi
\int _{0}^{a} {\rm d}x
\int _{0}^{b} \Theta(l \sin \phi - y) {\rm d}y \nonumber \\
& = a \int _{0}^{\pi} {\rm d}\phi
\int _{0}^{b} \Theta(l \sin \phi - y) {\rm d}y, \\
%Scoll-xy
\label{eq:Scollxy-2D}
S_{\rm coll} (x,y)
& = \int _{0}^{\frac{\pi}{2}} {\rm d}\phi
\int _{0}^{a} \Theta(l \cos \phi - x) {\rm d}x
\int _{0}^{b} \Theta(l \sin \phi - y) {\rm d}y
+ \int _{\frac{\pi}{2}}^{\pi} {\rm d}\phi
\int _{0}^{a} \Theta(- l \cos \phi + x - a) {\rm d}x
\int _{0}^{b} \Theta(l \sin \phi - y) {\rm d}y \nonumber \\
& = \int _{0}^{\pi} {\rm d}\phi
\int _{0}^{a} \Theta(|l \cos \phi| - x) {\rm d}x
\int _{0}^{b} \Theta(l \sin \phi - y) {\rm d}y.
\end{align}
The step function $\Theta (x)$ is defined as $\Theta (x) = 1$ if $x \geqslant 0$ and $0$ otherwise.
Combining a step function with the integrands on the parameters $(x_{\rm U}, y_{\rm U}, \phi)$ singles out proper ranges which permit an intersection between a needle and a grid cell.
On the righthand side of Eqs.(\ref{eq:Scollx-2D}) and (\ref{eq:Scollxy-2D}),
we make a change of variables as $x - a \to - \hat x$. We have
\begin{align}
\int _{0}^{a} \Theta (- l \cos \phi + x - a) {\rm d} x
& = \int _{a}^{0} \Theta (- l \cos \phi - \hat x) {\rm d} (- \hat{x}) \nonumber \\
& = \int _{0}^{a} \Theta (- l \cos \phi - \hat x) {\rm d} \hat x.
\end{align}
%

%an alternative way to calculation
An alternative way to calculate $S_{\rm coll}$ is
\begin{equation}
S_{\rm coll} = S_{\rm all} - S_{\rm non-coll},
\end{equation}
in which $S_{\rm non-coll}$ is the volume of the parameter space of $(x_{\rm U}, y_{\rm U}, \phi)$ where there is no intersection between a needle and a grid cell. Following the same logic in Eqs.(\ref{eq:Scollx-2D})-(\ref{eq:Scollxy-2D}), we have
\begin{equation}
S_{\rm non-coll}
= \int _{0}^{\pi} {\rm d}\phi
\int _{0}^{a} \Theta(x - |l \cos \phi|) {\rm d}x
\int _{0}^{b} \Theta(y - l \sin \phi) {\rm d}y.
\end{equation}
We should mention that these two frameworks to calculate $S_{\rm coll}$ are equivalent.
Yet the former is more direct in computation, which we will follow in this paper.

% define summations
To calculate Eqs.(\ref{eq:Scollx-2D})-(\ref{eq:Scollxy-2D}), we first define three integrals as
\begin{align}
\label{eq:Ala-2D}
A(l, a)
& \equiv \int _{0}^{\pi} {\rm d} \phi \int _{0}^{a} \Theta(| l \cos \phi | - x) {\rm d}x, \\
\label{eq:Blb-2D}
B(l, b)
& \equiv \int _{0}^{\pi} {\rm d} \phi \int _{0}^{b} \Theta(l \sin \phi - y) {\rm d}y, \\
\label{eq:ABlab-2D}
AB(l, a, b)
& \equiv \int _{0}^{\pi} {\rm d}\phi
\int _{0}^{a} \Theta(|l \cos \phi| - x) {\rm d}x
\int _{0}^{b} \Theta(l \sin \phi - y) {\rm d}y.
\end{align}
Correspondingly, we rewrite $S_{\rm coll}$ as
\begin{equation}
\label{eq:Scoll-2D-reform}
S_{\rm coll} = b A(l, a) + a B(l, b) - AB(l, a, b).
\end{equation}
%

%final equations
To finish the computation of Eqs.(\ref{eq:Ala-2D})-(\ref{eq:ABlab-2D}) and finally $P(l, a, b)$, we should examine the relative sizes among $l$, $a$, and $b$.
We leave details of calculation in Appendix A.
With Eqs.(\ref{eq:Plab-2D-lb}), (\ref{eq:Plab-2D-bla}), (\ref{eq:Plab-2D-al1}), and (\ref{eq:Plab-2D-al2}), we lay down our final equations as
\begin{equation}
\label{eq:Plab-2D-full}
P(l, a, b) =
\left\{ \begin{array}{ll}
% l < b
 \frac {2}{\pi} \left( \frac {l}{a} + \frac {l}{b} \right) - \frac {1}{\pi} \frac {l^2}{ab},
 & {\rm if}\ l \leqslant b; \\
 % b < l < a
 \frac {1}{\pi} \frac {b}{a}
+ \frac {2}{\pi} \frac {l}{b} \left( 1 - \sqrt{1 - \frac {b^2}{l^2} } \right)
+ \frac {2}{\pi} \left( \frac {\pi}{2} - \arcsin \frac {b}{l} \right),
& {\rm if}\ b < l \leqslant a; \\
% a < sqrt{a^2 + b^2}
 \frac {1}{\pi} \left( \frac {L^2}{ab} + \frac {l^2}{ab} \right)
- \frac {2}{\pi} \left( \frac {l}{a} \sqrt {1 - \frac {a^2}{l^2} } + \frac {l}{b} \sqrt {1 - \frac {b^2}{l^2} } \right) \\
+ \frac {2}{\pi} \left( \pi - \arcsin \frac {a}{l} - \arcsin \frac {b}{l} \right),
& {\rm if}\ a < l \leqslant L; \\
% sqrt{a^2 + b^2}  < l
1, & {\rm if}\ L< l,
\end{array}\right.
\end{equation}
in which we define $L \equiv \sqrt {a^2 + b^2}$.

%special cases: l > L and P = 1
We find that when $l > L$, we have $P(l, a, b) = 1$.
This result has a simple geometrical interpretation. The longest distance between two points in a grid cell with a width $a$ and a height $b$ is the diagonal line with a distance of $L$. Thus any needle with a length $l > L$ dropped on a grid is certain to intersect with grid lines, leading to $P(l, a, b) = 1$.

%check boundary cases
A necessary condition for the correctness of Eq.(\ref{eq:Plab-2D-full}) is the consistency of $P(l, a, b)$, and equivalently $S_{\rm coll}$, at the boundaries of the four cases of $l$. We leave details of calculation in Appendix B.

%BLNP --> BNP
When $a \to \infty$, then $l / a \to 0$. Eq.(\ref{eq:Plab-2D-full}) simplifies to the result for the BNP with $2$D needles as
\begin{equation}
\label{eq:Plb-2D-full}
P(l, \infty, b)
 = \left\{ \begin{array}{ll}
 \frac {2}{\pi} \frac {l}{b}, & {\rm if}\ l \leqslant b; \\
 \frac {2}{\pi} \frac {l}{b} \left( 1 - \sqrt {1 - \frac {b^2}{l^2} } \right)
 + \frac {2}{\pi} \left( \frac {\pi}{2} - \arcsin \frac {b}{l} \right), & {\rm if}\ b < l.
 \end{array}\right.
 \end{equation}
%

%%%%%%%%%%%%%%%%%%%%%%%%%%%%%%%%
\section{Theory of BLNP with $2$D Spherocylinders}
\label{sec:theory-2DSC}

Following the language in the BLNP with $2$D needles, we calculate the intersection probability for the BLNP with $2$D spherocylinders as
\begin{equation}
P(l, \sigma, a, b) = \frac {S_{\rm coll}}{S_{\rm all}}.
\end{equation}
With a slight abuse of notations, we define $S_{\rm all}$ as the volume of total parameter space of $(x_{\rm U}, y_{\rm U}, \phi)$, and $S_{\rm coll}$ as the volume of parameters space of $(x_{\rm U}, y_{\rm U}, \phi)$ which permit an intersection between a $2$D spherocylinder and a grid.
Rewriting $S_{\rm coll}$ as Eq.(\ref{eq:Scoll-2D}), we define $S_{\rm coll}(x)$, $S_{\rm coll}(y)$, and $S_{\rm coll}(x, y)$ here as the volume of the parameter space of $(x_{\rm U}, y_{\rm U}, \phi)$ which lead to an intersection between a $2$D spherocylinder and a grid cell on the boundaries of the grid cell along the $y$-axis, the $x$-axis, and the $x$- and $y$-axes at the same time, respectively.
Based on Eqs.(\ref{eq:configuration-2D-uppertip}), (\ref{eq:configuration-2D-lowertip}), and (\ref{eq:intersection-SC}), and following the derivation in Eqs.(\ref{eq:Scollx-2D})-(\ref{eq:Scollxy-2D}), we have
\begin{align}
%S_{coll}(x)
\label{eq:Scollx-2DSC}
S_{\rm coll}(x)
& = \int _{0}^{\frac{\pi}{2}} {\rm d} \phi \left[
\int _{0}^{\frac {\sigma}{2}} {\rm d}x
+ \int _{\frac {\sigma}{2}}^{a - \frac {\sigma}{2}} \Theta \left(l \cos \phi - x + \frac {\sigma}{2} \right) {\rm d}x + \int _{a - \frac {\sigma}{2}}^{a} {\rm d}x \right] \int_{0}^{b} {\rm d}y  \nonumber \\
& + \int _{\frac {\pi}{2}}^{\pi} {\rm d} \phi \left[
\int _{0}^{\frac {\sigma}{2}} {\rm d}x
+ \int _{\frac {\sigma}{2}}^{a - \frac {\sigma}{2}} \Theta \left(- l \cos \phi + x + \frac {\sigma}{2} - a\right) {\rm d}x
+ \int _{a - \frac {\sigma}{2}}^{a} {\rm d}x \right] \int_{0}^{b} {\rm d}y \nonumber \\
%equal sign 2
& = b \int _{0}^{\pi}  {\rm d}\phi
\int _{0}^{a - \sigma} \Theta \left(| l \cos \phi | - x \right) {\rm d}x + \sigma b \pi. \\
%S_{coll}(y)
\label{eq:Scolly-2DSC}
S_{\rm coll}(y) 
& = \int _{0}^{\pi} {\rm d} \phi
\int _{0}^{a} {\rm d} x
\left[ \int _{0}^{\frac {\sigma}{2}} {\rm d}y
+ \int _{\frac {\sigma}{2}}^{b - \frac {\sigma}{2}} \Theta \left(l \sin \phi - y + \frac {\sigma}{2} \right) {\rm d}y + \int _{b - \frac {\sigma}{2}}^{b} {\rm d}y \right] \nonumber \\
%equal sign 2
& = a \int _{0}^{\pi} {\rm d} \phi \int _{0}^{b - \sigma} \Theta(l \sin \phi - y) {\rm d}y + \sigma a \pi, \\
%S_{coll}(x,y)
\label{eq:Scollxy-2DSC}
S_{\rm coll}(x, y)
& = \int _{0}^{\frac{\pi}{2}} {\rm d} \phi
\left[ \int _{0}^{\frac {\sigma}{2}} {\rm d}x
+ \int _{\frac {\sigma}{2}}^{a - \frac {\sigma}{2}} \Theta \left(l \cos \phi - x + \frac {\sigma}{2} \right) {\rm d}x + \int _{a - \frac {\sigma}{2}}^{a} {\rm d}x \right]  \nonumber \\
& \times \left[ \int _{0}^{\frac {\sigma}{2}} {\rm d}y
+ \int _{\frac {\sigma}{2}}^{b - \frac {\sigma}{2}} \Theta \left(l \sin \phi - y + \frac {\sigma}{2} \right) {\rm d}y + \int _{b - \frac {\sigma}{2}}^{b} {\rm d}y \right] \nonumber \\
& + \int _{\frac {\pi}{2}}^{\pi} {\rm d} \phi
\left[ \int _{0}^{\frac {\sigma}{2}} {\rm d}x
+ \int _{\frac {\sigma}{2}}^{a - \frac {\sigma}{2}} \Theta \left( - l \cos \phi + x + \frac {\sigma}{2} - a \right) {\rm d}x + \int _{a - \frac {\sigma}{2}}^{a} {\rm d}x \right] \nonumber \\
& \times \left[ \int _{0}^{\frac {\sigma}{2}} {\rm d}y
+ \int _{\frac {\sigma}{2}}^{b - \frac {\sigma}{2}} \Theta \left( l \sin \phi - y + \frac {\sigma}{2} \right) {\rm d}y + \int _{b - \frac {\sigma}{2}}^{b} {\rm d}y \right] \nonumber \\
%equal sign 2
& = \int _{0}^{\pi} {\rm d}\phi
\int _{0}^{a - \sigma} \Theta (| l \cos \phi | - x) {\rm d}x
\int _{0}^{b - \sigma} \Theta (l \sin \phi - y) {\rm d}y \nonumber \\
& + \sigma \int _{0}^{\pi} {\rm d}\phi
\int _{0}^{a - \sigma} \Theta (| l \cos \phi | - x) {\rm d}x
+ \sigma \int _{0}^{\pi} {\rm d}\phi
\int _{0}^{b - \sigma} \Theta (l \sin \phi - y) {\rm d}y
+ \sigma ^2 \pi.
\end{align}
In the above equations, we expand the brackets and further make some rearrangement in integrals with a change of variables.
With $x - \frac{\sigma}{2} \to \hat x$, we have
\begin{equation}
\int _{\frac {\sigma}{2}}^{a - \frac {\sigma}{2}}
\Theta \left( l \cos \phi - x + \frac {\sigma}{2} \right) {\rm d}x
= \int _{0}^{a - \sigma} \Theta \left( l \cos \phi - \hat x \right) {\rm d} \hat x.
\end{equation}
With $x + \frac{\sigma}{2} - a \to - \hat x$, we have
\begin{align}
\int _{\frac {\sigma}{2}}^{a - \frac {\sigma}{2}}
\Theta \left( - l \cos \phi + x + \frac {\sigma}{2} - a \right) {\rm d}x
& = \int _{a - \sigma}^{0} \Theta \left( - l \cos \phi - \hat x \right) {\rm d} (- \hat x) \nonumber \\
& = \int _{0}^{a - \sigma} \Theta \left( - l \cos \phi - \hat x \right) {\rm d} \hat x.
\end{align}
With $y - \frac{\sigma}{2} \to \hat y$, we have
\begin{equation}
\int _{\frac {\sigma}{2}}^{b - \frac {\sigma}{2}}
\Theta \left( l \sin \phi - y + \frac {\sigma}{2} \right) {\rm d}y
= \int _{0}^{b - \sigma} \Theta \left( l \sin \phi - \hat y \right) {\rm d} \hat y.
\end{equation}

With the definitions in Eqs. (\ref{eq:Ala-2D})-(\ref{eq:ABlab-2D}),
we reformulate Eqs. (\ref{eq:Scollx-2DSC})-(\ref{eq:Scollxy-2DSC}) as
\begin{align}
S_{\rm coll} (x)
& = b A(l, a - \sigma) + \sigma b \pi, \\
S_{\rm coll} (y)
& = a B(l, b - \sigma) + \sigma a \pi, \\
S_{\rm coll} (x, y)
& = AB(l, a - \sigma, b - \sigma) + \sigma A(l, a - \sigma) + \sigma B(l, b - \sigma) + \sigma^2 \pi.
\end{align}
Correspondingly, we have
\begin{equation}
\label{eq:Scoll-2DSC-reform}
S_{\rm coll}
= (b - \sigma) A(l, a - \sigma)
+ (a - \sigma) B(l, b - \sigma)
- AB(l, a - \sigma, b - \sigma)
+ \left( \sigma a + \sigma b - \sigma ^2 \right) \pi. 
\end{equation}

With the explicit forms of Eqs.(\ref{eq:Ala-2D})-(\ref{eq:ABlab-2D}) in Appendix A, we can directly derive the result of Eq.(\ref{eq:Scoll-2DSC-reform}), and finally $P(l, \sigma, a, b)$.
We leave details of calculation in Appendix C.
With Eqs.(\ref{eq:Plsab-2DSC-lb}), (\ref{eq:Plsab-2DSC-bla}), (\ref{eq:Plsab-2DSC-al1}), and (\ref{eq:Plsab-2DSC-al2}), we lay down our final equations as
\begin{equation}
\label{eq:Plsab-2DSC-full}
P(l, \sigma, a, b) = 
\left\{ \begin{array}{ll}
%l < b
\frac {2}{\pi} \frac {b - \sigma}{b} \frac {l}{a}
+ \frac {2}{\pi} \frac {a - \sigma}{a} \frac {l}{b} 
- \frac {1}{\pi} \frac {l^2}{ab}
 + \left( \frac {\sigma}{a} + \frac {\sigma}{b} - \frac {\sigma ^2}{ab} \right),
 & {\rm if}\ l \leqslant b - \sigma; \\
 %b < l < a
\frac {1}{\pi} \frac {(b - \sigma)^2}{ab}
+ \frac {2}{\pi} \frac {a - \sigma}{a} \frac {l}{b} \left[ 1 - \sqrt{1 - \frac {(b - \sigma)^2}{l^2}} \right] \\
+ \frac {2}{\pi} \frac {a - \sigma}{a} \frac {b - \sigma}{b} \left( \frac {\pi}{2} - \arcsin \frac {b - \sigma}{l} \right) 
 + \left( \frac {\sigma}{a} + \frac {\sigma}{b} - \frac {\sigma ^2}{ab}  \right),
& {\rm if}\ b - \sigma < l \leqslant a - \sigma; \\
%a < l < sqrt
\frac {1}{\pi} \left( \frac{\hat{L}^2}{ab} + \frac{l^{2}}{ab} \right)
- \frac {2}{\pi} \left[ \frac {b - \sigma}{b} \frac {l}{a} \sqrt {1- \frac {(a - \sigma)^2}{l^2}} 
+ \frac {a - \sigma}{a} \frac {l}{b} \sqrt {1- \frac {(b - \sigma)^2}{l^2} } \right] \\
+ \frac {2}{\pi} \frac {a - \sigma}{a} \frac {b - \sigma}{b} \left( \pi - \arcsin \frac {a - \sigma}{l} - \arcsin \frac {b - \sigma}{l} \right)
+ \left( \frac {\sigma}{a} + \frac {\sigma}{b} - \frac {\sigma^2}{ab} \right),
& {\rm if}\ a - \sigma < l \leqslant \hat{L}; \\
%sqrt < l
1, & {\rm if}\ \hat{L} < l,
\end{array}\right.
\end{equation}
in which we define $\hat{L} \equiv \sqrt {(a - \sigma)^2 + (b - \sigma)^2}$.

%special case: l > L and P = 1
We find that when $l > \hat{L}$, $P(l, \sigma, a, b) = 1$. 
This result has a similar intuitive interpretation with the case of $l > L$ in the $2$D needle version.
When a $2$D spherocylinder can be fitted in a grid cell, its largest length corresponds to the situation when it aligns with the diagonal line of the grid cell and its two semicircles touch the grid boundaries at two opposite corners. This length of a spherocylinder is simply $\hat {L}$. Thus any $2$D spherocylinder with a length $l > \hat{L}$ dropped on a grid is certain to intersect with grid lines, leading to $P(l, \sigma, a, b) = 1$.

%check boundary cases
We then check the consistency of Eq.(\ref{eq:Plsab-2DSC-full}) for $P(l, \sigma, a, b)$, and equivalently $S_{\rm coll}$, at the boundaries of the four cases of $l$.
We leave the details of calculation in Appendix D.

% BLNP --> BNP
When $a \to \infty$, Eq.(\ref{eq:Plsab-2DSC-full}) simplifies to the result for the BNP of $2$D spherocylinders as
\begin{equation}
\label{eq:Plsb-2DSC-full}
P(l, \sigma, \infty, b)
= \left\{ \begin{array}{ll}
 \frac {2}{\pi} \frac {l}{b} + \frac {\sigma}{b}, & {\rm if}\ l \leqslant b - \sigma; \\
 \frac {2}{\pi} \frac {l}{b} \left[ 1 - \sqrt {1- \frac {(b - \sigma)^2}{l^2}} \right] 
 + \frac {2}{\pi} \frac {b - \sigma}{b} \left( \frac {\pi}{2} - \arcsin \frac {b - \sigma}{l} \right)
 +\frac {\sigma}{b}, & {\rm if}\ b - \sigma < l.
 \end{array}\right.
 \end{equation}
%

%%%%%%%%%%%%%%%%%%%%%%%%%%%%%%%%
\section{Theory of BLNP with $3$D Needles}
\label{sec:theory-3D}

We follow the analytical framework for the BLNP with $2$D needles.
We define $S_{\rm all}^{\rm 3D}$ as the volume of total parameter space of $(x_{\rm U}, y_{\rm U}, \psi, \phi)$ in the case of $3$D needles, and $S_{\rm coll}^{\rm 3D}$ as the volume of proper parameter space of $(x_{\rm U}, y_{\rm U}, \psi, \phi)$ where an intersection happens. The intersection probability is thus
\begin{equation}
P^{\rm 3D}(l, a, b) = \frac {S_{\rm coll}^{\rm 3D}}{S_{\rm all}^{\rm 3D}}.
\end{equation}
For $S_{\rm all}^{\rm 3D}$, we easily have
\begin{equation}
S_{\rm all}^{\rm 3D} = a b \pi \frac {\pi}{2}.
\end{equation}
$S_{\rm coll}^{\rm 3D}$ can be reformulated as
\begin{equation}
\label{eq:Scoll-3D}
S_{\rm coll}^{\rm 3D} = S_{\rm coll}^{\rm 3D}(x) + S_{\rm coll}^{\rm 3D}(y) - S_{\rm coll}^{\rm 3D} (x, y),
\end{equation}
while $S_{\rm coll}^{\rm 3D}(x)$, $S_{\rm coll}^{\rm 3D}(y)$, and $S_{\rm coll}^{\rm 3D}(x, y)$ denote the volume of the parameter space of $(x_{\rm U}, y_{\rm U}, \psi, \phi)$ when an intersection between the projection of a $3$D needle and a grid cell happens on the boundaries of the grid cell along the $y$-axis, the $x$-axis, and the $x$- and $y$-axes at the same time, respectively.
Based on Eqs.(\ref{eq:configuration-3D-uppertip})-(\ref{eq:intersection-needle}), we substitute $l$ with $l \sin \psi$ and put an extra $\int _{0}^{\pi/2} {\rm d} \psi$ on each term in Eqs.(\ref{eq:Scollx-2D})-(\ref{eq:Scollxy-2D}). We have
\begin{align}
%Scollx
S_{\rm coll}^{\rm 3D} (x)
& = b  \int _{0}^{\frac{\pi}{2}} {\rm d} \psi
\int _{0}^{\pi} {\rm d}\phi
\int _{0}^{a} \Theta(|l \sin \psi \cos \phi| - x) {\rm d}x, \\
%Scoll_y
S_{\rm coll}^{\rm 3D} (y)
& = a \int _{0}^{\frac{\pi}{2}} {\rm d} \psi
\int _{0}^{\pi} {\rm d}\phi
\int _{0}^{b} \Theta(l \sin \psi \sin \phi - y) {\rm d}y, \\
%Scoll_xy
S_{\rm coll}^{\rm 3D} (x,y)
& = \int _{0}^{\frac{\pi}{2}} {\rm d} \psi
\int _{0}^{\pi} {\rm d}\phi
\int _{0}^{a} \Theta(|l \sin \psi \cos \phi| - x) {\rm d}x
\int _{0}^{b} \Theta(l \sin \psi \sin \phi - y) {\rm d}y.
\end{align}

To compute the above equations, we define integrals like Eqs.(\ref{eq:Ala-2D})-(\ref{eq:ABlab-2D}) as
\begin{align}
\label{eq:Ala-3D}
%A(l, a)
A^{\rm 3D}(l, a)
& \equiv \int _{0}^{\frac {\pi}{2}} {\rm d} \psi 
\int _{0}^{\pi} {\rm d} \phi 
\int _{0}^{a} \Theta (|l \sin \psi \cos \phi | - x) {\rm d} x, \\
\label{eq:Blb-3D}
%B(l, b)
B^{\rm 3D}(l, b)
& \equiv \int _{0}^{\frac {\pi}{2}} {\rm d} \psi 
\int _{0}^{\pi} {\rm d} \phi 
\int _{0}^{b} \Theta (l \sin \psi \sin \phi - y) {\rm d} y, \\
\label{eq:ABlab-3D}
%AB(l, a, b)
AB^{\rm 3D} (l, a, b)
& \equiv  \int _{0}^{\frac{\pi}{2}} {\rm d} \psi
\int _{0}^{\pi} {\rm d}\phi
\int _{0}^{a} \Theta(|l \sin \psi \cos \phi| - x) {\rm d}x
\int _{0}^{b} \Theta(l \sin \psi \sin \phi - y) {\rm d}y.
\end{align}
Thus we reformulate $S_{\rm coll}^{\rm 3D}$ as
\begin{equation}
\label{eq:Scoll-3D-reform}
S_{\rm coll}^{\rm 3D} = b A^{\rm 3D}(l, a) + a B^{\rm 3D}(l, b) - AB^{\rm 3D}(l, a, b).
\end{equation}

To calculate above equations, a thorough analysis of the relative size among $l$, $a$, and $b$ should be carried out, in which $\psi$ is considered in different regions in integration.
We leave details of calculation in Appendix E.
With Eqs.(\ref{eq:Plab-3D-lb}), (\ref{eq:Plab-3D-bla}), (\ref{eq:Plab-3D-al1}), and (\ref{eq:Plab-3D-al2}),
we lay down our final equations as
\begin{equation}
\label{eq:Plab-3D-full}
P^{\rm 3D}(l, a, b) =
\left\{ \begin{array}{ll}
% l < b
\frac {4}{\pi ^2} \left( \frac {l}{a} + \frac {l}{b} \right) - \frac {1}{2 \pi} \frac {l^2}{ab},
 & {\rm if}\ l \leqslant b; \\
 % b < l < a
\frac {4}{\pi ^2} \frac {l}{a}
+ \left[ \frac {4}{\pi ^2} \frac {l}{b}
- \frac {4}{\pi ^2} \frac {l}{b} F \left( \arcsin \frac {b}{l}, \frac {\pi}{2}, \frac {b}{l} \right)
+ \frac {4}{\pi ^2} G\left( \arcsin \frac {b}{l}, \frac {\pi}{2}, \frac {b}{l} \right) \right] \\
- \left[ \frac {1}{\pi ^2} \frac {l^2}{ab} \arcsin \frac {b}{l}
+ \frac {3}{\pi ^2} \frac {l}{a} \sqrt {1 - \frac {b^2}{l^2}}
- \frac {2}{\pi ^2} \frac {b}{a} \left( \frac {\pi}{2} - \arcsin \frac {b}{l} \right) \right],
& {\rm if}\ b < l \leqslant a; \\
% a < l < L
\left[ \frac {4}{\pi ^2} \frac {l}{a}
- \frac {4}{\pi ^2} \frac {l}{a} F \left( \arcsin \frac {a}{l}, \frac {\pi}{2}, \frac {a}{l} \right)
+ \frac {4}{\pi ^2} G\left( \arcsin \frac {a}{l}, \frac {\pi}{2}, \frac {a}{l} \right)  \right] \\
+ \left[  \frac {4}{\pi ^2} \frac {l}{b}
-  \frac {4}{\pi ^2} \frac {l}{b} F \left(  \arcsin \frac {b}{l}, \frac {\pi}{2}, \frac {b}{l} \right)
+ \frac {4}{\pi ^2} G\left(  \arcsin \frac {b}{l}, \frac {\pi}{2}, \frac {b}{l} \right) \right] \\
- \left[ \frac {1}{\pi ^2} \frac {l^2}{ab} \arcsin \frac {a}{l}
+ \frac {3}{\pi ^2} \frac {l}{b} \sqrt {1 - \frac {a^2}{l^2}}
- \frac {2}{\pi ^2} \frac {a}{b} \left( \frac {\pi}{2} - \arcsin \frac {a}{l} \right) \right] \\
- \left[ \frac {1}{\pi ^2} \frac {l^2}{ab}\arcsin \frac {b}{l}
+ \frac {3}{\pi ^2} \frac {l}{a} \sqrt {1 - \frac {b^2}{l^2}}
- \frac {2}{\pi ^2} \frac {b}{a} \left( \frac {\pi}{2} - \arcsin \frac {b}{l} \right) \right]
+ \frac {1}{2 \pi} \frac {l^2}{ab},
& {\rm if}\ a < l \leqslant L ; \\
% L < l
\left[ \frac {4}{\pi ^2} \frac {l}{a}
- \frac {4}{\pi ^2} \frac {l}{a} F \left( \arcsin \frac {a}{l}, \arcsin \frac {L}{l}, \frac {a}{l} \right)
+ \frac {4}{\pi ^2} G \left( \arcsin \frac {a}{l}, \arcsin \frac {L}{l}, \frac {a}{l} \right) \right] \\
+ \left[ \frac {4}{\pi ^2} \frac {l}{b}
- \frac {4}{\pi ^2} \frac {l}{b} F \left( \arcsin \frac {b}{l}, \arcsin \frac {L}{l}, \frac {b}{l} \right)
+ \frac {4}{\pi ^2} G \left( \arcsin \frac {b}{l}, \arcsin \frac {L}{l}, \frac {b}{l} \right) \right] \\
- \left[ \frac {1}{\pi ^2} \frac {l^2}{ab} \arcsin \frac {a}{l}
+ \frac {3}{\pi ^2} \frac {l}{b} \sqrt {1 - \frac {a^2}{l^2}}
- \frac {2}{\pi ^ 2} \frac {a}{b} \left( \frac {\pi}{2} - \arcsin \frac {a}{l} \right) \right] \\
- \left[ \frac {1}{\pi ^2} \frac {l^2}{ab} \arcsin \frac {b}{l}
+ \frac {3}{\pi ^2} \frac {l}{a} \sqrt {1 - \frac {b^2}{l^2}}
- \frac {2}{\pi ^2} \frac {b}{a} \left( \frac {\pi}{2} - \arcsin \frac {b}{l} \right) \right] \\
+ \left[ \frac {1}{\pi ^2} \frac {l^2}{ab} \arcsin \frac {L}{l}
- \frac {1}{\pi ^2} \frac {l L}{ab} \sqrt {1 - \frac {L^2}{l^2}} 
- \frac {2}{\pi ^2} \frac {L^2}{ab} \left( \frac {\pi}{2} - \arcsin \frac {L}{l} \right) \right] 
+ \frac {2}{\pi} \left( \frac {\pi}{2} - \arcsin \frac {L}{l} \right),
& {\rm if}\ L < l,
\end{array}\right.
\end{equation}
in which two integrals $F(a, b, t)$ and $G(a, b, t)$ are defined as Eqs.(\ref{eq:Fabt}) and (\ref{eq:Gabt}), respectively.

%consistency on boundaries
We can check the consistency of Eq.(\ref{eq:Plab-3D-full}) for $P^{\rm 3D}(l, a, b)$, and equivalently $S^{\rm 3D}_{\rm coll}$, at the boundaries of the four cases of $l$. We leave details of calculation in Appendix F.

% l -> infty
We then consider the case of $l \to \infty$.
Calculation in Appendix G shows that, only when $l \to \infty$, we have $P^{\rm 3D}(\infty, a, b) = 1$.
This is intuitively understandable. The projection of a finite long needle onto a $2$D plane eludes the boundaries of a grid cell with a non-vanishing probability, since the orientation $\psi$ always has a finite, even narrow, range to trigger an intersection. Only an infinitely long needle is sure to intersect with the grid.

%BLNP --> BNP
When $a \to \infty$, Eq.(\ref{eq:Plab-3D-full}) simplifies to the result for the BNP with $3$D needles as
\begin{equation}
\label{eq:Plb-3D-full}
P^{\rm 3D}(l, \infty, b)
= \left\{ \begin{array}{ll}
\frac {4}{\pi ^2} \frac {l}{b}, & {\rm if}\ l \leqslant b; \\
\frac {4}{\pi ^2} \frac {l}{b}
-  \frac {4}{\pi ^2} \frac {l}{b} F \left( \arcsin \frac {b}{l}, \frac {\pi}{2}, \frac {b}{l}  \right)
+ \frac {4}{\pi ^2} G \left( \arcsin \frac {b}{l}, \frac {\pi}{2}, \frac {b}{l} \right), & {\rm if}\ b < l.
\end{array}\right.
\end{equation}
%

%%%%%%%%%%%%%%%%%%%%%%%%%%%%%%%%
\section{Theory of BLNP with $3$D Spherocylinders}
\label{sec:theory-3DSC}

We follow the framework for the BLNP with $3$D needles, and calculate the intersection probability for the BLNP with $3$D spherocylinders here as
\begin{equation}
P^{\rm 3D}(l, \sigma, a, b) = \frac {S_{\rm coll}^{\rm 3D}}{S_{\rm all}^{\rm 3D}}.
\end{equation}
$S_{\rm coll}^{\rm 3D}$ can be further reformulated as Eq.(\ref{eq:Scoll-3D}). $S^{\rm 3D}_{\rm coll}(x)$, $S^{\rm 3D}_{\rm coll}(y)$, and $S^{\rm 3D}_{\rm coll}(x, y)$ denote the volume of the parameter space of $(x_{\rm U}, y_{\rm U}, \psi, \phi)$ when an intersection between the projection of a $3$D spherocylinder and a grid cell happens on the boundaries of the grid cell along the $y$-axis, the $x$-axis, and the $x$- and $y$-axes at the same time, respectively.
Based on Eqs. (\ref{eq:configuration-3D-uppertip}), (\ref{eq:configuration-3D-lowertip}), and (\ref{eq:intersection-SC}), we move from $l$ to $l \sin \psi$ and add an extra $\int _{0}^{\pi / 2} {\rm d}\psi$ to each term of Eqs.(\ref{eq:Scollx-2DSC})-(\ref{eq:Scollxy-2DSC}). We have
\begin{align}
%S_{coll}(x)
S_{\rm coll}^{\rm 3D} (x)
& = b \int _{0}^{\frac {\pi}{2}} {\rm d} \psi
\int _{0}^{\pi}  {\rm d}\phi
\int _{0}^{a - \sigma} \Theta \left(| l \sin \psi \cos \phi | - x \right) {\rm d}x
+ \sigma b \frac {\pi ^2}{2}. \\
%S_{coll}(y)
S_{\rm coll}^{\rm 3D} (y)
& = a \int _{0}^{\frac {\pi}{2}} {\rm d} \psi
\int _{0}^{\pi} {\rm d} \phi \int _{0}^{b - \sigma} \Theta(l \sin \psi \sin \phi - y) {\rm d}y
+ \sigma a \frac {\pi ^2}{2}, \\
%S_{coll}(x,y)
S_{\rm coll}^{\rm 3D} (x, y)
& = \int _{0}^{\frac {\pi}{2}} {\rm d} \psi
\int _{0}^{\pi} {\rm d}\phi
\int _{0}^{a - \sigma} \Theta (| l \sin \psi \cos \phi | - x) {\rm d}x
\int _{0}^{b - \sigma} \Theta (l \sin \psi \sin \phi - y) {\rm d}y \nonumber \\
& + \sigma \int _{0}^{\frac {\pi}{2}} {\rm d} \psi
\int _{0}^{\pi} {\rm d}\phi
\int _{0}^{a - \sigma} \Theta (| l \sin \psi \cos \phi | - x) {\rm d}x
+ \sigma \int _{0}^{\frac {\pi}{2}} {\rm d} \psi
\int _{0}^{\pi} {\rm d}\phi
\int _{0}^{b - \sigma} \Theta (l \sin \psi \sin \phi - y) {\rm d}y \nonumber \\
& + \sigma ^2 \frac {\pi ^2}{2} .
\end{align}

With the definitions in Eqs.(\ref{eq:Ala-3D})-(\ref{eq:ABlab-3D}), we have
\begin{align}
S_{\rm coll}^{\rm 3D} (x)
& = b A^{\rm 3D}(l, a - \sigma) + \sigma b \frac {\pi ^2}{2}, \\
S_{\rm coll}^{\rm 3D} (y)
& = a B^{\rm 3D}(l, b - \sigma) + \sigma a \frac {\pi ^2}{2}, \\
S_{\rm coll}^{\rm 3D} (x, y)
& = AB^{\rm 3D}(l, a - \sigma, b - \sigma)
+ \sigma A^{\rm 3D}(l, a - \sigma) + \sigma B^{\rm 3D}(l, b - \sigma)
+ \sigma ^2 \frac {\pi ^2}{2}.
\end{align}
Equivalently, we have
\begin{equation}
\label{eq:Scoll-3DSC-reform}
S_{\rm coll}^{\rm 3D} 
= (b - \sigma) A^{\rm 3D}(l, a - \sigma) + (a - \sigma) B^{\rm 3D}(l, b - \sigma)
- AB^{\rm 3D}(l, a - \sigma, b - \sigma)
+ \left( \sigma a + \sigma b - \sigma ^2 \right) \frac {\pi ^2}{2}.
\end{equation}

With the explicit forms of $A^{\rm 3D}(l, a)$, $B^{\rm 3D}(l, b)$ and $AB^{\rm 3D}(l, a, b)$ in Appendix E, we can easily lay down $S_{\rm coll}^{\rm 3D} $ with Eq.(\ref{eq:Scoll-3DSC-reform}), and finally $P^{\rm 3D}(l, \sigma, a, b)$.
We leave details of calculation in Appendix H.
With Eqs.(\ref{eq:Plsab-3DSC-lb}), (\ref{eq:Plsab-3DSC-bla}), (\ref{eq:Plsab-3DSC-al1}), and (\ref{eq:Plsab-3DSC-al2}), we list our final equations as
\begin{eqnarray}
\label{eq:Plsab-3DSC-full}
&& P^{\rm 3D}(l, \sigma, a, b) =  \nonumber \\
&& \left\{ \begin{array}{ll}
% l < b
\frac {4}{\pi ^2} \frac {b - \sigma}{b} \frac {l}{a}
+ \frac {4}{\pi ^2} \frac {a - \sigma}{a} \frac {l}{b}
- \frac {1}{2 \pi} \frac {l^2}{ab}
+ \left( \frac {\sigma}{a} + \frac {\sigma}{b} - \frac {\sigma ^2}{ab} \right),
 & {\rm if}\ l \leqslant b - \sigma; \\
 % b < l < a
\frac {4}{\pi ^2} \frac {b - \sigma}{b} \frac {l}{a}
+ \left[ \frac {4}{\pi ^2} \frac {a - \sigma}{a} \frac {l}{b}
- \frac {4}{\pi ^2} \frac {a - \sigma}{a} \frac {l}{b}F \left( \arcsin \frac {b - \sigma}{l}, \frac {\pi}{2}, \frac {b - \sigma}{l} \right)
+ \frac {4}{\pi ^2} \frac{a - \sigma}{a} \frac{b - \sigma}{b} G \left( \arcsin \frac {b - \sigma}{l}, \frac {\pi}{2}, \frac {b - \sigma}{l} \right) \right] \\
 - \left[ \frac {1}{\pi ^2} \frac {l^2}{ab} \arcsin \frac {b - \sigma}{l}
+ \frac {3}{\pi ^2}  \frac {b - \sigma}{b} \frac {l}{a} \sqrt {1 - \frac {(b - \sigma)^2}{l^2}}
- \frac {2}{\pi ^2} \frac {(b - \sigma)^2}{ab} \left( \frac {\pi}{2} - \arcsin \frac {b - \sigma}{l} \right) \right]
+ \left( \frac {\sigma}{a} + \frac{\sigma}{b} - \frac {\sigma ^2}{ab} \right),
& {\rm if}\ b - \sigma < l \leqslant a - \sigma; \\
% a < l < L
\left[ \frac {4}{\pi ^2} \frac {b - \sigma}{b} \frac {l}{a}
- \frac {4}{\pi ^2} \frac {b - \sigma}{b} \frac {l}{a}F \left( \arcsin \frac {a - \sigma}{l}, \frac {\pi}{2}, \frac {a - \sigma}{l} \right)
+ \frac {4}{\pi ^2} \frac{a - \sigma}{a} \frac{b - \sigma}{b} G \left( \arcsin \frac {a - \sigma}{l}, \frac {\pi}{2}, \frac {a - \sigma}{l} \right) \right] \\
 + \left[ \frac {4}{\pi ^2} \frac {a - \sigma}{a} \frac {l}{b}
- \frac {4}{\pi ^2} \frac {a - \sigma}{a} \frac {l}{b}F \left( \arcsin \frac {b - \sigma}{l}, \frac {\pi}{2}, \frac {b - \sigma}{l} \right)
+ \frac {4}{\pi ^2} \frac{a - \sigma}{a} \frac{b - \sigma}{b} G \left( \arcsin \frac {b - \sigma}{l}, \frac {\pi}{2}, \frac {b - \sigma}{l} \right) \right] \\
 - \left[ \frac {1}{\pi ^2} \frac {l^2}{ab} \arcsin \frac {a - \sigma}{l}
+ \frac {3}{\pi ^2}  \frac {a - \sigma}{a} \frac {l}{b} \sqrt {1 - \frac {(a - \sigma)^2}{l^2}}
- \frac {2}{\pi ^2} \frac {(a - \sigma)^2}{ab} \left( \frac {\pi}{2} - \arcsin \frac {a - \sigma}{l} \right) \right] \\
 - \left[ \frac {1}{\pi ^2} \frac {l^2}{ab} \arcsin \frac {b - \sigma}{l}
+ \frac {3}{\pi ^2}  \frac {b - \sigma}{b} \frac {l}{a} \sqrt {1 - \frac {(b - \sigma)^2}{l^2}}
- \frac {2}{\pi ^2} \frac {(b - \sigma)^2}{ab} \left( \frac {\pi}{2} - \arcsin \frac {b - \sigma}{l} \right) \right] \\
+ \frac {1}{2\pi} \frac {l^2}{ab} + \left( \frac {\sigma}{a} + \frac{\sigma}{b} - \frac {\sigma ^2}{ab} \right),
& {\rm if}\ a - \sigma < l \leqslant \hat {L} ; \\
% L < l
\left[ \frac {4}{\pi ^2} \frac {b - \sigma}{b} \frac {l}{a}
- \frac {4}{\pi ^2} \frac {b - \sigma}{b} \frac {l}{a} F \left( \arcsin \frac {a - \sigma}{l}, \arcsin \frac {\hat{L}}{l}, \frac {a - \sigma}{l} \right)
+ \frac {4}{\pi ^2} \frac {a - \sigma}{a} \frac {b - \sigma}{b} G \left( \arcsin \frac {a - \sigma}{l}, \arcsin \frac {\hat{L}}{l}, \frac {a - \sigma}{l} \right) \right]  \\
+ \left[ \frac {4}{\pi ^2} \frac {a - \sigma}{a} \frac {l}{b}
- \frac {4}{\pi ^2} \frac {a - \sigma}{a} \frac {l}{b} F \left( \arcsin \frac {b - \sigma}{l}, \arcsin \frac {\hat{L}}{l}, \frac {b - \sigma}{l} \right)
+ \frac {4}{\pi ^2} \frac {a - \sigma}{a} \frac {b - \sigma}{b} G \left( \arcsin \frac {b - \sigma}{l}, \arcsin \frac {\hat{L}}{l}, \frac {b - \sigma}{l} \right) \right] \\
- \left[ \frac {1}{\pi ^2} \frac {l^2}{ab} \arcsin \frac {a - \sigma}{l}
+ \frac {3}{\pi ^2} \frac {a - \sigma}{a} \frac {l}{b} \sqrt {1 - \frac {(a - \sigma)^2}{l^2}}
- \frac {2}{\pi ^2} \frac {(a - \sigma)^2}{ab} \left( \frac {\pi}{2} - \arcsin \frac {a - \sigma}{l} \right) \right] \\
- \left[ \frac {1}{\pi ^2} \frac {l^2}{ab} \arcsin \frac {b - \sigma}{l}
+ \frac {3}{\pi ^2} \frac {b - \sigma}{b} \frac {l}{a} \sqrt {1 - \frac {(b - \sigma)^2}{l^2}}
- \frac {2}{\pi ^2} \frac {(b - \sigma)^2}{ab} \left( \frac {\pi}{2} - \arcsin \frac {b - \sigma}{l} \right) \right] \\
+ \left[ \frac {1}{\pi ^2} \frac {l^2}{ab} \arcsin \frac {\hat{L}}{l}
- \frac {1}{\pi ^2} \frac {l \hat {L}}{ab} \sqrt {1 - \frac {\hat{L}^2}{l^2}}
- \frac {2}{\pi ^2} \frac {\hat{L}^2}{ab} \left( \frac {\pi}{2} - \arcsin \frac {\hat{L}}{l} \right) \right] \\
+ \frac {2}{\pi} \frac {a - \sigma}{a} \frac {b - \sigma}{b} \left( \frac {\pi}{2} - \arcsin \frac {\hat{L}}{l} \right)
+ \left( \frac{\sigma}{a} + \frac {\sigma}{b} - \frac {\sigma ^2}{ab} \right),
& {\rm if}\ \hat{L} < l.
\end{array}\right.
\end{eqnarray}
%

%check consistency on boundaries
We further check the consistency of Eq.(\ref{eq:Plsab-3DSC-full})  for $P^{\rm 3D}(l, \sigma, a, b)$, and equivalently $S^{\rm 3D}_{\rm coll}$, on the boundaries of four cases of $l$. We leave details of calculation in Appendix I.

%case of infinitely large l
We then consider the limit case of $l \to \infty$.
In Appendix J, we show that, just like the case of infinitely long $3$D needles, only a $3$D spherocylinder with an infinitely large length $l$ intersect with a grid with a probability of $1$.

%BLNP --> BNP
When $a \to \infty$, Eq. (\ref{eq:Plsab-3DSC-full}) simplifies to the result for the BNP of $3$D spherocylinders as
\begin{equation}
\label{eq:Plsb-3DSC-full}
P^{\rm 3D}(l, \sigma, \infty, b)
 = \left\{ \begin{array}{ll}
 \frac {4}{\pi ^2} \frac {l}{b} + \frac {\sigma}{b}, & {\rm if}\ l \leqslant b - \sigma; \\
 \frac {4}{\pi ^2} \frac {l}{b}
-  \frac {4}{\pi ^2} \frac {l}{b} F \left( \arcsin \frac {b - \sigma}{l}, \frac {\pi}{2}, \frac {b - \sigma}{l}  \right) 
+ \frac {4}{\pi ^2} \frac {b - \sigma}{b} G \left( \arcsin \frac {b - \sigma}{l}, \frac {\pi}{2}, \frac {b - \sigma}{l} \right)
+ \frac {\sigma}{b}, & {\rm if}\ b - \sigma < l.
 \end{array} \right.
\end{equation}

\section{Comparison with Previous Results}
\label{sec:theory-comparison}

%(1)BLNP: l < min(a,b)
First, our Eq.(\ref{eq:Plab-2D-full}) with $l \leqslant b$ retrieves the result for the `short' needle case of $l < b$ for BLNP with $2$D needles in \cite{Laplace-1812, Laplace-1820, Uspensky-1937}.

%(1)BNP: any l
Then, we can see that Eq.(\ref{eq:Plb-2D-full}), which corresponds to the case of $a \to \infty$ in our theory for BLNP with $2$D needles, corresponds to the result for the original BNP with arbitrary length $l$ in \cite{Uspensky-1937}.

%(3)BLNP: l > 1, a=b=1
We further compare our results on BLNP with those in \cite{Dell.Franklin-JSTAT-2009} with a special case of $a = b = 1$.
Eq.(6) in \cite{Dell.Franklin-JSTAT-2009} shows the intersection probability in BLNP with $2$D needles under $b < l$. Based on Eq.(\ref{eq:Plab-2D-full}) in the case of $a < l \leqslant L$ and setting $a = b = 1$, we have
\begin{equation}
P(l, 1, 1) 
= 1 - \frac {2}{\pi} \left( \arcsin \frac {1}{l} - \arccos \frac {1}{l}
+ 2 \sqrt {l^2 - 1} - \frac {l^2}{2} - 1 \right).
\end{equation}
We can easily find that there is an error in Eq.(6) in \cite{Dell.Franklin-JSTAT-2009} as the coefficient $1 / \pi$ of the second term should be $2 / \pi$ as we show above.

%(4)BLNP with 2D spherocylinder: l < 1, a=b=1
Besides, Eq.(8) in \cite{Dell.Franklin-JSTAT-2009} presents the intersection probability in BLNP with $2$D spherocylinders in the case of small $l$. With Eq.(\ref{eq:Plsab-2DSC-full}) in the case of $l \leqslant b - \sigma$ and setting $a = b = 1$, we have
\begin{equation}
P(l, \sigma, 1, 1)
= \frac {4l - l^2}{\pi} + \frac {\sigma}{\pi} (- 4 l + 2 \pi - \sigma \pi).
\end{equation}
The difference between Eq.(8) in \cite{Dell.Franklin-JSTAT-2009} and our above prediction is
\begin{equation}
\Delta P = \frac {\sigma}{\pi} \left(- 4 + 3 \sigma + 2 l \right).
\end{equation}
It is easy to see that, when $\sigma \ll 1$, we have $\Delta P < 0$.
Yet when $\sigma \approx l \approx 1$, we have instead $\Delta P > 0$.
As all our analytical results can be validated in the following section, we believe that there is probably an error in the derivation of Eq.(8) in \cite{Dell.Franklin-JSTAT-2009}.

%%%%%%%%%%%%%%%%%%%%%%%%%%%%%%%%
\section{Result}
\label{sec:result}

\begin{figure}[htbp]
\begin{center}
 \includegraphics[width = 0.95 \linewidth]{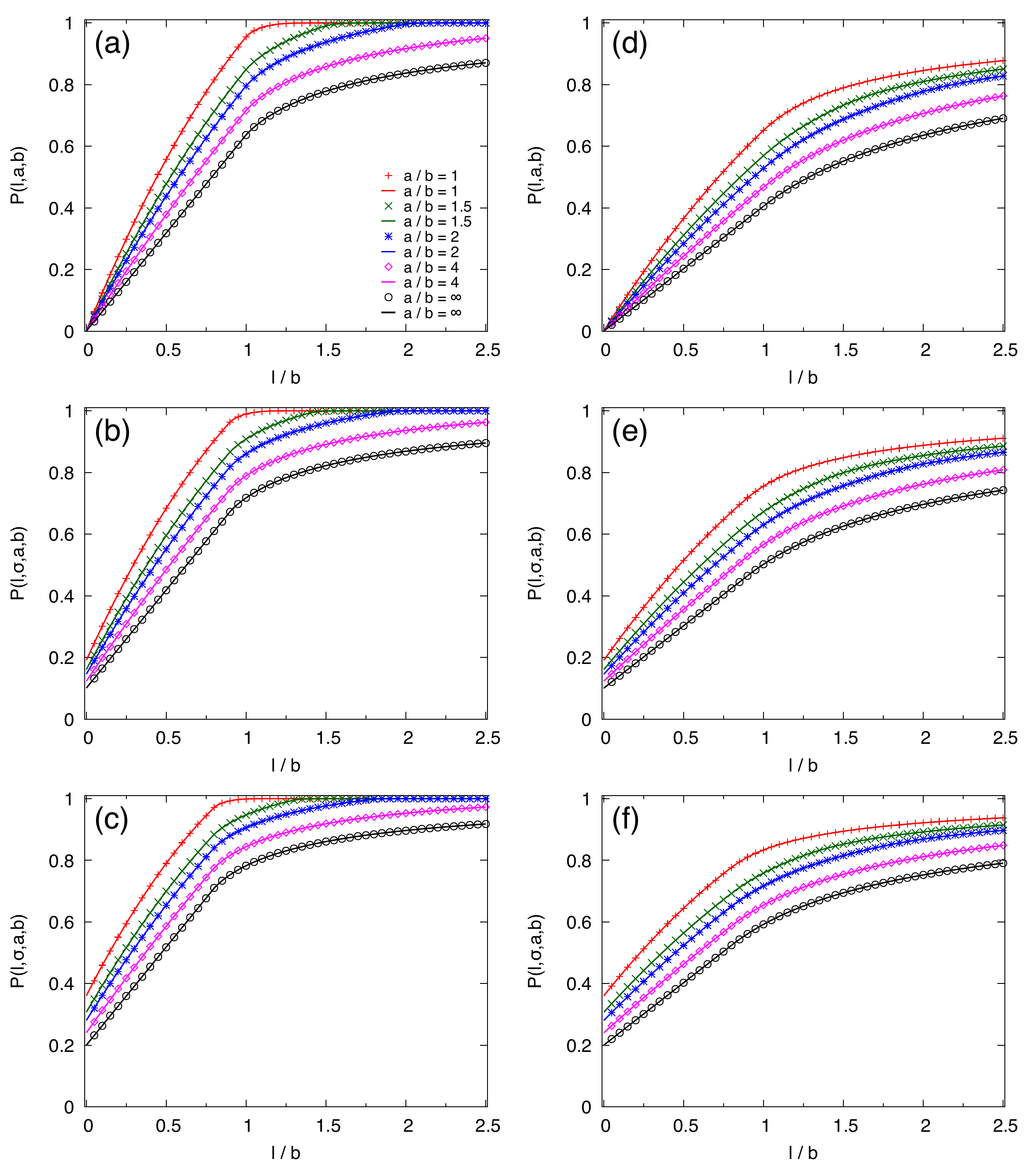}
\end{center}
\caption{
 \label{fig:result-lb}
Result of BLNP with needles and spherocylinders with fixed height $b$ of a grid cell.
(a) The intersection probability $P(l, a, b)$ versus the ratio $l / b$ is shown in the case of $2$D needles.
(b)-(c) The intersection probability $P(l, \sigma, a, b)$ versus $l / b$ is shown in the case of $2$D spherocylinders with $\sigma / b = 0.1$ and $0.2$, respectively.
(d)-(f) Results are shown in the same format as (a)-(c), yet for the case of $3$D needles and spherocylinders.
In each subfigure, we show results when $a/b = 1, 1.5, 2, 4, \infty$,
the last of which simply corresponds to the BNP version.
Each sign is from result of a single Monte Carlo simulation with $b = 3$ and a configuration size $N_{\rm all} = 10^6$ for the $2$D case and $N_{\rm all} = 10^7$ for the $3$D case.
Solid lines are from analytical results.}
\end{figure}
\begin{figure}[htbp]
\begin{center}
 \includegraphics[width = 0.95 \linewidth]{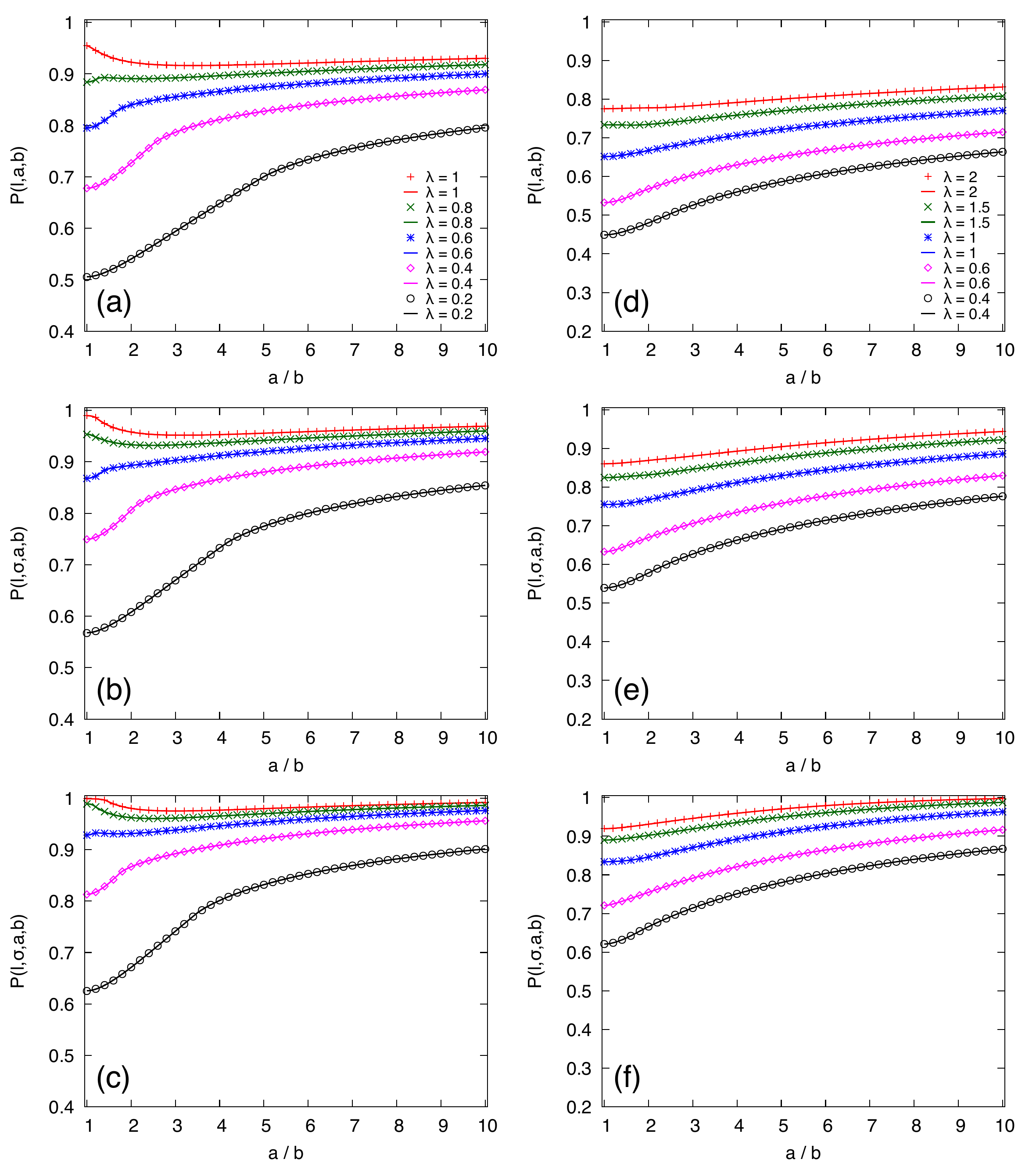}
\end{center}
\caption{
 \label{fig:result-ab}
Result of BLNP with needles and spherocylinders with fixed area $ab$ of a grid cell.
(a) The intersection probability $P(l, a, b)$ versus the ratio $a / b$ is shown in the case of $2$D needles.
(b)-(c) The intersection probability $P(l, \sigma, a, b)$ versus $a / b$ is shown in the case of $2$D spherocylinders with $\sigma / l = 0.1$ and $0.2$, respectively.
In each subfigure of (a)-(c), we show results when $\lambda = 1, 0.8, 0.6, 0.4, 0.2$.
(d)-(f) Results are shown in the same format as (a)-(c), yet for the cases of $3$D needles and spherocylinders.
In each subfigure of (d)-(f), we show results when $\lambda = 2, 1.5, 1, 0.6, 0.4$.
Each sign is from result of a single Monte Carlo simulation with $l = 3$ and a configuration size $N_{\rm all} = 10^6$ for the $2$D case and $N_{\rm all} = 10^7$ for the $3$D case.
Solid lines are from analytical results.}
\end{figure}

Our final analytical equations are summed in Eqs.(\ref{eq:Plab-2D-full}), (\ref{eq:Plb-2D-full}), (\ref{eq:Plsab-2DSC-full}), (\ref{eq:Plsb-2DSC-full}),
(\ref{eq:Plab-3D-full}),  (\ref{eq:Plb-3D-full}),
(\ref{eq:Plsab-3DSC-full}), and (\ref{eq:Plsb-3DSC-full}).
There are two points we should mention.
The first one is that, the intersection probabilities depend only on the relative sizes of dropped objects and grids, say $l$, $\sigma$, $a$, and $b$, not their sizes per se.
The second one is that, to calculate an intersection probability there are four regimes separated by three length scales, each with a specific analytical form.
To verify the correctness of our analytical predictions, we perform Monte Carlo simulation of dropping objects for given finite $l$, $\sigma$, $a$, and $b$ explained in Section \ref{sec:model}, and compare these empirical intersection probability with theoretical predictions from our analytical framework.

%Simpson's rule
To numerically calculate integrals, there is a large toolbox in existing literature we can refer to. Here for simplicity we adopt the classical Simpson's rule \cite{Press.etal-2007-3e} for Eqs. (\ref{eq:Fabt}) and (\ref{eq:Gabt}). We can see that the integrands in both integrals have a closed form. We divide the range $[a, b]$ into equally spaced intervals with a size $\lceil (b - a) \times N_{\rm unit} \times 2 i \rceil$ with $i = 1, 2, \cdots$, and calculate a sequence of numerical summations $F(i)$ based on Simpson's rule. When $| F(i +1) - F(i)| < \varepsilon$, we consider $F(i + 1)$ as our numerical approximation of an integral.
In this section, we set $N_{\rm unit} = 10^4$ and $\varepsilon = 10^{-9}$.

%result 1: varying l/b
We first consider how the sizes of needles and spherocylinders affect an intersection probability.
Fig.\ref{fig:result-lb} shows both results from numerical simulation and analytical theory.
%observation 1: larger size
As we can see, an intersection probability increases monotonously with a larger length for a needle or a larger length and radius for a spherocylinder once the other size parameters are fixed.
%observation 2: 2D to 3D
Besides, a comparison between Fig.\ref{fig:result-lb} (a)-(c) with (d)-(f) shows that, the introduction of a new degree of freedom for a dropped object generally pushes down intersection probability. The profile of intersection probability further changes drastically in the case of long objects:  in the $2$D cases, objects with a length beyond a finite critical value is sure to intersect with grid lines, while in the $3$D cases only infinitely long objects intersect with grid lines with probability $1$.

%a naive extrapolation of equations
We should also mention that, for our equations of intersection probabilities, a naive extrapolation of an analytical form in a certain region to a different one with a larger $l$ only leads to an overestimation of intersection probability, until there is an unphysical result as the probability goes beyond $1$.

%result 2: varying aspect ratio a/b
From a different perspective, we further consider how an intersection probability evolves with the aspect ratio of grid cell $a / b$. This context is highly relevant when we consider irregular filter in a filtration process.
With our theory, we can locate those parameters leading to the maximal flux (correspondingly the minimal intersection probability), which can be further tested in a controlled experiment.
We define a dimensionless control parameter $\lambda \equiv l^2 / ab$, and calculate intersection probability with tuned $\lambda$. Given the control coefficient $\lambda$, the radius ratio $\sigma/l \ (\equiv \sigma _{l})$, and the aspect ratio $a/b \ (\equiv t)$, a correspondence between parameters can be established as
\begin{equation}
\left( \frac {l}{a}, \frac {l}{b}, \frac {\sigma}{a}, \frac {\sigma}{b} \right)
= \left( \sqrt {\frac {\lambda}{t}}, \sqrt {\lambda t}, \sigma _{l} \sqrt {\frac {\lambda}{t}}, \sigma _{l} \sqrt {\lambda t}  \right).
\end{equation}
Fig.\ref{fig:result-ab} shows both simulation and analytical result. Unlike the pattern of monotonously increasing intersection probability in Fig.\ref{fig:result-lb}, we can find a much more complicated picture.

%2D results
%observation: scenario of moving minimum
Results in Fig.\ref{fig:result-ab} (a)-(c) for $2$D cases follow a quite similar pattern, as we can find a local minimum moving from $t = 1$ to some $t > 1$ with increasing $\lambda$. 
We present a possible qualitative scenario here:
(1) when $\lambda \leqslant \lambda _{1}^{\ast}$, there is only one local minimum of intersection probability, which is exactly at $t_{1} = 1$;
(2) when $\lambda _{1}^{\ast} < \lambda \leqslant \lambda _{2}^{\ast}$, a second local minimum emerges at $t_{2} > 1$, and a hump-like shape with a local maximum separates the two minima, yet the global minimum still happens at $t_{1} = 1$;
(3) when $\lambda _{2}^{\ast} < \lambda \leqslant \lambda _{3}^{\ast}$, the two local minima still coexists as in (2), yet the global minimum moves from $t_{1} = 1$ to $t_{2} > 1$ ;
(4) when $\lambda _{3}^{\ast} < \lambda$, the local minimum at $t_{1} = 1$ vanishes, and the local minimum at $t_{2} > 1$ becomes the global one.
For Fig.\ref{fig:result-ab} (a), we test $\lambda$ numerically and find that $\lambda _{1}^{\ast} \approx 0.771$, $\lambda _{2}^{\ast} \approx 0.830$, and $\lambda _{3}^{\ast} \approx 0.999$.
%an intuitive understanding
An intuitive understanding of the scenario is that, with an increasing $a/b$ and a grid cell becoming more elongated, there is a competition between a decrease in the intersection probability due to a larger width $a$ and an increase in the intersection probability due to a smaller height $b$ of the grid cell. Yet the dominant term in the competition simply depends on the control parameter $\lambda$.
%correspondence to critical transitions
The above scenario of moving minimum is much like the one in the analysis of critical transitions and hysteresis with a stability landscape in a multiple stable system (see Fig.(2.6) in \cite{Scheffer-2009}), while the two minima here correspond to the two stable equilibrium states of the system and the control parameter $\lambda$ acts as the strength of an external perturbation on the system.

%3D results
Results in Fig.\ref{fig:result-ab} (d)-(f) for $3$D cases follow a much more smooth way than those in the $2$D cases. In most cases in Fig.\ref{fig:result-ab} (d)-(f), $t = 1$ is the global minimum for intersection probability. Yet we can still numerically find a different global minimum at $t > 1$. For example, in Fig.\ref{fig:result-ab} (d) for the case of $3$D needles with $\lambda = 1.5$, we have $P(l, a, b) = 0.733559$ at $t = 1$ and a smaller $P(l, a, b) = 0.732816$ at $t = 1.605$.

%unsolved questions
For Fig.\ref{fig:result-ab}, there still remains some fundamental questions to resolve, such as to theoretically ascertain the scenario of moving minimum in $2$D cases, to analytically locate the global minimum with any given $\lambda$ in the $2$D cases, and to prove whether the results in the $3$D cases still follow the scenario in the $2$D cases. To fully resolve these questions, we need to carry out a detailed analysis of landscapes of intersection probabilities with respect to $t$. Due to the current amount of theoretical analysis and results in this paper, we would like to leave them to a future work.

%%%%%%%%%%%%%%%%%%%%%%%%%%%%%%%%
\section{Conclusion}
\label{sec:conclusion}

%conclusion
Here we analytically study BLNP in different versions on their intersection probabilities.
With a single probabilistic framework, we solve problem versions with dropped objects, with increasing analytical complexity, from a $2$D needle  with arbitrary length $l$ to a $3$D spherocylinder with arbitrary length $l$ and diameter $\sigma$ onto a grid with any width $a$ and height $b$.
Our analytical predictions are further validated by Monte Carlo simulation.

%connection to filtration proces
For the filtration process, our analytical results here only describe the contact probability between particles and a mesh upon the arrival of particles in flow. To fully construct a probabilistic picture of a filtration process \cite{Roussel.etal-PRL-2007}, a number of considerations should be combined into current  models.
First, multiple contacts between a particle with a mesh with the center of mass between contacts leads to a caught particle on a mesh, and an analytical theory for the distribution of number of intersections is highly relevant for an analysis of caking \cite{Gerber.etal-PRL-2018}.
Then, after a contact with a mesh, a particle slightly rotates due to a torque from the mesh and finally leads to a collective effect with other particles to induce clogging \cite{Dell.Franklin-JSTAT-2009}. Additional terms can be incorporated into our model to account the dynamical effect from the interaction between particles and a mesh.
Although highly mathematically involved, the above considerations represent essential steps towards a full analytical picture of filtration process. We leave the related treatment in future works.

%discussion
Besides, there are already many extensions of BLNP with physical backgrounds. Our analytical framework here, which applies to any shape size of needles and spherocylinders, helps to pave the way for principled theoretical approaches into untested parameter regions for these generalizations and variants.

%%%%%%%%%%%%%%%%%%%%%%%%%%%%%%%%
\section*{Acknowledgements}

This work is supported by
Guangdong Basic and Applied Basic Research Foundation (Grant No. 2022A1515011765),
Guangdong Major Project of Basic and Applied Basic Research No. 2020B0301030008, and
National Natural Science Foundation of China (Grant No. 12171479).

%%%%%%%%%%%%%%%%%%%%%%%%%%%%%%%%
\section{Appendix A: Calculation for Theory of BLNP with $2$D Needles}

%reset equation numbers in Appendix A
\setcounter{equation}{0}
\renewcommand{\theequation}{A.\arabic{equation}}

%relative size among l, a, b
%
\begin{figure}[htbp]
\begin{center}
 \includegraphics[width = 0.95 \linewidth]{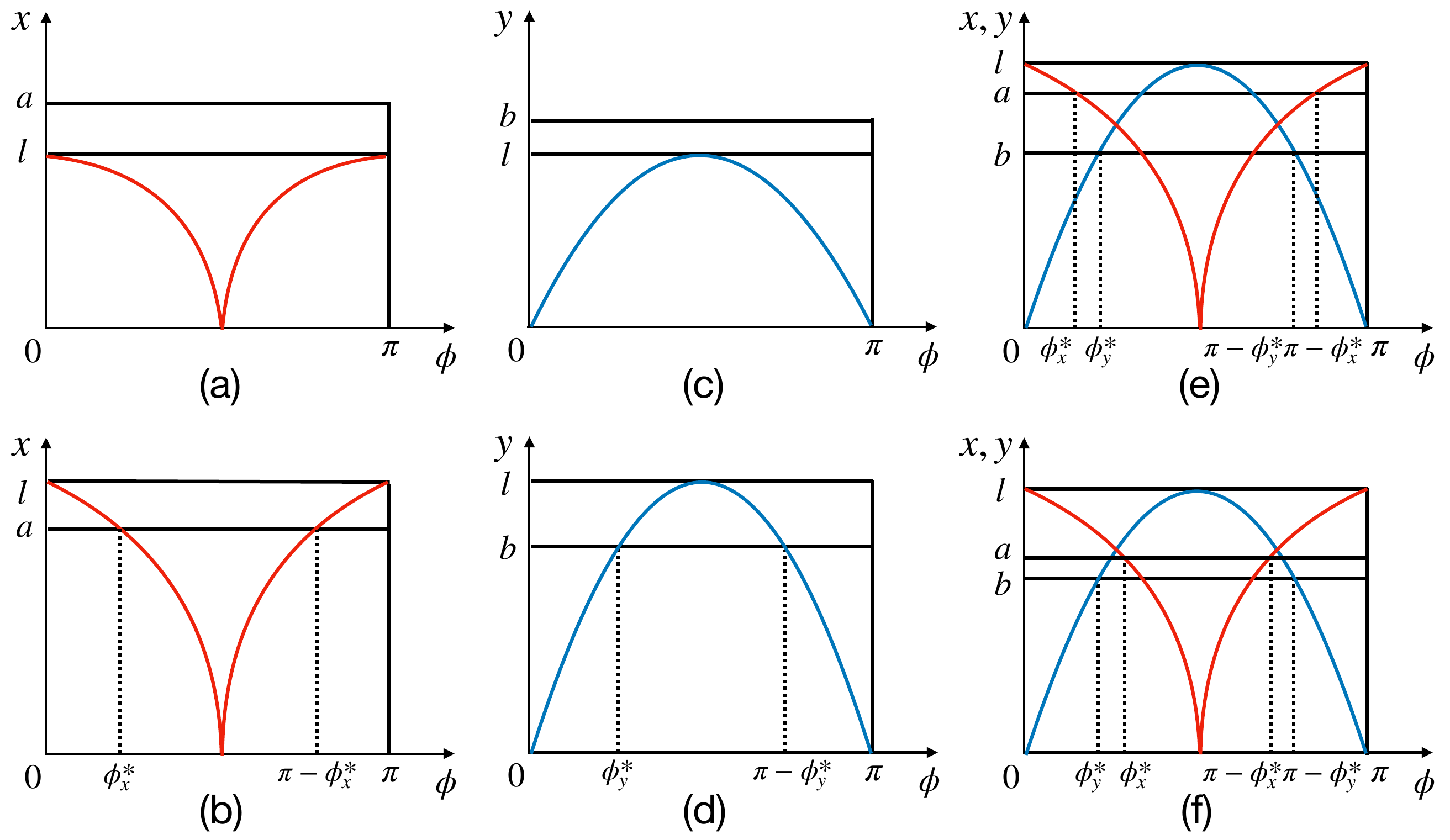}
\end{center}
\caption{
 \label{fig:diagram-2D}
Relative sizes among $l$, $a$, and $b$ for $A(l, a)$, $B(l, b)$, and $AB(l, a, b)$.
(a)-(b) The curve $x = |l \cos \phi|$ is shown under the conditions of $l \leqslant a$ and $a < l$, respectively.
$\phi _{x}^{\ast} \in (0, \pi / 2)$ and $\pi - \phi _{x}^{\ast}$ are the $\phi$-coordinates at which the lines $x = |l \cos \phi|$ and $x = a$ overlap.
(c)-(d) The curve $y = l \sin \phi$ is shown under the conditions of $l \leqslant b$ and $b < l$, respectively.
$\phi _{y}^{\ast} \in (0, \pi / 2)$ and $\pi - \phi _{y}^{\ast}$ are the $\phi$-coordinates at which the lines $y = l \sin \phi$ and $y = b$ overlap.
(e)-(f) The curves $x = |l \cos \phi|$ and $y = l \sin \phi$ are shown at the same time under the conditions of $\phi _{x}^{\ast} \leqslant \phi _{y}^{\ast}$, equivalently $l \leqslant L$, and $\phi _{y}^{\ast} < \phi _{x}^{\ast}$, equivalently $L < l$, respectively.}
\end{figure}
%

%basic logic to calculate integrals
We can find that the integrals in Eqs.(\ref{eq:Ala-2D})-(\ref{eq:ABlab-2D}) has a simple form as
\begin{equation}
\int _{0}^{t} \Theta (m - z) {\rm d}z = \min \{t, m\},
\end{equation}
in which $(t, m, z) = (a, |l \cos \phi|,x)$ and $(b, l \sin \phi, y)$.

%calculate A(l, a)
We write $A(l, a)$ in the case of $l \leqslant a$ and $a < l$ respectively as $A(l, a | l \leqslant a)$ and $A(l, a | a < l)$.
When $l \leqslant a$, $| l \cos \phi | \leqslant a$ always holds for $\phi \in [0, \pi)$.
Yet when $a < l$, $| l \cos \phi | \leqslant a$ holds only when $\phi \in [\phi _{x}^{\ast}, \pi - \phi _{x}^{\ast}]$, in which we define
\begin{equation}
\label{eq:phi-x-2D}
\phi _{x}^{\ast} = \arccos \frac {a}{l}.
\end{equation}
See Fig.\ref{fig:diagram-2D} (a) and (b) for an illustration.
We have
\begin{align}
A(l, a | l \leqslant a)
& = \int _{0}^{\frac {\pi}{2}} {\rm d} \phi \cdot l \cos \phi
+ \int _{\frac {\pi}{2}}^{\pi} {\rm d} \phi \cdot (- l \cos \phi) \nonumber \\
& = 2 l,\\
A(l, a | a < l)
& = \int _{0}^{\phi _{x}^{\ast}} {\rm d} \phi \cdot a
+ \int _{\phi _{x}^{\ast}}^{\frac {\pi}{2}} {\rm d} \phi \cdot l \cos \phi
+ \int _{\frac {\pi}{2}}^{\pi - \phi _{x}^{\ast}} {\rm d} \phi \cdot (- l \cos \phi)
+ \int _{\pi - \phi _{x}^{\ast}}^{\pi} {\rm d} \phi \cdot a \nonumber \\
& = 2l (1 - \sin \phi _{x}^{\ast}) + 2a \phi _{x}^{\ast}.
\end{align}
%

%calculate B(l, b)
We then write $B(l, b)$ in the case of $l \leqslant b$ and $b < l$ respectively as $B(l, b | l \leqslant b)$ and $B(l, b | b < l)$.
When $l \leqslant b$, $l \sin \phi \leqslant b$ always holds for $\phi \in [0, \pi)$.
Yet when $b < l$, $l \sin \phi \leqslant b$ holds only when $\phi \in [0, \phi _{y}^{\ast}] \cup [\pi - \phi _{y}^{\ast}, \pi)$, in which we define
\begin{equation}
\label{eq:phi-y-2D}
\phi _{y}^{\ast} = \arcsin \frac {b}{l}.
\end{equation}
See Fig.\ref{fig:diagram-2D} (c) and (d) for an illustration. We have
\begin{align}
B(l, b | l \leqslant b)
& = \int _{0}^{\pi} {\rm d}\phi \cdot l \sin \phi \nonumber \\
& = 2 l, \\
B(l, b | b < l)
& = \int _{0}^{\phi _{y}^{\ast}} {\rm d} \phi \cdot l \sin \phi
+ \int _{\phi _{y}^{\ast}}^{\pi - \phi _{y}^{\ast}} {\rm d} \phi \cdot b 
+ \int _{\pi - \phi _{y}^{\ast}}^{\pi} {\rm d} \phi \cdot l \sin \phi \nonumber \\
& = 2l \left(1 - \cos \phi _{y}^{\ast} \right) + 2b \left( \frac {\pi}{2} - \phi _{y}^{\ast} \right).
\end{align}
%

%calculate AB(l, a, b)
To calculate $AB(l, a, b)$ and finally $P(l, a, b)$, the relative sizes among $(l, a, b)$ are more involved. We will detail four cases below.

\subsection{The case of $l \leqslant b$}

We write $AB(l, a, b)$ when $l \leqslant b$ as $AB(l, a, b | l \leqslant b)$.
With Eq.(\ref{eq:ABlab-2D}) we have
\begin{align}
AB(l, a, b | l \leqslant b)
& = \int _{0}^{\frac {\pi}{2}} {\rm d} \phi \cdot l \cos \phi \cdot l \sin \phi 
 + \int _{\frac {\pi}{2}}^{\pi} {\rm d} \phi \cdot (- l \cos \phi) \cdot l \sin \phi \nonumber \\
& = l^2.
\end{align}
%
%Scoll
With Eq.(\ref{eq:Scoll-2D-reform}) we then have
\begin{align}
S_{\rm coll}
& = b A(l, a | l \leqslant a) + a B(l, b | l \leqslant b) - AB(l, a, b | l \leqslant b) \nonumber \\ 
& = 2 lb + 2la - l^2.
\end{align}
We thus have
\begin{equation}
\label{eq:Plab-2D-lb}
P(l, a, b)
= \frac {2}{\pi} \left( \frac {l}{a} + \frac {l}{b} \right) - \frac {1}{\pi} \frac {l^2}{ab}.
\end{equation}

\subsection{The case of $b < l \leqslant a$}

We write $AB(l, a, b)$ when $b < l \leqslant a$ as $AB(l, a, b | b < l \leqslant a)$.
We have
\begin{align}
AB(l, a, b | b < l \leqslant a)
& = \int _{0}^{\phi _{y}^{\ast}} {\rm d} \phi \cdot l \cos \phi \cdot l \sin \phi
+ \int _{\phi _{y}^{\ast}}^{\frac {\pi}{2}} {\rm d} \phi \cdot l \cos \phi \cdot b \nonumber \\
& + \int _{\frac {\pi}{2}}^{\pi - \phi _{y}^{\ast}} {\rm d} \phi \cdot (- l \cos \phi) \cdot b 
+ \int _{\pi - \phi _{y}^{\ast}}^{\pi} {\rm d} \phi \cdot (- l \cos \phi) \cdot l \sin \phi \nonumber \\
& = 2lb - b^2.
\end{align}
Then,
\begin{align}
S_{\rm coll}
& = b A(l, a | l \leqslant a) + a B(l, b | b < l) - AB(l, a, b | b < l \leqslant a) \nonumber \\ 
& = b^2 + 2 la (1 - \cos \phi _{y}^{\ast}) + 2 a b \left( \frac {\pi}{2} - \phi _{y}^{\ast} \right) \nonumber \\
& = b^2 + 2 la \left( 1 - \sqrt{1 - \frac {b^2}{l^2} } \right)
+ 2ab \left( \frac{\pi}{2} - \arcsin \frac {b}{l} \right).
\end{align}
We thus have
\begin{equation}
\label{eq:Plab-2D-bla}
P(l, a, b)
= \frac {1}{\pi} \frac {b}{a}
+ \frac {2}{\pi} \frac {l}{b} \left( 1 - \sqrt{1 - \frac {b^2}{l^2} } \right)
+  \frac {2}{\pi} \left( \frac {\pi}{2} - \arcsin \frac {b}{l} \right).
\end{equation}

\subsection{The cases in $a < l$}

To further calculate $AB(l, a, b)$, we should check the relative size between $\phi _{x}^{\ast}$ and $\phi _{y}^{\ast}$.
When $a < l \leqslant L$, we have $1 - a^2 / l^2 \leqslant b^2 / l^2$. Equivalently we have $\sin ^2 \phi _{x}^{\ast} \leqslant \sin ^2 \phi _{y}^{\ast}$, correspondingly $\phi _{x}^{\ast} \leqslant \phi _{y}^{\ast}$.
See Fig.\ref{fig:diagram-2D} (e) for an illustration.
Here we write $AB(l, a, b)$ when $a < l \leqslant L$ as $AB(l, a, b | a < l \leqslant L)$.
We have
\begin{align}
AB(l, a, b | a < l \leqslant L)
& = \int _{0}^{\phi _{x}^{\ast}} {\rm d} \phi \cdot a \cdot l \sin \phi
+ \int _{\phi _{x}^{\ast}}^{\phi _{y}^{\ast}} {\rm d} \phi \cdot l \cos \phi \cdot l \sin \phi 
+ \int _{\phi _{y}^{\ast}}^{\frac {\pi}{2}} {\rm d} \phi \cdot l \cos \phi \cdot b \nonumber \\
& + \int _{\frac {\pi}{2}}^{\pi - \phi _{y}^{\ast}} {\rm d} \phi \cdot (- l \cos \phi) \cdot b
+ \int _{\pi - \phi _{y}^{\ast}}^{\pi - \phi _{x}^{\ast}} {\rm d} \phi \cdot (- l \cos \phi) \cdot l \sin \phi
+ \int _{\pi - \phi _{x}^{\ast}}^{\pi} {\rm d} \phi \cdot a \cdot l \sin \phi \nonumber \\
& = 2la + 2lb - L^2 - l^2.
\end{align}
Then,
\begin{align}
S_{\rm coll}
& = b A(l, a | a < l) + a B(l, b | b < l) - AB(l, a, b | a < l \leqslant L) \nonumber \\ 
%result
& = L^2 + l^2
- 2 l b \sin \phi _{x}^{\ast} - 2 l a \cos \phi _{y}^{\ast} 
+ 2 a b \left( \phi _{x}^{\ast} + \frac {\pi}{2} - \phi _{y}^{\ast} \right) \nonumber \\
%reformulated
& = L^2 + l^2
- 2 \left( l b \sqrt {1 - \frac {a^2}{l^2} } + l a \sqrt {1 - \frac {b^2}{l^2} }  \right)
+ 2 ab \left( \pi - \arcsin \frac {a}{l} - \arcsin \frac {b}{l} \right).
\end{align}
Correspondingly,
\begin{equation}
\label{eq:Plab-2D-al1}
P(l, a, b)
= \frac {1}{\pi} \left( \frac {L^2}{ab} + \frac {l^2}{ab} \right)
- \frac {2}{\pi} \left( \frac {l}{a} \sqrt {1 - \frac {a^2}{l^2} } + \frac {l}{b} \sqrt {1 - \frac {b^2}{l^2} } \right)
+ \frac {2}{\pi} \left( \pi - \arcsin \frac {a}{l} - \arcsin \frac {b}{l} \right).
\end{equation}

When $L < l$, we have $\phi _{y}^{\ast} < \phi _{x}^{\ast}$ equivalently.
We write $AB(l, a, b)$ when $L < l$ as $AB(l, a, b | L < l)$.
We have
\begin{align}
AB(l, a, b | L < l)
& = \int _{0}^{\phi _{y}^{\ast}} {\rm d} \phi \cdot a \cdot l \sin \phi
+ \int _{\phi _{y}^{\ast}}^{\phi _{x}^{\ast}} {\rm d} \phi \cdot a \cdot b
+ \int _{\phi _{x}^{\ast}}^{\frac {\pi}{2}} {\rm d} \phi \cdot l \cos \phi \cdot b \nonumber \\
& + \int _{\frac {\pi}{2}}^{\pi - \phi _{x}^{\ast}} {\rm d} \phi \cdot (- l \cos \phi) \cdot b
+ \int _{\pi - \phi _{x}^{\ast}}^{\pi - \phi _{y}^{\ast}} {\rm d} \phi \cdot a \cdot b
+ \int _{\pi - \phi _{y}^{\ast}}^{\pi} {\rm d} \phi \cdot a \cdot l \sin \phi \nonumber \\
%result
& = 2 la (1 - \cos \phi _{y}^{\ast}) + 2 lb (1 - \sin \phi _{x}^{\ast})
+ 2ab (\phi _{x}^{\ast} - \phi _{y}^{\ast}).
\end{align}
Then,
\begin{align}
\label{eq:Scoll-2D-al2}
S_{\rm coll} 
& = b A(l, a | a < l) + a B(l, b | b < l) - AB(l, a, b | L < l) \nonumber \\ 
& = ab \pi,
\end{align}
after most terms cancel out. Correspondingly, we have
\begin{equation}
\label{eq:Plab-2D-al2}
P(l, a, b) = 1.
\end{equation}
%

%%%%%%%%%%%%%%%%%%%%%%%%%%%%%%%%
\section{Appendix B: Consistency of Equations in Theory of BLNP with $2$D Needles}

%reset equation numbers in Appendix B
\setcounter{equation}{0}
\renewcommand{\theequation}{B.\arabic{equation}}

%l = 0
When $l = 0$, in the case of $l \leqslant b$, we have
\begin{equation}
S_{\rm coll} = 0,
\end{equation}
which is quite intuitive.

%l = b
When $l = b$, both $S_{\rm coll}$ in the cases of $l \leqslant b$ and $b < l \leqslant a$ reduce to
\begin{equation}
S_{\rm coll} = b^2 + 2ab.
\end{equation}
%

% l = a
When $l = a$, both $S_{\rm coll}$ in the cases of $b < l \leqslant a$ and $a < l \leqslant L$ reduce to
\begin{equation}
S_{\rm coll}
= b^2 + 2 a^2 \left(1 - \sqrt{1 - \frac {b^2}{a^2}}  \right)
+ 2 ab \left( \frac{\pi}{2} - \arcsin \frac {b}{a} \right).
\end{equation}
%

%l = L
When $l = L$, $S_{\rm coll}$ in the case of $a < l \leqslant L$ reduces to
\begin{equation}
S_{\rm coll} = ab \pi,
\end{equation}
which equals $S_{\rm coll}$ in the case of $L < l$ in Eq.(\ref{eq:Scoll-2D-al2}).

%%%%%%%%%%%%%%%%%%%%%%%%%%%%%%%%
\section{Appendix C: Calculation for Theory of BLNP with $2$D Spherocylinders }

%reset equation numbers in Appendix C
\setcounter{equation}{0}
\renewcommand{\theequation}{C.\arabic{equation}}

As Eqs.(\ref{eq:phi-x-2D}) and (\ref{eq:phi-y-2D}), we first define here
\begin{align}
\hat{\phi} _{x}^{\ast}
& = \arccos \frac {a - \sigma}{l}, \\
\hat{\phi} _{y}^{\ast}
& = \arcsin \frac {b - \sigma}{l}.
\end{align}
%

%l < b
When $l \leqslant b - \sigma$, we have
\begin{align}
S_{\rm coll}
& = (b - \sigma) A(l, a - \sigma | l \leqslant a - \sigma)
+ (a - \sigma) B(l, b - \sigma | l \leqslant b - \sigma) \nonumber \\
& - AB(l, a - \sigma, b - \sigma | l \leqslant b - \sigma)
+ \left( \sigma a + \sigma b - \sigma ^2 \right) \pi \nonumber \\
%result
& = 2 l (b - \sigma) + 2 l (a - \sigma) - l^2 + \left( \sigma a + \sigma b - \sigma^2 \right) \pi.
\end{align}
Thus,
\begin{equation}
\label{eq:Plsab-2DSC-lb}
P(l, \sigma, a, b)
= \frac {2}{\pi} \frac {b - \sigma}{b} \frac {l}{a}
+ \frac {2}{\pi}  \frac {a - \sigma}{a} \frac {l}{b}
- \frac {1}{\pi} \frac {l^2}{ab}
+ \left( \frac {\sigma}{a} + \frac {\sigma}{b} - \frac {\sigma ^2}{ab} \right).
\end{equation}
%

% b < l < a
When $b - \sigma < l \leqslant a - \sigma$, we have
\begin{align}
S_{\rm coll}
& = (b - \sigma) A(l, a - \sigma | l \leqslant a - \sigma)
+ (a - \sigma) B(l, b - \sigma | b - \sigma < l) \nonumber \\
& - AB(l, a - \sigma, b - \sigma | b - \sigma < l \leqslant a - \sigma)
+ \left( \sigma a + \sigma b - \sigma ^2 \right) \pi \nonumber \\
%result
& = (b - \sigma)^2 + 2 l (a - \sigma)(1 - \cos \hat{\phi} _{y}^{\ast})
+ 2 (a - \sigma) (b - \sigma) \left( \frac {\pi}{2} - \hat{\phi} _{y}^{\ast} \right)
+ \left( \sigma a + \sigma b - \sigma ^2 \right) \pi \nonumber \\
%reform
& = (b - \sigma)^2 + 2l (a - \sigma) \left[ 1 - \sqrt {1 - \frac {(b - \sigma)^2}{l^2}} \right] 
+ 2 (a - \sigma) (b - \sigma) \left( \frac {\pi}{2} - \arcsin \frac {b - \sigma}{l} \right)
+ \left( \sigma a + \sigma b - \sigma ^2 \right) \pi.
\end{align}
Thus,
\begin{align}
\label{eq:Plsab-2DSC-bla}
P(l, \sigma, a, b)
& = \frac {1}{\pi} \frac {(b - \sigma)^2}{ab}
+ \frac {2}{\pi} \frac {a - \sigma}{a} \frac {l}{b}
\left[ 1 - \sqrt{1 - \frac {(b - \sigma)^2}{l^2}} \right] \nonumber \\
& + \frac {2}{\pi} \frac {a - \sigma}{a} \frac {b - \sigma}{b}
\left( \frac {\pi}{2} - \arcsin \frac {b - \sigma}{l} \right)
+ \left( \frac {\sigma}{a} + \frac {\sigma}{b} - \frac {\sigma ^2}{ab} \right).
\end{align}
%

%a < l < L
When $a - \sigma < l \leqslant \hat {L}$,
we have $\hat{\phi} _{x}^{\ast} \leqslant \hat{\phi} _{y}^{\ast}$ correspondingly. Then
\begin{align}
S_{\rm coll}
& = (b - \sigma) A(l, a - \sigma | a - \sigma < l)
+ (a - \sigma) B(l, b - \sigma | b - \sigma < l) \nonumber \\
& - AB(l, a - \sigma, b - \sigma | a - \sigma < l \leqslant \hat{L})
+ \left( \sigma a + \sigma b - \sigma ^2 \right) \pi \nonumber \\
%result
& = \hat{L}^2 + l^{2}
- 2 l (b - \sigma) \sin \hat{\phi} _{x}^{\ast} - 2 l (a - \sigma) \cos \hat{\phi} _{y}^{\ast} \nonumber \\
& + 2 (a - \sigma) (b - \sigma) \left( \hat{\phi} _{x}^{\ast} + \frac {\pi}{2} - \hat{\phi} _{y}^{\ast} \right)
+ \left( \sigma a + \sigma b - \sigma ^2 \right) \pi \nonumber \\
%reform
& = \hat{L}^2 + l^{2}
- 2 l (b - \sigma) \sqrt {1- \frac {(a - \sigma)^2}{l^2}}
- 2 l (a - \sigma) \sqrt {1- \frac {(b - \sigma)^2}{l^2}} \nonumber \\
& + 2 (a - \sigma) (b - \sigma)
\left( \pi - \arcsin \frac {a - \sigma}{l} - \arcsin \frac {b - \sigma}{l} \right)
+ \left( \sigma a + \sigma b - \sigma ^2 \right) \pi.
\end{align}
Thus,
\begin{align}
\label{eq:Plsab-2DSC-al1}
P(l, \sigma, a, b)
& = \frac {1}{\pi} \left( \frac{\hat{L}^2}{ab}
+ \frac{l^{2}}{ab} \right)
- \frac {2}{\pi} \left[ \frac {b - \sigma}{b} \frac {l}{a} \sqrt {1 - \frac {(a - \sigma)^2}{l^2}}
+ \frac {a - \sigma}{a} \frac {l}{b} \sqrt {1 - \frac {(b - \sigma)^2}{l^2}} \right] \nonumber \\
& + \frac {2}{\pi} \frac {a - \sigma}{a} \frac {b - \sigma}{b} \left( \pi - \arcsin \frac {a - \sigma}{l} - \arcsin \frac {b - \sigma}{l} \right)
+ \left( \frac {\sigma}{a} + \frac {\sigma}{b} - \frac {\sigma^2}{ab} \right).
\end{align}
%

% L < l
When $\hat{L} < l$, we have $\hat{\phi} _{y}^{\ast} < \hat{\phi} _{x}^{\ast}$.
Then,
\begin{align}
\label{eq:Scoll-2DSC-al2}
S_{\rm coll}
& = (b - \sigma) A(l, a - \sigma | a - \sigma < l)
+ (a - \sigma) B(l, b - \sigma | b - \sigma < l) \nonumber \\
& - AB(l, a - \sigma, b - \sigma | \hat{L} < l)
+ \left( \sigma a + \sigma b - \sigma ^2 \right) \pi \nonumber \\
& = a b \pi.
\end{align}
Thus,
\begin{equation}
\label{eq:Plsab-2DSC-al2}
P(l, \sigma, a, b) = 1.
\end{equation}
%

%%%%%%%%%%%%%%%%%%%%%%%%%%%%%%%%
\section{Appendix D: Consistency of Equations in Theory of BLNP with $2$D Spherocylinders}

%reset equation numbers in Appendix D
\setcounter{equation}{0}
\renewcommand{\theequation}{D.\arabic{equation}}

%l -->0 
When $l = 0$, in the case of $l \leqslant b - \sigma$ we have
\begin{equation}
\label{eq:Scoll-2DSC-l0}
S_{\rm coll} = \left( \sigma a + \sigma b - \sigma ^2 \right) \pi.
\end{equation}
Intuitively, a $2$D spherocylinder with $l = 0$ reduces to a circle. An intersection between the circle and a grid cell happens when the center of the circle is within a distance of its radius $\sigma / 2$ from the four boundaries of the grid cell. We easily have
\begin{equation}
S_{\rm coll}
=  [a b - (a - \sigma) (b - \sigma) ] \times \pi
= \left( \sigma a + \sigma b - \sigma ^2 \right) \pi.
\end{equation}
%

%l = b
When $l = b - \sigma$, both $S_{\rm coll}$ in the cases of $l \leqslant b - \sigma$ and $b - \sigma < l \leqslant a - \sigma$ reduce to 
\begin{equation}
S_{\rm coll}
= (b - \sigma)^{2} + 2 (a - \sigma) (b - \sigma) + \left( \sigma a + \sigma b - \sigma ^2 \right) \pi.
\end{equation}
%

%l = a
When $l = a - \sigma$, both $S_{\rm coll}$ in the cases of $b - \sigma < l \leqslant a - \sigma$ and $a - \sigma < l \leqslant \hat{L}$ correspond to
\begin{align}
S_{\rm coll}
& = (b - \sigma)^{2}
+ 2 (a - \sigma)^{2} \left[ 1 - \sqrt {1 - \frac {(b - \sigma)^2}{(a - \sigma)^2} } \right] \nonumber \\
& + 2 (a - \sigma) (b - \sigma) \left( \frac {\pi}{2} - \arcsin \frac {b - \sigma}{a -  \sigma} \right)
+ \left( \sigma a + \sigma b - \sigma ^2 \right) \pi.
\end{align}
%

% l = sqrt
When $l = \hat{L}$ in the case of $a - \sigma < l \leqslant \hat{L}$, we have $\cos ^2 \hat{\phi} _{x}^{\ast} + \sin ^2 \hat{\phi} _{y}^{\ast} = 1$, equivalently $\hat{\phi} _{x}^{\ast} = \hat{\phi} _{y}^{\ast}$. After some calculation we have
\begin{equation}
S_{\rm coll} = a b \pi,
\end{equation}
which coincides with the uniform $S_{\rm coll}$ in the case of $\hat{L} < l$ in Eq.(\ref{eq:Scoll-2DSC-al2}).

%%%%%%%%%%%%%%%%%%%%%%%%%%%%%%%%
\section{Appendix E: Calculation for Theory of BLNP with $3$D Needles}

%reset equation numbers in Appendix E
\setcounter{equation}{0}
\renewcommand{\theequation}{E.\arabic{equation}}

% figure: calculation in 3D
%
\begin{figure}[htbp]
\begin{center}
 \includegraphics[width = 0.95 \linewidth]{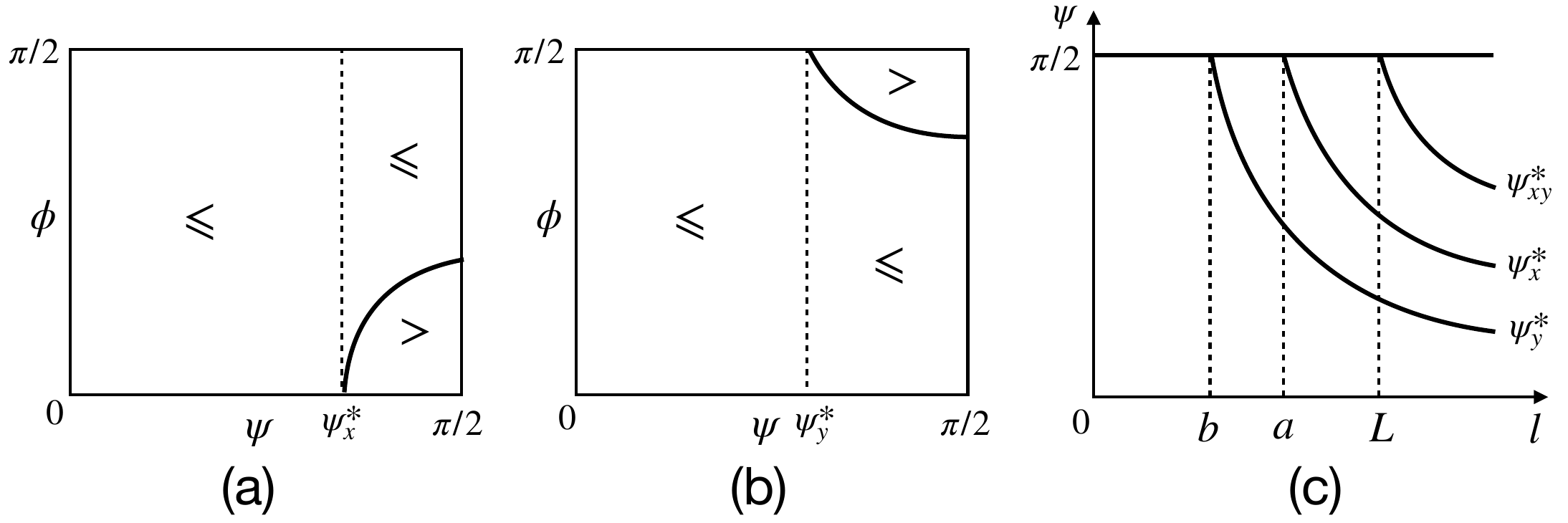}
\end{center}
\caption{
 \label{fig:diagram-3D}
Schematic diagrams of calculation in the BLNP with $3$D needles.
(a) The dashed and solid lines correspond to $\psi = \psi _{x}^{\ast}$ and $\phi = \phi _{x}^{\ast} (\psi)$, respectively.
The notions $\leqslant$ and $>$ within diagrams respectively denotes the regions of parameters $(\psi, \phi)$ in which $l \sin \psi \cos \phi \leqslant a$ and $l \sin \psi \cos \phi > a$ hold.
(b) The dashed and solid lines show $\psi = \psi _{y}^{\ast}$ and $\phi = \phi _{y}^{\ast} (\psi)$, respectively.
The notions $\leqslant$ and $>$ respectively denote the regions of parameters $(\psi, \phi)$ in which $l \sin \psi \sin \phi \leqslant b$ and $l \sin \psi \sin \phi > b$ hold.
(c) Boundaries of regions $\psi _{x}^{\ast}$, $\psi _{y}^{\ast}$, and $\psi _{xy}^{\ast}$, when Eqs.(\ref{eq:ABlab-3D}) is calculated by integrating on $\psi$, are shown in the cases of $l \leqslant b$, $b < l \leqslant a$, $a < l \leqslant L$, and $L < l$, respectively.}
\end{figure}
%

%calculate A(l, a)
We first calculate $A(l, a)$ and $B(l, b)$.
We write $A^{\rm 3D}(l, a)$ when $l \leqslant a$ and $a < l$ as $A^{\rm 3D}(l, a | l \leqslant a)$ and $A^{\rm 3D}(l, a | a < l)$, respectively. 
% l <= a
When $l \leqslant a$, $|l \sin \psi \cos \phi | \leqslant a$ always holds for any $\psi \in [0, \pi/2]$ and $\phi \in [0, \pi)$. We have
\begin{align}
A^{\rm 3D}(l, a | l \leqslant a)
& = \int _{0}^{\frac {\pi}{2}} {\rm d} \psi 
\int _{0}^{\frac {\pi}{2}} {\rm d} \phi \cdot
l \sin \psi \cos \phi 
+ \int _{0}^{\frac {\pi}{2}} {\rm d} \psi 
\int _{\frac {\pi}{2}}^{\pi} {\rm d} \phi \cdot
(- l \sin \psi \cos \phi) \nonumber \\
& = 2l,
\end{align}
%
%a < l
When $a < l$, 
 the picture is slightly more complicated. We first define
\begin{align}
\label{eq:psi-x-3D}
\psi _{x}^{\ast}
& = \arcsin \frac{a}{l}, \\
\phi _{x}^{\ast} (\psi)
& = \arccos \frac{a}{l \sin \psi}.
\end{align}
We can see that, when $\psi \leqslant \psi _{x}^{\ast}$, we have $l \sin \psi \cos \phi \leqslant a$ satisfied for any $\phi \in [0, \pi / 2]$.
Yet when $\psi > \psi _{x}^{\ast}$, two situations happen depending on the value of $\phi$:
when $\phi \in [\phi _{x}^{\ast}(\psi), \pi / 2]$,  $l \sin \psi \cos \phi \leqslant a$ still holds;
yet when $\phi \in [0, \phi _{x}^{\ast}(\psi))$, $l \sin \psi \cos \phi > a$ holds instead.
A simple schematics is in Fig.\ref{fig:diagram-3D} (a).
Later we simply write $\phi _{x}^{\ast} (\psi)$ as $\phi _{x}^{\ast}$ for short in this section.
%result
We have
\begin{align}
A^{\rm 3D}(l, a | a < l)
& = \int _{0}^{\psi _{x} ^{\ast}} {\rm d} \psi
\int _{0}^{\frac {\pi}{2}} {\rm d} \phi \cdot
l \sin \psi \cos \phi
+ \int _{0}^{\psi _{x} ^{\ast}} {\rm d} \psi
\int _{\frac {\pi}{2}}^{\pi} {\rm d} \phi \cdot
(- l \sin \psi \cos \phi) \nonumber \\
& + \int _{\psi _{x} ^{\ast}}^{\frac{\pi}{2} } {\rm d} \psi
\int _{0}^{\phi _{x} ^{\ast}} {\rm d} \phi \cdot a
+ \int _{\psi _{x} ^{\ast}}^{\frac {\pi}{2}} {\rm d} \psi
\int _{\phi _{x} ^{\ast}}^{\frac {\pi}{2}} {\rm d} \phi \cdot
l \sin \psi \cos \phi \nonumber \\
& + \int _{\psi _{x} ^{\ast}}^{\frac{\pi}{2} } {\rm d} \psi
\int _{\frac {\pi}{2}}^{\pi - \phi _{x} ^{\ast}} {\rm d} \phi \cdot
(- l \sin \psi \cos \phi) 
 + \int _{\psi _{x} ^{\ast}}^{\frac {\pi}{2}} {\rm d} \psi
\int _{\pi - \phi _{x} ^{\ast}}^{\pi} {\rm d} \phi \cdot a \nonumber \\
%result
& = 2l - 2l \int _{\psi _{x}^{\ast}}^{\frac {\pi}{2}} {\rm d} \psi \sqrt{\sin ^{2} \psi - \frac {a^2}{l^2}}
+ 2 a \int _{\psi _{x}^{\ast}}^{\frac {\pi}{2}} {\rm d} \psi \left(\frac {\pi}{2} - \arcsin \frac {a}{l \sin \psi} \right) \nonumber \\
%in short-handed form
& = 2l  - 2l F \left( \psi _{x}^{\ast}, \frac {\pi}{2}, \frac {a}{l} \right)
+ 2a G \left( \psi _{x}^{\ast}, \frac {\pi}{2}, \frac {a}{l} \right),
\end{align}
in whose last step we define two integrals respectively as
\begin{align}
\label{eq:Fabt}
F(a, b, t) 
& \equiv \int _{a}^{b} {\rm d} \psi \sqrt {\sin^{2} \psi - t^2}, \\
\label{eq:Gabt}
G(a, b, t)
& \equiv \int _{a}^{b} {\rm d} \psi \left( \frac {\pi}{2} - \arcsin \frac {t}{\sin \psi} \right).
\end{align}
These two integrals can be calculated with standard numerical integration methods \cite{Press.etal-2007-3e}.

%l <= b
We then write $B^{\rm 3D}(l, b)$ when $l \leqslant b$ and $b < l$ as $B^{\rm 3D}(l, b | l \leqslant b)$ and $B^{\rm 3D}(l, b | b < l)$, respectively.
When $l \leqslant b$, $l \sin \psi \sin \phi \leqslant b$ always hold for any $\psi \in [0, \pi/2]$ and $\phi \in [0, \pi)$. 
We have
\begin{align}
B^{\rm 3D}(l, b | l \leqslant b)
& = \int _{0}^{\frac {\pi}{2}} {\rm d} \psi 
\int _{0}^{\pi} {\rm d} \phi \cdot
l \sin \psi \sin \phi \nonumber \\
& = 2l.
\end{align}
%
%b < l
When $b < l$, as in the case of $a < l$ for $A^{\rm 3D}(l, a)$, we first define
\begin{align}
\label{eq:psi-y-3D}
\psi _{y}^{\ast}
& = \arcsin \frac{b}{l}, \\
\phi _{y}^{\ast} (\psi)
& = \arcsin \frac{b}{l \sin \psi}.
\end{align}
We can find that, when $\psi \leqslant \psi _{y}^{\ast}$, $l \sin \psi \sin \phi \leqslant b$ is satisfied for any $\phi \in [0, \pi / 2]$.
Yet when $\psi > \psi _{y}^{\ast}$, two situations happen depending on the value of $\phi$:
when $\phi \in [0, \phi _{y}^{\ast}(\psi)]$,  $l \sin \psi \sin \phi \leqslant b$ still holds;
yet when $\phi \in (\phi _{y}^{\ast}(\psi), \pi / 2]$, $l \sin \psi \sin \phi > b$ holds instead.
A simple schematics is in Fig.\ref{fig:diagram-3D} (b).
Later we simply write $\phi _{y}^{\ast} (\psi)$ as $\phi _{y}^{\ast}$ for short in this section.
We  have
\begin{align}
B^{\rm 3D}(l, b | b < l)
& = \int _{0}^{\psi _{y} ^{\ast}} {\rm d} \psi
\int _{0}^{\pi} {\rm d} \phi \cdot l \sin \psi \sin \phi \nonumber \\
& + \int _{\psi _{y} ^{\ast}}^{\frac {\pi}{2}} {\rm d} \psi
\int _{0}^{\phi _{y}^{\ast}} {\rm d} \phi \cdot l \sin \psi \sin \phi
+ \int _{\psi _{y} ^{\ast}}^{\frac {\pi}{2}} {\rm d} \psi
\int _{\phi _{y}^{\ast}}^{\pi - \phi _{y}^{\ast}} {\rm d} \phi \cdot b
+ \int _{\psi _{y} ^{\ast}}^{\frac {\pi}{2}} {\rm d} \psi
\int _{\pi - \phi _{y}^{\ast}}^{\pi} {\rm d} \phi \cdot l \sin \psi \sin \phi \nonumber \\
%result
& = 2l - 2l \int _{\psi _{y}^{\ast}}^{\frac {\pi}{2}} {\rm d} \psi \sqrt{\sin ^{2} \psi - \frac {b^2}{l^2}}
+ 2 b \int _{\psi _{y}^{\ast}}^{\frac {\pi}{2}} {\rm d} \psi \left(\frac {\pi}{2} - \arcsin \frac {b}{l \sin \psi} \right) \nonumber \\
%in short-handed form
& = 2l  - 2l F \left( \psi _{y}^{\ast}, \frac {\pi}{2}, \frac {b}{l} \right)
+ 2b G \left( \psi _{y}^{\ast}, \frac {\pi}{2}, \frac {b}{l} \right).
\end{align}

To further calculate $AB^{\rm 3D}(l, a, b)$ and finally $P^{\rm 3D}(l, a, b)$, we consider the following four cases.

\subsection{The case of $l \leqslant b$}

Here the range $\psi \in [0, \pi / 2]$ represents a simple case in $AB(l, a, b)$, in which $|l \sin \psi \cos \phi | \leqslant a$ and $l \sin \psi \sin \phi \leqslant b$ both hold for any $\phi \in [0, \pi)$. See Fig.\ref{fig:diagram-3D}(c) for an illustration.
We write $AB^{\rm 3D}(l, a, b)$ when $l \leqslant b$ as $AB^{\rm 3D}(l, a, b| l \leqslant b)$.
We have
\begin{align}
& AB^{\rm 3D}(l, a, b | l \leqslant b) \nonumber \\
& = \int _{0}^{\frac {\pi}{2}} {\rm d} \psi
\int _{0}^{\frac {\pi}{2}} {\rm d} \phi
\cdot l \sin \psi \cos \phi \cdot l \sin \psi \sin \phi
+ \int _{0}^{\frac {\pi}{2}} {\rm d} \psi
\int _{\frac {\pi}{2}}^{\pi} {\rm d} \phi
\cdot (- l \sin \psi \cos \phi) \cdot l \sin \psi \sin \phi \nonumber \\
& = S^{(1)} \left( 0, \frac {\pi}{2} \right).
\end{align}

In the above equations, we define and calculate a summation as
\begin{align}
S^{(1)}(0, \psi ^{\ast})
& \equiv \int _{0}^{\psi ^{\ast}} {\rm d} \psi
\int _{0}^{\frac {\pi}{2}} {\rm d} \phi \cdot l \sin \psi \cos \phi \cdot l \sin \psi \sin \phi
+ \int _{0}^{\psi ^{\ast}} {\rm d} \psi
\int _{\frac {\pi}{2}}^{\pi} {\rm d} \phi \cdot (- l \sin \psi \cos \phi) \cdot l \sin \psi \sin \phi \nonumber \\
%result
& = \frac {l^2}{2} \left( \psi ^{\ast} - \sin \psi ^{\ast} \cos \psi^{\ast} \right).
\end{align}
Thus we have
\begin{align}
S^{(1)} \left( 0, \frac {\pi}{2} \right)
& = \frac {\pi}{4} l^2, \\
S^{(1)}(0, \psi _{y}^{\ast})
& = \frac {l^2}{2} \psi _{y}^{\ast} - \frac {lb}{2} \cos \psi _{y}^{\ast}.
\end{align}

With Eq.(\ref{eq:Scoll-3D-reform}) we have
\begin{align}
S_{\rm coll}^{\rm 3D}
& = b A^{\rm 3D}(l, a | l \leqslant a) + a B^{\rm 3D}(l, b | l \leqslant b) - AB^{\rm 3D}(l, a, b | l \leqslant b) \nonumber \\
& = 2lb + 2la - \frac {\pi}{4} l^2.
\end{align}
Then,
\begin{equation}
\label{eq:Plab-3D-lb}
P^{\rm 3D}(l, a, b) = \frac {4}{\pi ^2} \left( \frac {l}{a} + \frac {l}{b} \right) - \frac {1}{2 \pi} \frac {l^2}{ab}.
\end{equation}

\subsection{The case of $b < l \leqslant a$}

Here we consider the integral on $\psi$ in $AB(l, a, b)$ in two ranges of $[0, \psi _{y}^{\ast}]$ and $(\psi _{y}^{\ast}, \pi / 2]$ in a separate way.
See Fig.\ref{fig:diagram-3D} (c) for a schematics.
We write $AB^{\rm 3D}(l, a, b)$ when $b < l \leqslant a$ as $AB^{\rm 3D}(l, a, b| b < l \leqslant a)$.
We have
\begin{align}
& AB^{\rm 3D}(l, a, b | b < l \leqslant a) \nonumber \\
%S^{(1)}
& = \int _{0}^{\psi _{y}^{\ast}} {\rm d} \psi
\int _{0}^{\frac {\pi}{2}} {\rm d} \phi
\cdot l \sin \psi \cos \phi \cdot l \sin \psi \sin \phi
+ \int _{0}^{\psi _{y}^{\ast}} {\rm d} \psi
\int _{\frac {\pi}{2}}^{\pi} {\rm d} \phi
\cdot (- l \sin \psi \cos \phi) \cdot l \sin \psi \sin \phi \nonumber \\
%S^{(2)}
& + \int _{\psi _{y}^{\ast}}^{\frac {\pi}{2}} {\rm d} \psi
\int _{0}^{\phi _{y}^{\ast}} {\rm d} \phi
\cdot l \sin \psi \cos \phi \cdot l \sin \psi \sin \phi
+ \int _{\psi _{y}^{\ast}}^{\frac {\pi}{2}} {\rm d} \psi
\int _{\phi _{y}^{\ast}}^{\frac {\pi}{2}} {\rm d} \phi
\cdot l \sin \psi \cos \phi \cdot b \nonumber \\
& + \int _{\psi _{y}^{\ast}}^{\frac {\pi}{2}} {\rm d} \psi
\int _{\frac {\pi}{2}}^{\pi - \phi _{y}^{\ast}} {\rm d} \phi
\cdot (- l \sin \psi \cos \phi) \cdot b
+ \int _{\psi _{y}^{\ast}}^{\frac {\pi}{2}} {\rm d} \psi
\int _{\pi - \phi _{y}^{\ast}}^{\pi} {\rm d} \phi
\cdot (- l \sin \psi \cos \phi) \cdot l \sin \psi \sin \phi  \nonumber \\
%sums
& = S^{(1)}(0, \psi _{y}^{\ast}) + S^{(2)} \left( \psi _{y}^{\ast}, \frac {\pi}{2} \right).
\end{align}

In the above equations, we define and calculate a summation as
\begin{align}
S^{(2)}(\psi _{y}^{\ast}, \psi ^{\ast})
& \equiv \int _{\psi _{y}^{\ast}}^{\psi ^{\ast}} {\rm d} \psi
\int _{0}^{\phi _{y}^{\ast}} {\rm d} \phi \cdot
l \sin \psi \cos \phi \cdot l \sin \psi \sin \phi
+ \int _{\psi _{y}^{\ast}}^{\psi ^{\ast}} {\rm d} \psi
\int _{\phi _{y}^{\ast}}^{\frac{\pi}{2}} {\rm d} \phi \cdot 
l \sin \psi \cos \phi \cdot b \nonumber \\
& + \int _{\psi _{y}^{\ast}}^{\psi ^{\ast}} {\rm d} \psi
\int _{\frac {\pi}{2}}^{\pi - \phi _{y}^{\ast}} {\rm d} \phi \cdot
(- l \sin \psi \cos \phi) \cdot b
+ \int _{\psi _{y}^{\ast}}^{\psi ^{\ast}} {\rm d} \psi
\int _{\pi - \phi _{y}^{\ast}}^{\pi} {\rm d} \phi \cdot
(- l \sin \psi \cos \phi) \cdot l \sin \psi \sin \phi \nonumber \\
%result
& = 2lb(\cos \psi _{y}^{\ast} - \cos \psi ^{\ast}) - b^2 (\psi ^{\ast} - \psi _{y}^{\ast}).
\end{align}
Thus we have
\begin{align}
S^{(2)} \left(\psi _{y}^{\ast}, \frac {\pi}{2} \right)
& = 2lb\cos \psi _{y}^{\ast} - b^2 \left( \frac {\pi}{2} - \psi _{y}^{\ast} \right), \\
S^{(2)} \left(\psi _{y}^{\ast}, \psi _{x}^{\ast} \right)
& = 2lb(\cos \psi _{y}^{\ast} - \cos \psi _{x}^{\ast}) - b^2 (\psi _{x}^{\ast} - \psi _{y}^{\ast}).
\end{align}
Then,
\begin{equation}
AB^{\rm 3D}(l, a, b | b < l \leqslant a)
= \frac {l^2}{2} \psi _{y}^{\ast}
+ \frac {3}{2} l b \cos \psi _{y}^{\ast}
- b^2 \left( \frac {\pi}{2} - \psi _{y}^{\ast} \right).
\end{equation}

We have
\begin{align}
S_{\rm coll}^{\rm 3D}
& = b A^{\rm 3D}(l, a | l \leqslant a) + a B^{\rm 3D}(l, b | b < l) - AB^{\rm 3D}(l, a, b | b < l \leqslant a) \nonumber \\
%result
& = 2lb + \left[ 2la - 2la F \left(\psi _{y}^{\ast}, \frac {\pi}{2}, \frac {b}{l} \right)
+ 2abG\left( \psi _{y}^{\ast}, \frac {\pi}{2}, \frac {b}{l} \right) \right] \nonumber \\
& - \left[ \frac {l^2}{2} \psi _{y}^{\ast}
+ \frac {3}{2} lb \cos \psi _{y}^{\ast} - b^2 \left( \frac {\pi}{2} - \psi _{y}^{\ast} \right) \right] \nonumber \\
%reformed
& = 2lb + \left[ 2la - 2la F \left( \arcsin \frac {b}{l}, \frac {\pi}{2}, \frac {b}{l} \right)
+ 2abG\left( \arcsin \frac {b}{l}, \frac {\pi}{2}, \frac {b}{l} \right) \right] \nonumber \\
& - \left[ \frac {l^2}{2} \arcsin \frac {b}{l}
+ \frac {3}{2} lb \sqrt {1 - \frac {b^2}{l^2}} - b^2 \left( \frac {\pi}{2} - \arcsin \frac {b}{l} \right) \right].
\end{align}
Then,
\begin{align}
\label{eq:Plab-3D-bla}
P^{\rm 3D}(l, a, b)
& = \frac {4}{\pi ^2} \frac {l}{a}
+ \left[ \frac {4}{\pi ^2} \frac {l}{b}
- \frac {4}{\pi ^2} \frac {l}{b} F \left( \arcsin \frac {b}{l}, \frac {\pi}{2}, \frac {b}{l} \right)
+ \frac {4}{\pi ^2} G\left( \arcsin \frac {b}{l}, \frac {\pi}{2}, \frac {b}{l} \right) \right] \nonumber \\
& - \left[ \frac {1}{\pi ^2} \frac {l^2}{ab} \arcsin \frac {b}{l}
+ \frac {3}{\pi ^2} \frac {l}{a} \sqrt {1 - \frac {b^2}{l^2}}
- \frac {2}{\pi ^2} \frac {b}{a} \left( \frac {\pi}{2} - \arcsin \frac {b}{l} \right) \right].
\end{align}

\subsection{The cases in $a < l$}

We consider here two subcases as $a < l \leqslant L$ and $L < l$.
%a < l < L
When $a < l \leqslant L$, we calculate the integral on $\psi$ in $AB(l, a, b)$ in three ranges of $[0, \psi _{y}^{\ast}]$, $(\psi _{y}^{\ast}, \psi _{x}^{\ast}]$, and $(\psi _{x}^{\ast}, \pi / 2]$ in a separate way.
See Fig.\ref{fig:diagram-3D} (c) for an illustration.
We write $AB^{\rm 3D}(l, a, b)$ when $a < l \leqslant L$ as $AB^{\rm 3D}(l, a, b| a < l \leqslant L)$.
We have
\begin{align}
& AB^{\rm 3D}(l, a, b | a < l \leqslant L) \nonumber \\
%S^{(1)}
& = \int _{0}^{\psi _{y}^{\ast}} {\rm d} \psi
\int _{0}^{\frac {\pi}{2}} {\rm d} \phi
\cdot l \sin \psi \cos \phi \cdot l \sin \psi \sin \phi
+ \int _{0}^{\psi _{y}^{\ast}} {\rm d} \psi
\int _{\frac {\pi}{2}}^{\pi} {\rm d} \phi
\cdot (- l \sin \psi \cos \phi) \cdot l \sin \psi \sin \phi \nonumber \\
%S^{(2)}
& + \int _{\psi _{y}^{\ast}}^{\psi _{x}^{\ast}} {\rm d} \psi
\int _{0}^{\phi _{y}^{\ast}} {\rm d} \phi
\cdot l \sin \psi \cos \phi \cdot l \sin \psi \sin \phi
+ \int _{\psi _{y}^{\ast}}^{\psi _{x}^{\ast}} {\rm d} \psi
\int _{\phi _{y}^{\ast}}^{\frac {\pi}{2}} {\rm d} \phi
\cdot l \sin \psi \cos \phi \cdot b \nonumber \\
& + \int _{\psi _{y}^{\ast}}^{\psi _{x}^{\ast}} {\rm d} \psi
\int _{\frac {\pi}{2}}^{\pi - \phi _{y}^{\ast}} {\rm d} \phi
\cdot (- l \sin \psi \cos \phi) \cdot b
+ \int _{\psi _{y}^{\ast}}^{\psi _{x}^{\ast}} {\rm d} \psi
\int _{\pi - \phi _{y}^{\ast}}^{\pi} {\rm d} \phi
\cdot (- l \sin \psi \cos \phi) \cdot l \sin \psi \sin \phi  \nonumber \\
%S^{(3)}
& + \int _{\psi _{x}^{\ast}}^{\frac {\pi}{2}} {\rm d} \psi
\int _{0}^{\phi _{x}^{\ast}} {\rm d} \phi
\cdot a \cdot l \sin \psi \sin \phi
+ \int _{\psi _{x}^{\ast}}^{\frac {\pi}{2}} {\rm d} \psi
\int _{\phi _{x}^{\ast}}^{\phi _{y}^{\ast}} {\rm d} \phi
\cdot l \sin \psi \cos \phi \cdot l \sin \psi \sin \phi \nonumber \\
& + \int _{\psi _{x}^{\ast}}^{\frac {\pi}{2}} {\rm d} \psi
\int _{\phi _{y}^{\ast}}^{\frac {\pi}{2}} {\rm d} \phi
\cdot l \sin \psi \cos \phi \cdot b \nonumber \\
& + \int _{\psi _{x}^{\ast}}^{\frac {\pi}{2}} {\rm d} \psi
\int _{\frac {\pi}{2}}^{\pi - \phi _{y}^{\ast}} {\rm d} \phi
\cdot (- l \sin \psi \cos \phi) \cdot b
+ \int _{\psi _{x}^{\ast}}^{\frac {\pi}{2}} {\rm d} \psi
\int _{\pi - \phi _{y}^{\ast}}^{\pi - \phi _{x}^{\ast}} {\rm d} \phi
\cdot (- l \sin \psi \cos \phi) \cdot l \sin \psi \sin \phi \nonumber \\
& + \int _{\psi _{x}^{\ast}}^{\frac {\pi}{2}} {\rm d} \psi
\int _{\pi - \phi _{x}^{\ast}}^{\pi} {\rm d} \phi
\cdot a \cdot l \sin \psi \sin \phi \nonumber \\
%sums
& = S^{(1)}(0, \psi _{y}^{\ast})
+ S^{(2)} ( \psi _{y}^{\ast}, \psi _{x}^{\ast} ) 
+ S^{(3)} \left( \psi _{x}^{\ast}, \frac {\pi}{2} \right).
\end{align}

In the above equations, we define and calculate a summation as
\begin{align}
S^{(3)}(\psi _{x}^{\ast}, \psi ^{\ast})
& \equiv \int _{\psi _{x}^{\ast}}^{\psi ^{\ast}} {\rm d} \psi
\int _{0}^{\phi _{x}^{\ast}} {\rm d} \phi
\cdot a \cdot l \sin \psi \sin \phi
+ \int _{\psi _{x}^{\ast}}^{\psi ^{\ast}} {\rm d} \psi
\int _{\phi _{x}^{\ast}}^{\phi _{y}^{\ast}} {\rm d} \phi
\cdot l \sin \psi \cos \phi \cdot l \sin \psi \sin \phi \nonumber \\
& + \int _{\psi _{x}^{\ast}}^{\psi ^{\ast}} {\rm d} \psi
\int _{\phi _{y}^{\ast}}^{\frac {\pi}{2}} {\rm d} \phi
\cdot l \sin \psi \cos \phi \cdot b \nonumber \\
& + \int _{\psi _{x}^{\ast}}^{\psi ^{\ast}} {\rm d} \psi
\int _{\frac {\pi}{2}}^{\pi - \phi _{y}^{\ast}} {\rm d} \phi
\cdot (- l \sin \psi \cos \phi) \cdot b
+ \int _{\psi _{x}^{\ast}}^{\psi ^{\ast}} {\rm d} \psi
\int _{\pi - \phi _{y}^{\ast}}^{\pi - \phi _{x}^{\ast}} {\rm d} \phi
\cdot (- l \sin \psi \cos \phi) \cdot l \sin \psi \sin \phi \nonumber \\
& + \int _{\psi _{x}^{\ast}}^{\psi ^{\ast}} {\rm d} \psi
\int _{\pi - \phi _{x}^{\ast}}^{\pi} {\rm d} \phi
\cdot a \cdot l \sin \psi \sin \phi \nonumber \\
%result
& = 2 l a (\cos \psi _{x}^{\ast} - \cos \psi ^{\ast}) - a^2 (\psi ^{\ast} - \psi _{x}^{\ast})
+ 2 l b (\cos \psi _{x}^{\ast} - \cos \psi ^{\ast}) - b^2 (\psi ^{\ast} - \psi _{x}^{\ast}) \nonumber \\
& - \left[ \left( \frac {l^2}{2} \psi ^{\ast} - \frac {l^2}{2} \sin \psi ^{\ast} \cos \psi ^{\ast} \right)
- \left( \frac {l^2}{2} \psi _{x}^{\ast} - \frac {l^2}{2} \sin \psi _{x}^{\ast} \cos \psi _{x}^{\ast} \right) \right].
\end{align}
Thus we have
\begin{align}
S^{(3)} \left( \psi _{x}^{\ast}, \frac {\pi}{2} \right)
& = 2 l a \cos \psi _{x}^{\ast}- a^2 \left( \frac {\pi}{2} - \psi _{x}^{\ast} \right)
+ 2 l b \cos \psi _{x}^{\ast} - b^2 \left( \frac {\pi}{2} - \psi _{x}^{\ast} \right) \nonumber \\
& - \left[ \frac {\pi}{4} l^2 
- \left( \frac {l^2}{2} \psi _{x}^{\ast} - \frac {l a}{2} \cos \psi _{x}^{\ast} \right) \right], \\
S^{(3)} \left( \psi _{x}^{\ast}, \psi _{xy}^{\ast} \right)
& = 2 l a (\cos \psi _{x}^{\ast} - \cos \psi _{xy}^{\ast}) - a^2 (\psi _{xy}^{\ast} - \psi _{x}^{\ast})
+ 2 l b (\cos \psi _{x}^{\ast} - \cos \psi _{xy}^{\ast}) - b^2 (\psi _{xy}^{\ast} - \psi _{x}^{\ast}) \nonumber \\
& - \left[ \left( \frac {l^2}{2} \psi _{xy}^{\ast} - \frac {l L}{2} \cos \psi _{xy}^{\ast} \right)
- \left( \frac {l^2}{2} \psi _{x}^{\ast} - \frac {l a}{2}  \cos \psi _{x}^{\ast} \right) \right],
\end{align}
in which we define
\begin{equation}
\label{eq:psi-xy-3D}
\psi _{xy}^{\ast} = \arcsin \frac {L}{l}.
\end{equation}
Then,
\begin{equation}
AB^{\rm 3D} (l, a, b | a < l \leqslant L)
= \left[ \frac {l^2}{2} \psi _{x}^{\ast}
+ \frac {3}{2} la \cos \psi _{x}^{\ast}
- a^2 \left( \frac {\pi}{2} - \psi _{x}^{\ast} \right) \right]
+ \left[ \frac {l^2}{2} \psi _{y}^{\ast}
+ \frac {3}{2} lb \cos \psi _{y}^{\ast}
- b^2 \left( \frac {\pi}{2} - \psi _{y}^{\ast} \right) \right]
- \frac {\pi}{4} l^2.
\end{equation}

We have
\begin{align}
S_{\rm coll}^{\rm 3D}
& = b A^{\rm 3D}(l, a | a < l) + a B^{\rm 3D}(l, b | b < l) - AB^{\rm 3D}(l, a, b | a < l \leqslant L) \nonumber \\
%result
& = \left[ 2lb - 2lb F \left( \psi _{x}^{\ast}, \frac {\pi}{2}, \frac {a}{l} \right)
+ 2ab G\left( \psi _{x}^{\ast}, \frac {\pi}{2}, \frac {a}{l} \right)  \right] 
+ \left[ 2la - 2la F \left( \psi _{y}^{\ast}, \frac {\pi}{2}, \frac {b}{l} \right)
+ 2ab G\left( \psi _{y}^{\ast}, \frac {\pi}{2}, \frac {b}{l} \right) \right] \nonumber \\
& - \left[ \frac {l^2}{2} \psi _{x}^{\ast}
+ \frac {3}{2} la \cos \psi _{x}^{\ast}
- a^2 \left( \frac {\pi}{2} - \psi _{x}^{\ast} \right) \right]
- \left[ \frac {l^2}{2} \psi _{y}^{\ast}
+ \frac {3}{2} lb \cos \psi _{y}^{\ast}
- b^2 \left( \frac {\pi}{2} - \psi _{y}^{\ast} \right) \right]
+ \frac {\pi}{4} l^2 \nonumber \\
%reformulated
& = \left[ 2lb - 2lb F \left( \arcsin \frac {a}{l}, \frac {\pi}{2}, \frac {a}{l} \right)
+ 2ab G\left( \arcsin \frac {a}{l}, \frac {\pi}{2}, \frac {a}{l} \right)  \right] \nonumber \\
& + \left[ 2la - 2la F \left(  \arcsin \frac {b}{l}, \frac {\pi}{2}, \frac {b}{l} \right)
+ 2ab G\left(  \arcsin \frac {b}{l}, \frac {\pi}{2}, \frac {b}{l} \right) \right] \nonumber \\
& - \left[ \frac {l^2}{2} \arcsin \frac {a}{l}
+ \frac {3}{2} la \sqrt {1 - \frac {a^2}{l^2}}
- a^2 \left( \frac {\pi}{2} - \arcsin \frac {a}{l} \right) \right] \nonumber \\
& - \left[ \frac {l^2}{2} \arcsin \frac {b}{l}
+ \frac {3}{2} lb \sqrt {1 - \frac {b^2}{l^2}}
- b^2 \left( \frac {\pi}{2} - \arcsin \frac {b}{l} \right) \right]
+ \frac {\pi}{4} l^2.
\end{align}
Then,
\begin{align}
\label{eq:Plab-3D-al1}
P^{\rm 3D}(l, a, b)
& = \left[ \frac {4}{\pi ^2} \frac {l}{a}
- \frac {4}{\pi ^2} \frac {l}{a} F \left( \arcsin \frac {a}{l}, \frac {\pi}{2}, \frac {a}{l} \right)
+ \frac {4}{\pi ^2} G\left( \arcsin \frac {a}{l}, \frac {\pi}{2}, \frac {a}{l} \right)  \right] \nonumber \\
& + \left[  \frac {4}{\pi ^2} \frac {l}{b}
-  \frac {4}{\pi ^2} \frac {l}{b} F \left(  \arcsin \frac {b}{l}, \frac {\pi}{2}, \frac {b}{l} \right)
+ \frac {4}{\pi ^2} G\left(  \arcsin \frac {b}{l}, \frac {\pi}{2}, \frac {b}{l} \right) \right] \nonumber \\
& - \left[ \frac {1}{\pi ^2} \frac {l^2}{ab} \arcsin \frac {a}{l}
+ \frac {3}{\pi ^2} \frac {l}{b} \sqrt {1 - \frac {a^2}{l^2}}
- \frac {2}{\pi ^2} \frac {a}{b} \left( \frac {\pi}{2} - \arcsin \frac {a}{l} \right) \right] \nonumber \\
& - \left[ \frac {1}{\pi ^2} \frac {l^2}{ab}\arcsin \frac {b}{l}
+ \frac {3}{\pi ^2} \frac {l}{a} \sqrt {1 - \frac {b^2}{l^2}}
- \frac {2}{\pi ^2} \frac {b}{a} \left( \frac {\pi}{2} - \arcsin \frac {b}{l} \right) \right]
+ \frac {1}{2 \pi} \frac {l^2}{ab}.
\end{align}
%

% L < l
When $L < l$, we calculate the integral on $\psi$ in $AB(l, a, b)$ in four ranges of $[0, \psi _{y}^{\ast}]$, $(\psi _{y}^{\ast}, \psi _{x}^{\ast}]$, $(\psi _{y}^{\ast}, \psi _{xy}^{\ast}]$, and $(\psi _{xy}^{\ast}, \pi / 2]$ in a separate way.
See Fig.\ref{fig:diagram-3D} (c) for a schematics.
We write $AB^{\rm 3D}(l, a, b)$ when $L < l$ as $AB^{\rm 3D}(l, a, b| L  < l)$, and have
\begin{align}
& AB^{\rm 3D} (l, a, b | L < l) \nonumber \\
%S^{(1)}
& = \int _{0}^{\psi _{y}^{\ast}} {\rm d} \psi
\int _{0}^{\frac {\pi}{2}} {\rm d} \phi
\cdot l \sin \psi \cos \phi \cdot l \sin \psi \sin \phi 
+ \int _{0}^{\psi _{y}^{\ast}} {\rm d} \psi
\int _{\frac {\pi}{2}}^{\pi} {\rm d} \phi
\cdot (- l \sin \psi \cos \phi) \cdot l \sin \psi \sin \phi \nonumber \\
%S^{(2)}
& + \int _{\psi _{y}^{\ast}}^{\psi _{x}^{\ast}} {\rm d} \psi
\int _{0}^{\phi _{y}^{\ast}} {\rm d} \phi
\cdot l \sin \psi \cos \phi \cdot l \sin \psi \sin \phi
+ \int _{\psi _{y}^{\ast}}^{\psi _{x}^{\ast}} {\rm d} \psi
\int _{\phi _{y}^{\ast}}^{\frac {\pi}{2}} {\rm d} \phi
\cdot l \sin \psi \cos \phi \cdot b \nonumber \\
& + \int _{\psi _{y}^{\ast}}^{\psi _{x}^{\ast}} {\rm d} \psi
\int _{\frac {\pi}{2}}^{\pi - \phi _{y}^{\ast}} {\rm d} \phi
\cdot (- l \sin \psi \cos \phi) \cdot b
+ \int _{\psi _{y}^{\ast}}^{\psi _{x}^{\ast}} {\rm d} \psi
\int _{\pi - \phi _{y}^{\ast}}^{\pi} {\rm d} \phi
\cdot (- l \sin \psi \cos \phi) \cdot l \sin \psi \sin \phi  \nonumber \\
%S^{(3)}
& + \int _{\psi _{x}^{\ast}}^{\psi _{xy}^{\ast}} {\rm d} \psi
\int _{0}^{\phi _{x}^{\ast}} {\rm d} \phi
\cdot a \cdot l \sin \psi \sin \phi
+ \int _{\psi _{x}^{\ast}}^{\psi _{xy}^{\ast}} {\rm d} \psi
\int _{\phi _{x}^{\ast}}^{\phi _{y}^{\ast}} {\rm d} \phi
\cdot l \sin \psi \cos \phi \cdot l \sin \psi \sin \phi \nonumber \\
& + \int _{\psi _{x}^{\ast}}^{\psi _{xy}^{\ast}} {\rm d} \psi
\int _{\phi _{y}^{\ast}}^{\frac {\pi}{2}} {\rm d} \phi
\cdot l \sin \psi \cos \phi \cdot b \nonumber \\
& + \int _{\psi _{x}^{\ast}}^{\psi _{xy}^{\ast}} {\rm d} \psi
\int _{\frac {\pi}{2}}^{\pi - \phi _{y}^{\ast}} {\rm d} \phi
\cdot (- l \sin \psi \cos \phi) \cdot b
+ \int _{\psi _{x}^{\ast}}^{\psi _{xy}^{\ast}} {\rm d} \psi
\int _{\pi - \phi _{y}^{\ast}}^{\pi - \phi _{x}^{\ast}} {\rm d} \phi
\cdot (- l \sin \psi \cos \phi) \cdot l \sin \psi \sin \phi \nonumber \\
& + \int _{\psi _{x}^{\ast}}^{\psi _{xy}^{\ast}} {\rm d} \psi
\int _{\pi - \phi _{x}^{\ast}}^{\pi} {\rm d} \phi
\cdot a \cdot l \sin \psi \sin \phi \nonumber \\
%S^{(4)}
& + \int _{\psi _{xy}^{\ast}}^{\frac {\pi}{2}} {\rm d} \psi
\int _{0}^{\phi _{y}^{\ast}} {\rm d} \phi
\cdot a \cdot l \sin \psi \sin \phi
+ \int _{\psi _{xy}^{\ast}}^{\frac {\pi}{2}} {\rm d} \psi
\int _{\phi _{y}^{\ast}}^{\phi _{x}^{\ast}} {\rm d} \phi
\cdot a \cdot b
 + \int _{\psi _{xy}^{\ast}}^{\frac {\pi}{2}} {\rm d} \psi
\int _{\phi _{x}^{\ast}}^{\frac {\pi}{2}} {\rm d} \phi
\cdot l \sin \psi \cos \phi \cdot b \nonumber \\
& + \int _{\psi _{xy}^{\ast}}^{\frac {\pi}{2}} {\rm d} \psi
\int _{\frac {\pi}{2}}^{\pi - \phi _{x}^{\ast}} {\rm d} \phi
\cdot (- l \sin \psi \cos \phi) \cdot b
+ \int _{\psi _{xy}^{\ast}}^{\frac {\pi}{2}} {\rm d} \psi
\int _{\pi - \phi _{x}^{\ast}}^{\pi - \phi _{y}^{\ast}} {\rm d} \phi
\cdot a \cdot b
+ \int _{\psi _{xy}^{\ast}}^{\frac {\pi}{2}} {\rm d} \psi
\int _{\pi - \phi _{y}^{\ast}}^{\pi} {\rm d} \phi
\cdot a \cdot l \sin \psi \sin \phi \nonumber \\
%sum
& = S^{(1)}(0, \psi _{y}^{\ast})
+ S^{(2)} ( \psi _{y}^{\ast}, \psi _{x}^{\ast} ) 
+ S^{(3)} ( \psi _{x}^{\ast}, \psi _{xy}^{\ast} ) 
+ S^{(4)} \left( \psi _{xy}^{\ast}, \frac {\pi}{2} \right).
\end{align}

In the above equations, we define and calculate a summation as
\begin{align}
S^{(4)} \left( \psi _{xy}^{\ast}, \frac {\pi}{2} \right)
& = \int _{\psi _{xy}^{\ast}}^{\frac {\pi}{2}} {\rm d} \psi
\int _{0}^{\phi _{y}^{\ast}} {\rm d} \phi
\cdot a \cdot l \sin \psi \sin \phi
+  \int _{\psi _{xy}^{\ast}}^{\frac {\pi}{2}} {\rm d} \psi
\int _{\phi _{y}^{\ast}}^{\phi _{x}^{\ast}} {\rm d} \phi
\cdot a \cdot b
+  \int _{\psi _{xy}^{\ast}}^{\frac {\pi}{2}} {\rm d} \psi
\int _{\phi _{x}^{\ast}}^{\frac {\pi}{2}} {\rm d} \phi
\cdot l \sin \psi \cos \phi \cdot b \nonumber \\
& +  \int _{\psi _{xy}^{\ast}}^{\frac {\pi}{2}} {\rm d} \psi
\int _{\frac {\pi}{2}}^{\pi - \phi _{x}^{\ast}} {\rm d} \phi
\cdot (- l \sin \psi \cos \phi) \cdot b
+  \int _{\psi _{xy}^{\ast}}^{\frac {\pi}{2}} {\rm d} \psi
\int _{\pi - \phi _{x}^{\ast}}^{\pi - \phi _{y}^{\ast}} {\rm d} \phi
\cdot a \cdot b
+  \int _{\psi _{xy}^{\ast}}^{\frac {\pi}{2}} {\rm d} \psi
\int _{\pi - \phi _{y}^{\ast}}^{\pi} {\rm d} \phi
\cdot a \cdot l \sin \psi \sin \phi \nonumber \\
%simplified form
& = \left[ 2 l b \cos \psi _{xy}^{\ast}
- 2 l b F\left( \psi _{xy}^{\ast}, \frac {\pi}{2}, \frac {a}{l} \right)
+ 2 a b G\left( \psi _{xy}^{\ast}, \frac {\pi}{2}, \frac {a}{l} \right) \right] \nonumber \\
& + \left[ 2 l a \cos \psi _{xy}^{\ast}
- 2 l a F\left( \psi _{xy}^{\ast}, \frac {\pi}{2}, \frac {b}{l} \right)
+ 2 a b G\left( \psi _{xy}^{\ast}, \frac {\pi}{2}, \frac {b}{l} \right) \right]
- a b \pi \left( \frac {\pi}{2} - \psi _{xy}^{\ast} \right).
\end{align}

Thus we have
\begin{align}
AB^{\rm 3D} (l, a, b | L < l)
& = \left[ - 2lb F \left( \psi _{xy}^{\ast}, \frac {\pi}{2}, \frac {a}{l} \right)
+ 2ab G \left( \psi _{xy}^{\ast}, \frac {\pi}{2}, \frac {a}{l} \right) \right] 
+ \left[ - 2la F \left( \psi _{xy}^{\ast}, \frac {\pi}{2}, \frac {b}{l} \right)
+ 2ab G \left( \psi _{xy}^{\ast}, \frac {\pi}{2}, \frac {b}{l} \right) \right] \nonumber \\
& + \left[ \frac {l^2}{2} \psi _{x}^{\ast}
+ \frac {3}{2} la \cos \psi _{x}^{\ast}
- a^2 \left( \frac {\pi}{2} - \psi _{x}^{\ast} \right) \right] 
+ \left[ \frac {l^2}{2} \psi _{y}^{\ast}
+ \frac {3}{2} lb \cos \psi _{y}^{\ast}
- b^2 \left( \frac {\pi}{2} - \psi _{y}^{\ast} \right) \right] \nonumber \\
& - \left[ \frac {l^2}{2} \psi _{xy}^{\ast}
- \frac {l L}{2} \cos \psi _{xy}^{\ast}
- L^2 \left( \frac {\pi}{2} - \psi _{xy}^{\ast} \right) \right]
- ab\pi \left( \frac {\pi}{2} - \psi _{xy}^{\ast} \right).
\end{align}

We have
\begin{align}
\label{eq:Scoll-3D-al2}
S_{\rm coll}^{\rm 3D}
& = b A^{\rm 3D}(l, a | a < l) + a B^{\rm 3D}(l, b | b < l) - AB^{\rm 3D}(l, a, b | L < l) \nonumber \\
%result
& = \left[ 2lb - 2lb F \left( \psi _{x}^{\ast}, \psi _{xy}^{\ast}, \frac {a}{l} \right)
+ 2ab G \left( \psi _{x}^{\ast}, \psi _{xy}^{\ast}, \frac {a}{l} \right) \right] 
+ \left[ 2la - 2la F \left( \psi _{y}^{\ast}, \psi _{xy}^{\ast}, \frac {b}{l} \right)
+ 2ab G \left( \psi _{y}^{\ast}, \psi _{xy}^{\ast}, \frac {b}{l} \right) \right] \nonumber \\
& - \left[ \frac {l^2}{2} \psi _{x}^{\ast}
+ \frac {3}{2} la \cos \psi _{x}^{\ast}
- a^2 \left( \frac {\pi}{2} - \psi _{x}^{\ast} \right) \right]
- \left[ \frac {l^2}{2} \psi _{y}^{\ast}
+ \frac {3}{2} lb \cos \psi _{y}^{\ast}
- b^2 \left( \frac {\pi}{2} - \psi _{y}^{\ast} \right) \right] \nonumber \\
& + \left[ \frac {l^2}{2} \psi _{xy}^{\ast}
- \frac {l L}{2} \cos \psi _{xy}^{\ast}
- L^2 \left( \frac {\pi}{2} - \psi _{xy}^{\ast} \right) \right]
+ ab\pi \left( \frac {\pi}{2} - \psi _{xy}^{\ast} \right) \nonumber \\
%reformulated
& = \left[ 2lb - 2lb F \left( \arcsin \frac {a}{l}, \arcsin \frac {L}{l}, \frac {a}{l} \right)
+ 2ab G \left( \arcsin \frac {a}{l}, \arcsin \frac {L}{l}, \frac {a}{l} \right) \right] \nonumber \\
& + \left[ 2la - 2la F \left( \arcsin \frac {b}{l}, \arcsin \frac {L}{l}, \frac {b}{l} \right)
+ 2ab G \left( \arcsin \frac {b}{l}, \arcsin \frac {L}{l}, \frac {b}{l} \right) \right] \nonumber \\
& - \left[ \frac {l^2}{2} \arcsin \frac {a}{l}
+ \frac {3}{2} la \sqrt {1 - \frac {a^2}{l^2}}
- a^2 \left( \frac {\pi}{2} - \arcsin \frac {a}{l} \right) \right] \nonumber \\
& - \left[ \frac {l^2}{2} \arcsin \frac {b}{l}
+ \frac {3}{2} lb \sqrt {1 - \frac {b^2}{l^2}}
- b^2 \left( \frac {\pi}{2} - \arcsin \frac {b}{l} \right) \right] \nonumber \\
& + \left[ \frac {l^2}{2} \arcsin \frac {L}{l}
- \frac {l L}{2} \sqrt {1 - \frac {L^2}{l^2}}
- L^2 \left( \frac {\pi}{2} - \arcsin \frac {L}{l} \right) \right] 
 + ab\pi \left( \frac {\pi}{2} - \arcsin \frac {L}{l} \right).
\end{align}
Then,
\begin{align}
\label{eq:Plab-3D-al2}
P^{\rm 3D}(l, a, b)
& = \left[ \frac {4}{\pi ^2} \frac {l}{a}
- \frac {4}{\pi ^2} \frac {l}{a} F \left( \arcsin \frac {a}{l}, \arcsin \frac {L}{l}, \frac {a}{l} \right)
+ \frac {4}{\pi ^2} G \left( \arcsin \frac {a}{l}, \arcsin \frac {L}{l}, \frac {a}{l} \right) \right] \nonumber \\
& + \left[ \frac {4}{\pi ^2} \frac {l}{b}
- \frac {4}{\pi ^2} \frac {l}{b} F \left( \arcsin \frac {b}{l}, \arcsin \frac {L}{l}, \frac {b}{l} \right)
+ \frac {4}{\pi ^2} G \left( \arcsin \frac {b}{l}, \arcsin \frac {L}{l}, \frac {b}{l} \right) \right] \nonumber \\
& - \left[ \frac {1}{\pi ^2} \frac {l^2}{ab} \arcsin \frac {a}{l}
+ \frac {3}{\pi ^2} \frac {l}{b} \sqrt {1 - \frac {a^2}{l^2}}
- \frac {2}{\pi ^ 2} \frac {a}{b} \left( \frac {\pi}{2} - \arcsin \frac {a}{l} \right) \right] \nonumber \\
& - \left[ \frac {1}{\pi ^2} \frac {l^2}{ab} \arcsin \frac {b}{l}
+ \frac {3}{\pi ^2} \frac {l}{a} \sqrt {1 - \frac {b^2}{l^2}}
- \frac {2}{\pi ^2} \frac {b}{a} \left( \frac {\pi}{2} - \arcsin \frac {b}{l} \right) \right] \nonumber \\
& + \left[ \frac {1}{\pi ^2} \frac {l^2}{ab} \arcsin \frac {L}{l}
- \frac {1}{\pi ^2} \frac {l L}{ab} \sqrt {1 - \frac {L^2}{l^2}}
- \frac {2}{\pi ^2} \frac {L^2}{ab} \left( \frac {\pi}{2} - \arcsin \frac {L}{l} \right) \right]
+ \frac {2}{\pi} \left( \frac {\pi}{2} - \arcsin \frac {L}{l} \right).
\end{align}
%

%%%%%%%%%%%%%%%%%%%%%%%%%%%%%%%%
\section{Appendix F: Consistency of Equations in Theory of BLNP with $3$D Needles}

%reset equation numbers in Appendix F
\setcounter{equation}{0}
\renewcommand{\theequation}{F.\arabic{equation}}

% l = 0
When $l = 0$, in the case of $l \leqslant b$, we have
\begin{equation}
P^{\rm 3D}(0, a, b) = 0,
\end{equation}
which is intuitively correct.

%l = b
When $l = b$, both the cases of $l \leqslant b$ and $b < l \leqslant a$ lead to
\begin{equation}
S_{\rm coll}^{\rm 3D} = \left( 2 - \frac {\pi}{4} \right) b^2 + 2ab.
\end{equation}
%

%l = a
When $l = a$, both the cases of $b < l \leqslant a$ and $a < l \leqslant L$ correspond to
\begin{align}
S_{\rm coll}^{\rm 3D}
& = 2ab
+ \left[ 2 a^2 - 2 a^2 F \left( \arcsin \frac {b}{a}, \frac {\pi}{2}, \frac {b}{a} \right) 
+ 2 a b G \left( \arcsin \frac {b}{a}, \frac {\pi}{2}, \frac {b}{a} \right) \right] \nonumber \\
& - \left[ \frac {a^2}{2} \arcsin \frac {b}{a}
+ \frac {3}{2} a b \sqrt {1- \frac {b^2}{a^2}}
- b^2 \left( \frac {\pi}{2} - \arcsin \frac {b}{a} \right) \right].
\end{align}
%

%l = L
When $l = L$, both the cases of $a < l \leqslant L$ and $L < l$ result in
\begin{align}
S_{\rm coll}^{\rm 3D}
& = \left[ 2 b L
-  2 b L F \left( \arcsin \frac {a}{L}, \frac {\pi}{2}, \frac {a}{L} \right) 
+ 2 a b G \left( \arcsin \frac {a}{L}, \frac {\pi}{2}, \frac {a}{L} \right) \right] \nonumber \\
& + \left[ 2 a L
-  2 a L F \left( \arcsin \frac {b}{L}, \frac {\pi}{2}, \frac {b}{L} \right) 
+ 2 a b G \left( \arcsin \frac {b}{L}, \frac {\pi}{2}, \frac {b}{L} \right) \right] \nonumber \\
& - \left[ 3ab - a^2 \left( \frac {\pi}{2} - \arcsin \frac {a}{L}\right)
- b^2 \left( \frac {\pi}{2} - \arcsin \frac {b}{L} \right) \right].
\end{align}
%

%%%%%%%%%%%%%%%%%%%%%%%%%%%%%%%%
\section{Appendix G: Case of Infinitely Large $l$ in Theory of BLNP with $3$D Needles}

%reset equation numbers in Appendix G
\setcounter{equation}{0}
\renewcommand{\theequation}{G.\arabic{equation}}

For a general integral $H(a, b, t) \equiv \int _{a}^{b} {\rm d}\psi h(\psi, t)$, if $a, b \to 0$ and both $h(a,t)$ and $h(b, t)$ are well-defined and finite, we have
\begin{equation}
H(a,b,t) \sim \frac {1}{2} (b - a) [h(a,t) + h(b,t)],
\end{equation}
based on the Trapezoidal rule in numerical integration methods \cite{Press.etal-2007-3e}.
Thus we have
\begin{align}
\label{eq:lbFa-3D}
lbF \left( \arcsin \frac {a}{l}, \arcsin \frac {L}{l}, \frac {a}{l} \right)
& \sim lb \cdot \frac {1}{2} \left( \arcsin \frac {L}{l} - \arcsin \frac {a}{l} \right)
\left[ \sqrt {\frac {L^2}{l^2} - \frac {a^2}{l^2}} + \sqrt {\frac {a^2}{l^2} - \frac {a^2}{l^2} }\right] \nonumber \\
& = \frac {1}{2} b^2  \left( \arcsin \frac {L}{l} - \arcsin \frac {a}{l} \right)
\to 0, \\
laF \left( \arcsin \frac {b}{l}, \arcsin \frac {L}{l}, \frac {b}{l} \right)
& \sim \frac {1}{2} a^2  \left( \arcsin \frac {L}{l} - \arcsin \frac {b}{l} \right)
\to 0.
\end{align}
We further have
\begin{align}
G \left( \arcsin \frac {a}{l}, \arcsin \frac {L}{l}, \frac {a}{l} \right)
& \sim \frac {1}{2} \left( \arcsin \frac {L}{l} - \arcsin \frac {a}{l} \right)
\left[ \left( \frac {\pi}{2} - \arcsin \frac {a/l}{L / l} \right)
+ \left( \frac {\pi}{2} - \arcsin \frac {a/l}{a / l} \right) \right] \nonumber \\
& = \frac {1}{2} \left( \arcsin \frac {L}{l} - \arcsin \frac {a}{l} \right)
\left( \frac {\pi}{2} - \arcsin \frac {a}{L} \right)
\to 0, \\
\label{eq:Gba-3D}
G \left( \arcsin \frac {b}{l}, \arcsin \frac {L}{l}, \frac {b}{l} \right)
& \sim \frac {1}{2} \left( \arcsin \frac {L}{l} - \arcsin \frac {b}{l} \right)
\left( \frac {\pi}{2} - \arcsin \frac {b}{L} \right)
\to 0.
\end{align}
Summing up Eqs.(\ref{eq:lbFa-3D})-(\ref{eq:Gba-3D}), Eq.(\ref{eq:Scoll-3D-al2}) reduces to
\begin{align}
S_{\rm coll}^{\rm 3D}
& = \left[ 2lb - 2 \cdot 0 + 2ab \cdot 0 \right]
+ \left[ 2la - 2 \cdot 0 + 2ab \cdot 0 \right] \nonumber \\
& - \left[ \frac {l^2}{2} \cdot \frac {a}{l} + \frac {3}{2} la \cdot \sqrt {1 - 0} - a^2 \cdot \left( \frac {\pi}{2} - 0 \right)\right] 
- \left[ \frac {l^2}{2} \cdot \frac {b}{l} + \frac {3}{2} lb \cdot \sqrt {1 - 0} - b^2 \cdot \left( \frac {\pi}{2} - 0 \right)\right] \nonumber \\
& + \left[ \frac {l^2}{2} \cdot \frac {L}{l} - \frac {l L}{2} \cdot \sqrt {1 - 0}
- L^2 \cdot \left( \frac {\pi}{2} - 0 \right) \right]
+ ab\pi \cdot \left( \frac {\pi}{2} - 0 \right) \nonumber \\
& = ab \frac {\pi ^2}{2}.
\end{align}
Correspondingly, $P^{\rm 3D}(\infty, a, b) = 1$.

%%%%%%%%%%%%%%%%%%%%%%%%%%%%%%%%
\section{Appendix H: Calculation for Theory of BLNP with $3$D Spherocylinders}

%reset equation numbers in Appendix H
\setcounter{equation}{0}
\renewcommand{\theequation}{H.\arabic{equation}}

As Eqs.(\ref{eq:psi-x-3D}), (\ref{eq:psi-y-3D}), and (\ref{eq:psi-xy-3D}) defined for the BLNP with $3$D needles, we define their respective version for the BLNP with $3$D spherocylinders as
\begin{align}
\hat {\psi} _{x}^{\ast}
& = \arcsin \frac {a - \sigma}{l}, \\
\hat {\psi} _{y}^{\ast}
& = \arcsin \frac {b - \sigma}{l}, \\
\hat {\psi} _{xy}^{\ast}
& = \arcsin \frac {\hat{L}}{l}.
\end{align}
%

%l < b
When $l \leqslant b - \sigma$, with Eq. (\ref{eq:Scoll-3DSC-reform}) we have
\begin{align}
S_{\rm coll}^{\rm 3D}
& = (b - \sigma) A^{\rm 3D}(l, a - \sigma | l \leqslant a - \sigma)
+ (a - \sigma) B^{\rm 3D}(l, b - \sigma | l \leqslant b - \sigma) \nonumber \\
& - AB^{\rm 3D}(l, a - \sigma, b - \sigma | l \leqslant b - \sigma)
+ \left( \sigma a + \sigma b - \sigma ^2 \right) \frac {\pi ^2}{2} \nonumber \\
%result
& = 2 l (b - \sigma) + 2 l (a - \sigma) - \frac {\pi}{4} l^2
+ \left( \sigma a + \sigma b - \sigma ^2 \right)  \frac {\pi ^2}{2}.
\end{align}
Then,
\begin{equation}
\label{eq:Plsab-3DSC-lb}
P^{\rm 3D}(l, \sigma, a, b)
= \frac {4}{\pi ^2} \frac {b - \sigma}{b} \frac {l}{a}
+ \frac {4}{\pi ^2} \frac {a - \sigma}{a} \frac {l}{b}
- \frac {1}{2 \pi} \frac {l^2}{ab}
+ \left( \frac {\sigma}{a} + \frac {\sigma}{b} - \frac {\sigma ^2}{ab} \right).
\end{equation}
%

%b < l < a
When $b - \sigma < l \leqslant a - \sigma$, we have
\begin{align}
S_{\rm coll}^{\rm 3D}
%result
& = (b - \sigma) A^{\rm 3D}(l, a - \sigma | l \leqslant a - \sigma)
+ (a - \sigma) B^{\rm 3D}(l, b - \sigma | b - \sigma < l) \nonumber \\
& - AB^{\rm 3D}(l, a - \sigma, b - \sigma | b - \sigma < l \leqslant a - \sigma) 
+ \left( \sigma a + \sigma b - \sigma ^2 \right) \frac {\pi ^2}{2}  \nonumber \\
& = 2 l (b - \sigma)
+ \left[ 2 l (a - \sigma)
- 2 l (a - \sigma) F \left( \hat{\psi} _{y}^{\ast}, \frac {\pi}{2}, \frac {b - \sigma}{l} \right)
+ 2 (a - \sigma) (b - \sigma) G \left( \hat{\psi} _{y}^{\ast}, \frac {\pi}{2}, \frac {b - \sigma}{l} \right) \right] \nonumber \\
& - \left[ \frac {l^2}{2} \hat{\psi}_{y}^{\ast} + \frac {3}{2} l (b - \sigma) \cos \hat{\psi}_{y}^{\ast}
- (b - \sigma)^2 \left( \frac {\pi}{2} - \hat{\psi}_{y}^{\ast} \right) \right]
+ \left( \sigma a + \sigma b - \sigma ^2 \right)  \frac {\pi^2}{2} \nonumber \\
%reformulated
& = 2 l (b - \sigma)
+ \left[ 2 l (a - \sigma)
- 2 l (a - \sigma) F \left( \arcsin \frac {b - \sigma}{l}, \frac {\pi}{2}, \frac {b - \sigma}{l} \right)
+ 2 (a - \sigma) (b - \sigma) G \left( \arcsin \frac {b - \sigma}{l}, \frac {\pi}{2}, \frac {b - \sigma}{l} \right) \right] \nonumber \\
& - \left[ \frac {l^2}{2} \arcsin \frac {b - \sigma}{l}
+ \frac {3}{2} l (b - \sigma) \sqrt {1 - \frac {(b - \sigma)^2}{l^2}}
- (b - \sigma)^2 \left( \frac {\pi}{2} - \arcsin \frac {b - \sigma}{l} \right) \right]
+ \left( \sigma a + \sigma b - \sigma ^2 \right) \frac {\pi^2}{2}.
\end{align}
Then,
\begin{align}
\label{eq:Plsab-3DSC-bla}
& P^{\rm 3D}(l, \sigma, a, b) \nonumber \\
& = \frac {4}{\pi ^2} \frac {b - \sigma}{b} \frac {l}{a}
+ \left[ \frac {4}{\pi ^2} \frac {a - \sigma}{a} \frac {l}{b}
- \frac {4}{\pi ^2} \frac {a - \sigma}{a} \frac {l}{b}F \left( \arcsin \frac {b - \sigma}{l}, \frac {\pi}{2}, \frac {b - \sigma}{l} \right)
+ \frac {4}{\pi ^2} \frac{a - \sigma}{a} \frac{b - \sigma}{b} G \left( \arcsin \frac {b - \sigma}{l}, \frac {\pi}{2}, \frac {b - \sigma}{l} \right) \right] \nonumber \\
& - \left[ \frac {1}{\pi ^2} \frac {l^2}{ab} \arcsin \frac {b - \sigma}{l}
+ \frac {3}{\pi ^2}  \frac {b - \sigma}{b} \frac {l}{a} \sqrt {1 - \frac {(b - \sigma)^2}{l^2}}
- \frac {2}{\pi ^2} \frac {(b - \sigma)^2}{ab} \left( \frac {\pi}{2} - \arcsin \frac {b - \sigma}{l} \right) \right]
+ \left( \frac {\sigma}{a} + \frac{\sigma}{b} - \frac {\sigma ^2}{ab} \right).
\end{align}
%

%a < l < L
When $a - \sigma < l \leqslant \hat {L}$, we have
\begin{align}
S_{\rm coll}^{\rm 3D}
& = (b - \sigma) A^{\rm 3D}(l, a - \sigma | a - \sigma < l)
+ (a - \sigma) B^{\rm 3D}(l, b - \sigma | b - \sigma < l) \nonumber \\
& - AB^{\rm 3D}(l, a - \sigma, b - \sigma | a - \sigma < l \leqslant \hat{L})
+ \left( \sigma a + \sigma b - \sigma ^2 \right) \frac {\pi ^2}{2} \nonumber \\
%result
& = \left[ 2 l (b - \sigma)
- 2 l (b - \sigma) F \left( \hat{\psi} _{x}^{\ast}, \frac {\pi}{2}, \frac {a - \sigma}{l} \right)
+ 2 (a - \sigma) (b - \sigma) G \left( \hat{\psi} _{x}^{\ast}, \frac {\pi}{2}, \frac {a - \sigma}{l} \right) \right]  \nonumber \\
& + \left[ 2 l (a - \sigma)
- 2 l (a - \sigma) F \left( \hat{\psi} _{y}^{\ast}, \frac {\pi}{2}, \frac {b - \sigma}{l} \right)
+ 2 (a - \sigma) (b - \sigma) G \left( \hat{\psi} _{y}^{\ast}, \frac {\pi}{2}, \frac {b - \sigma}{l} \right) \right] \nonumber \\
& - \left[ \frac {l^2}{2} \hat{\psi}_{x}^{\ast} + \frac {3}{2} l (a - \sigma) \cos \hat{\psi}_{x}^{\ast}
- (a - \sigma)^2 \left( \frac {\pi}{2} - \hat{\psi}_{x}^{\ast} \right) \right] \nonumber \\
& - \left[ \frac {l^2}{2} \hat{\psi}_{y}^{\ast} + \frac {3}{2} l (b - \sigma) \cos \hat{\psi}_{y}^{\ast}
- (b - \sigma)^2 \left( \frac {\pi}{2} - \hat{\psi}_{y}^{\ast} \right) \right] \nonumber \\
& + \frac {\pi}{4} l^2 + \left( \sigma a + \sigma b - \sigma ^2 \right) \frac {\pi^2}{2} \nonumber \\
%reformulated
& = \left[ 2 l (b - \sigma)
- 2 l (b - \sigma) F \left( \arcsin \frac {a - \sigma}{l}, \frac {\pi}{2}, \frac {a - \sigma}{l} \right)
+ 2 (a - \sigma) (b - \sigma) G \left( \arcsin \frac {a - \sigma}{l}, \frac {\pi}{2}, \frac {a - \sigma}{l} \right) \right] \nonumber \\
& + \left[ 2 l (a - \sigma)
- 2 l (a - \sigma) F \left( \arcsin \frac {b - \sigma}{l}, \frac {\pi}{2}, \frac {b - \sigma}{l} \right)
+ 2 (a - \sigma) (b - \sigma) G \left( \arcsin \frac {b - \sigma}{l}, \frac {\pi}{2}, \frac {b - \sigma}{l} \right) \right] \nonumber \\
& - \left[ \frac {l^2}{2} \arcsin \frac {a - \sigma}{l}
+ \frac {3}{2} l (a - \sigma) \sqrt {1 - \frac {(a - \sigma)^2}{l^2}}
- (a - \sigma)^2 \left( \frac {\pi}{2} - \arcsin \frac {a - \sigma}{l} \right) \right] \nonumber \\
& - \left[ \frac {l^2}{2} \arcsin \frac {b - \sigma}{l}
+ \frac {3}{2} l (b - \sigma) \sqrt {1 - \frac {(b - \sigma)^2}{l^2}}
- (b - \sigma)^2 \left( \frac {\pi}{2} - \arcsin \frac {b - \sigma}{l} \right) \right] \nonumber \\
& + \frac {\pi}{4} l^2 + \left( \sigma a + \sigma b - \sigma ^2 \right) \frac {\pi^2}{2} .
\end{align}
Then,
\begin{align}
\label{eq:Plsab-3DSC-al1}
P^{\rm 3D}(l, \sigma, a, b)
& = \left[ \frac {4}{\pi ^2} \frac {b - \sigma}{b} \frac {l}{a}
- \frac {4}{\pi ^2} \frac {b - \sigma}{b} \frac {l}{a}F \left( \arcsin \frac {a - \sigma}{l}, \frac {\pi}{2}, \frac {a - \sigma}{l} \right)
+ \frac {4}{\pi ^2} \frac{a - \sigma}{a} \frac{b - \sigma}{b} G \left( \arcsin \frac {a - \sigma}{l}, \frac {\pi}{2}, \frac {a - \sigma}{l} \right) \right] \nonumber \\
& + \left[ \frac {4}{\pi ^2} \frac {a - \sigma}{a} \frac {l}{b}
- \frac {4}{\pi ^2} \frac {a - \sigma}{a} \frac {l}{b}F \left( \arcsin \frac {b - \sigma}{l}, \frac {\pi}{2}, \frac {b - \sigma}{l} \right)
+ \frac {4}{\pi ^2} \frac{a - \sigma}{a} \frac{b - \sigma}{b} G \left( \arcsin \frac {b - \sigma}{l}, \frac {\pi}{2}, \frac {b - \sigma}{l} \right) \right] \nonumber \\
& - \left[ \frac {1}{\pi ^2} \frac {l^2}{ab} \arcsin \frac {a - \sigma}{l}
+ \frac {3}{\pi ^2}  \frac {a - \sigma}{a} \frac {l}{b} \sqrt {1 - \frac {(a - \sigma)^2}{l^2}}
- \frac {2}{\pi ^2} \frac {(a - \sigma)^2}{ab} \left( \frac {\pi}{2} - \arcsin \frac {a - \sigma}{l} \right) \right] \nonumber \\
& - \left[ \frac {1}{\pi ^2} \frac {l^2}{ab} \arcsin \frac {b - \sigma}{l}
+ \frac {3}{\pi ^2}  \frac {b - \sigma}{b} \frac {l}{a} \sqrt {1 - \frac {(b - \sigma)^2}{l^2}}
- \frac {2}{\pi ^2} \frac {(b - \sigma)^2}{ab} \left( \frac {\pi}{2} - \arcsin \frac {b - \sigma}{l} \right) \right] \nonumber \\
& + \frac {1}{2\pi} \frac {l^2}{ab} + \left( \frac {\sigma}{a} + \frac{\sigma}{b} - \frac {\sigma ^2}{ab} \right).
\end{align}
%

%L < l
When $\hat{L} < l$, we have
\begin{align}
\label{eq:Scoll-3DSC-al2}
S_{\rm coll}^{\rm 3D}
& = (b - \sigma) A^{\rm 3D}(l, a - \sigma | a - \sigma < l)
+ (a - \sigma) B^{\rm 3D}(l, b - \sigma | b - \sigma < l) \nonumber \\
& - AB^{\rm 3D}(l, a - \sigma, b - \sigma | \hat{L} < l)
+ \left( \sigma a + \sigma b - \sigma ^2 \right) \frac {\pi ^2}{2} \nonumber \\
%result
& = \left[ 2 l (b - \sigma)
- 2 l (b - \sigma) F \left( \hat{\psi} _{x}^{\ast}, \hat{\psi}_{xy}^{\ast}, \frac {a - \sigma}{l} \right)
+ 2 (a - \sigma) (b - \sigma) G \left( \hat{\psi} _{x}^{\ast}, \hat{\psi}_{xy}^{\ast}, \frac {a - \sigma}{l} \right) \right]  \nonumber \\
& + \left[ 2 l (a - \sigma)
- 2 l (a - \sigma) F \left( \hat{\psi} _{y}^{\ast}, \hat{\psi}_{xy}^{\ast}, \frac {b - \sigma}{l} \right)
+ 2 (a - \sigma) (b - \sigma) G \left( \hat{\psi} _{y}^{\ast}, \hat{\psi}_{xy}^{\ast}, \frac {b - \sigma}{l} \right) \right] \nonumber \\
& - \left[ \frac {l^2}{2} \hat{\psi}_{x}^{\ast} + \frac {3}{2} l (a - \sigma) \cos \hat{\psi}_{x}^{\ast}
- (a - \sigma)^2 \left( \frac {\pi}{2} - \hat{\psi}_{x}^{\ast} \right) \right] \nonumber \\
& - \left[ \frac {l^2}{2} \hat{\psi}_{y}^{\ast} + \frac {3}{2} l (b - \sigma) \cos \hat{\psi}_{y}^{\ast}
- (b - \sigma)^2 \left( \frac {\pi}{2} - \hat{\psi}_{y}^{\ast} \right) \right] \nonumber \\
& + \left[ \frac {l^2}{2} \hat{\psi}_{xy}^{\ast}
- \frac {l \hat {L}}{2} \cos \hat{\psi}_{xy}^{\ast}
- \hat {L}^2 \left( \frac {\pi}{2} - \hat {\psi}_{xy}^{\ast} \right) \right] \nonumber \\
& + (a - \sigma) (b - \sigma) \pi \left( \frac {\pi}{2} - \hat{\psi}_{xy}^{\ast} \right)
+ \left( \sigma a + \sigma b - \sigma ^2 \right) \frac {\pi^2}{2} \nonumber \\
%reformulated
& = \left[ 2 l (b - \sigma)
- 2 l (b - \sigma) F \left( \arcsin \frac {a - \sigma}{l}, \arcsin \frac {\hat{L}}{l}, \frac {a - \sigma}{l} \right)
+ 2 (a - \sigma) (b - \sigma) G \left( \arcsin \frac {a - \sigma}{l}, \arcsin \frac {\hat{L}}{l}, \frac {a - \sigma}{l} \right) \right] \nonumber \\
& + \left[ 2 l (a - \sigma)
- 2 l (a - \sigma) F \left( \arcsin \frac {b - \sigma}{l}, \arcsin \frac {\hat{L}}{l}, \frac {b - \sigma}{l} \right)
+ 2 (a - \sigma) (b - \sigma) G \left( \arcsin \frac {b - \sigma}{l}, \arcsin \frac {\hat{L}}{l}, \frac {b - \sigma}{l} \right) \right] \nonumber \\
& - \left[ \frac {l^2}{2} \arcsin \frac {a - \sigma}{l}
+ \frac {3}{2} l (a - \sigma) \sqrt {1 - \frac {(a - \sigma)^2}{l^2}}
- (a - \sigma)^2 \left( \frac {\pi}{2} - \arcsin \frac {a - \sigma}{l} \right) \right] \nonumber \\
& - \left[ \frac {l^2}{2} \arcsin \frac {b - \sigma}{l}
+ \frac {3}{2} l (b - \sigma) \sqrt {1 - \frac {(b - \sigma)^2}{l^2}}
- (b - \sigma)^2 \left( \frac {\pi}{2} - \arcsin \frac {b - \sigma}{l} \right) \right] \nonumber \\
& + \left[ \frac {l^2}{2} \arcsin \frac {\hat{L}}{l}
- \frac {l \hat{L}}{2} \sqrt {1 - \frac {\hat{L}^2}{l^2}}
- \hat{L}^2 \left( \frac {\pi}{2} - \arcsin \frac {\hat{L}}{l} \right) \right] \nonumber \\
& + (a - \sigma) (b - \sigma) \pi \left( \frac {\pi}{2} - \arcsin \frac {\hat{L}}{l} \right)
+ \left( \sigma a + \sigma b - \sigma ^2 \right) \frac {\pi^2}{2}.
\end{align}
Then,
\begin{align}
\label{eq:Plsab-3DSC-al2}
& P^{\rm 3D}(l, \sigma, a, b) \nonumber \\
& = \left[ \frac {4}{\pi ^2} \frac {b - \sigma}{b} \frac {l}{a}
- \frac {4}{\pi ^2} \frac {b - \sigma}{b} \frac {l}{a} F \left( \arcsin \frac {a - \sigma}{l}, \arcsin \frac {\hat{L}}{l}, \frac {a - \sigma}{l} \right)
+ \frac {4}{\pi ^2} \frac {a - \sigma}{a} \frac {b - \sigma}{b} G \left( \arcsin \frac {a - \sigma}{l}, \arcsin \frac {\hat{L}}{l}, \frac {a - \sigma}{l} \right) \right] \nonumber \\
& + \left[ \frac {4}{\pi ^2} \frac {a - \sigma}{a} \frac {l}{b}
- \frac {4}{\pi ^2} \frac {a - \sigma}{a} \frac {l}{b} F \left( \arcsin \frac {b - \sigma}{l}, \arcsin \frac {\hat{L}}{l}, \frac {b - \sigma}{l} \right)
+ \frac {4}{\pi ^2} \frac {a - \sigma}{a} \frac {b - \sigma}{b} G \left( \arcsin \frac {b - \sigma}{l}, \arcsin \frac {\hat{L}}{l}, \frac {b - \sigma}{l} \right) \right] \nonumber \\
& - \left[ \frac {1}{\pi ^2} \frac {l^2}{ab} \arcsin \frac {a - \sigma}{l}
+ \frac {3}{\pi ^2} \frac {a - \sigma}{a} \frac {l}{b} \sqrt {1 - \frac {(a - \sigma)^2}{l^2}}
- \frac {2}{\pi ^2} \frac {(a - \sigma)^2}{ab} \left( \frac {\pi}{2} - \arcsin \frac {a - \sigma}{l} \right) \right] \nonumber \\
& - \left[ \frac {1}{\pi ^2} \frac {l^2}{ab} \arcsin \frac {b - \sigma}{l}
+ \frac {3}{\pi ^2} \frac {b - \sigma}{b} \frac {l}{a} \sqrt {1 - \frac {(b - \sigma)^2}{l^2}}
- \frac {2}{\pi ^2} \frac {(b - \sigma)^2}{ab} \left( \frac {\pi}{2} - \arcsin \frac {b - \sigma}{l} \right) \right] \nonumber \\
& + \left[ \frac {1}{\pi ^2} \frac {l^2}{ab} \arcsin \frac {\hat{L}}{l}
- \frac {1}{\pi ^2} \frac {l \hat {L}}{ab} \sqrt {1 - \frac {\hat{L}^2}{l^2}}
- \frac {2}{\pi ^2} \frac {\hat{L}^2}{ab} \left( \frac {\pi}{2} - \arcsin \frac {\hat{L}}{l} \right) \right] \nonumber \\
& + \frac {2}{\pi} \frac {a - \sigma}{a} \frac {b - \sigma}{b} \left( \frac {\pi}{2} - \arcsin \frac {\hat{L}}{l} \right)
+ \left( \frac{\sigma}{a} + \frac {\sigma}{b} - \frac {\sigma ^2}{ab} \right).
\end{align}
%

%%%%%%%%%%%%%%%%%%%%%%%%%%%%%%%%
\section{Appendix I: Consistency of Equations in Theory of BLNP with $3$D Spherocylinders}

%reset equation numbers in Appendix I
\setcounter{equation}{0}
\renewcommand{\theequation}{I.\arabic{equation}}

%l = 0
When $l = 0$, based on the case of $l \leqslant b$, we have
\begin{equation}
S_{\rm coll}^{\rm 3D} = \left( \sigma a + \sigma b - \sigma ^2 \right) \frac {\pi ^2}{2},
\end{equation}
which is simply $S_{\rm coll}$  in Eq. (\ref{eq:Scoll-2DSC-l0}) for the problem of $2$D spherocylinders with $l = 0$, multiplying by the range of the new degree of freedom $\psi$, or $\pi / 2$.

%l = b
When $l = b - \sigma$, both the cases of $l \leqslant b - \sigma$ and $b - \sigma < l \leqslant a - \sigma$ reduce to
\begin{equation}
S_{\rm coll}^{\rm 3D}
= \left( 2 - \frac {\pi}{4} \right) (b - \sigma)^2 + 2 (a - \sigma) (b - \sigma)
+ \left( \sigma a + \sigma b - \sigma ^2 \right) \frac {\pi ^2}{2}.
\end{equation}
%

%l = a
When $l = a - \sigma$, both the cases of $b - \sigma < l \leqslant a - \sigma$ and $a - \sigma < l \leqslant \hat {L}$ correspond to
\begin{align}
S_{\rm coll}^{\rm 3D}
& = 2 (a - \sigma) (b - \sigma) \nonumber \\
& + \left[ 2 (a - \sigma)^2
- 2 (a - \sigma)^2 F \left( \arcsin \frac {b - \sigma}{a - \sigma}, \frac {\pi}{2}, \frac {b - \sigma}{a - \sigma} \right)
+ 2 (a - \sigma) (b - \sigma) G \left( \arcsin \frac {b - \sigma}{a - \sigma}, \frac {\pi}{2}, \frac {b - \sigma}{a - \sigma} \right) \right] \nonumber \\
& - \left[ \frac {(a - \sigma)^2}{2} \arcsin \frac {b - \sigma}{a - \sigma}
+ \frac {3}{2} (a - \sigma) (b - \sigma) \sqrt {1 - \frac {(b - \sigma)^2}{(a - \sigma)^2}} 
- (b - \sigma)^2 \left( \frac {\pi}{2} - \arcsin \frac {b - \sigma}{a - \sigma} \right) \right] \nonumber \\
& + \left( \sigma a + \sigma b - \sigma ^2 \right)  \frac {\pi ^2}{2}.
\end{align}
%

%l = L
When $l = \hat{L}$, both the cases of $a - \sigma < l \leqslant \hat{L}$ and $\hat{L} < l$ result in
\begin{align}
S_{\rm coll}^{\rm 3D}
& = \left[ 2 \hat{L} (b - \sigma)
- 2 \hat{L} (b - \sigma) F \left( \arcsin \frac {a - \sigma}{\hat{L}}, \frac {\pi}{2}, \frac {a - \sigma}{\hat{L}} \right) 
+ 2 (a - \sigma) (b - \sigma) G \left( \arcsin \frac {a - \sigma}{\hat{L}}, \frac {\pi}{2}, \frac {a - \sigma}{\hat{L}} \right) \right] \nonumber \\
& + \left[ 2 \hat{L} (a - \sigma)
- 2 \hat{L} (a - \sigma) F \left( \arcsin \frac {b - \sigma}{\hat{L}}, \frac {\pi}{2}, \frac {b - \sigma}{\hat{L}} \right) 
+ 2 (a - \sigma) (b - \sigma) G \left( \arcsin \frac {b - \sigma}{\hat{L}}, \frac {\pi}{2}, \frac {b - \sigma}{\hat{L}} \right) \right] \nonumber \\
& - \left[ 3 (a - \sigma) (b - \sigma)
- (a - \sigma)^2 \left( \frac {\pi}{2} - \arcsin \frac {a - \sigma}{\hat{L}} \right) 
- (b - \sigma)^2 \left( \frac {\pi}{2} - \arcsin \frac {b - \sigma}{\hat{L}} \right) \right] \nonumber \\
& + \left( \sigma a + \sigma b - \sigma ^2 \right) \frac {\pi ^2}{2}.
\end{align}
%

%%%%%%%%%%%%%%%%%%%%%%%%%%%%%%%%
\section{Appendix J: Case of Infinitely Large $l$ in Theory of BLNP with $3$D Spherocylinders}

%reset equation numbers in Appendix J
\setcounter{equation}{0}
\renewcommand{\theequation}{J.\arabic{equation}}

%l -> infty
With Eqs.(\ref{eq:lbFa-3D})-(\ref{eq:Gba-3D}) and substituting $(a, b, L)$ with $(a - \sigma, b - \sigma, \hat{L})$ respectively, we have
\begin{align}
l (b - \sigma) F \left( \arcsin \frac {a - \sigma}{l}, \arcsin \frac {\hat{L}}{l}, \frac {a - \sigma}{l} \right)
& \sim \frac {1}{2} (b - \sigma)^2  \left( \arcsin \frac {\hat{L}}{l} - \arcsin \frac {a - \sigma}{l} \right)
\to 0, \\
l (a - \sigma) F \left( \arcsin \frac {b - \sigma}{l}, \arcsin \frac {\hat{L}}{l}, \frac {b - \sigma}{l} \right)
& \sim \frac {1}{2} (a - \sigma)^2  \left( \arcsin \frac {\hat{L}}{l} - \arcsin \frac {b - \sigma}{l} \right)
\to 0, \\
G \left( \arcsin \frac {a - \sigma}{l}, \arcsin \frac {\hat{L}}{l}, \frac {a - \sigma}{l} \right)
& \sim \frac {1}{2} \left( \arcsin \frac {\hat{L}}{l} - \arcsin \frac {a - \sigma}{l} \right)
\left( \frac {\pi}{2} - \arcsin \frac {a - \sigma}{\hat{L}} \right)
\to 0, \\
G \left( \arcsin \frac {b - \sigma}{l}, \arcsin \frac {\hat{L}}{l}, \frac {b - \sigma}{l} \right) 
& \sim \frac {1}{2} \left( \arcsin \frac {\hat{L}}{l} - \arcsin \frac {b - \sigma}{l} \right)
\left( \frac {\pi}{2} - \arcsin \frac {b - \sigma}{\hat{L}} \right)
\to 0.
\end{align}
Thus Eq.(\ref{eq:Scoll-3DSC-al2}) reduces to
\begin{align}
S_{\rm coll}^{\rm 3D}
& = \left[ 2l (b - \sigma) - 2 \cdot 0 + 2 (a - \sigma) (b - \sigma) \cdot 0 \right] \nonumber \\
& + \left[ 2l (a - \sigma) - 2 \cdot 0 + 2 (a - \sigma) (b - \sigma) \cdot 0 \right] \nonumber \\
& - \left[ \frac {l^2}{2} \cdot \frac {a - \sigma}{l} + \frac {3}{2} l (a - \sigma) \cdot \sqrt {1 - 0} - (a - \sigma)^2 \cdot \left( \frac {\pi}{2} - 0 \right)\right]  \nonumber \\
& - \left[ \frac {l^2}{2} \cdot \frac {b - \sigma}{l} + \frac {3}{2} l (b - \sigma) \cdot \sqrt {1 - 0} - (b - \sigma)^2 \cdot \left( \frac {\pi}{2} - 0 \right)\right] \nonumber \\
& + \left[ \frac {l^2}{2} \cdot \frac {\hat {L}}{l} - \frac {l \hat{L}}{2} \cdot \sqrt {1 - 0}
- \hat{L}^2 \cdot \left( \frac {\pi}{2} - 0 \right) \right] \nonumber \\
& + (a - \sigma) (b - \sigma) \pi \cdot \left( \frac {\pi}{2} - 0 \right)
+ (\sigma a + \sigma b - \sigma^2) \frac {\pi ^2}{2} \nonumber \\
& = ab \frac {\pi ^2}{2}.
\end{align}
Equivalently, $P^{\rm 3D}(\infty, \sigma, a, b) = 1$.

%%%%%%%%%%%%%%%%%%%%%%%%%%%%%%%%

\end{document}